%% file: main.tex
\documentclass[12pt,a4paper,reqno]{amsart}
\usepackage{amscd}
\usepackage{amssymb}
\usepackage{amsthm}
\usepackage[centering,text={14.5cm,22cm}]{geometry}
\usepackage{graphicx}
\usepackage{color}
\usepackage[all]{xy}
\usepackage{mathrsfs}
\usepackage{marvosym}
\usepackage{stmaryrd}
\usepackage{srcltx}
\usepackage{mathtools}
\usepackage{scalerel,amssymb}
\usepackage{tikz-cd}
\usepackage{adjustbox}
\usepackage{placeins}
\usepackage{floatrow}
\newfloatcommand{capbtabbox}{table}[][\FBwidth]
\definecolor{amaranth}{rgb}{0.9, 0.17, 0.31}
\usepackage[colorlinks=true,citecolor=amaranth,linkcolor=black]{hyperref}

\definecolor{shadecolor}{rgb}{1,0.9,0.7}

\newtheorem{theorem}{Theorem}[section]
\newtheorem{lemma}[theorem]{Lemma}

\newtheorem{conjecture}[theorem]{Conjecture}

\newtheoremstyle{defstyle}
  {.6em} 
  {.1em} 
  {} 
  {} 
  {\bfseries} 
  {.} 
  {.5em} 
  {} 
\theoremstyle{defstyle} \newtheorem{definition}[theorem]{Definition}

\newtheorem{example}[theorem]{Example}

\theoremstyle{remark}
\newtheorem{remark}[theorem]{Remark}

\numberwithin{equation}{section}
\numberwithin{figure}{section}

\newcommand {\lfor} {\llbracket}
\newcommand {\rfor} {\rrbracket}

\newcommand{\NN} {\mathbb{N}}

\newcommand{\RR} {\mathbb{R}}

\newcommand {\foD}  {\mathfrak{D}}

\newcommand {\fod}  {\mathfrak{d}}

\newcommand {\gp}  {{\operatorname{gp}}}

\newcommand {\Hom}  {\operatorname{Hom}}

\newcommand {\lra}  {\longrightarrow}

\newcommand {\N}  {\operatorname{N}}

\newcommand {\NE}  {\operatorname{NE}}

\newcommand {\out}  {\mathrm{out}}

\renewcommand{\P}  {\mathscr{P}}

\newcommand {\scrH}  {\mathscr{H}}

\newcommand {\Sing} {\operatorname{Sing}}

\newcommand {\Spec} {\operatorname{Spec}}

\newcommand {\Supp} {\operatorname{Supp}}

\newcommand {\sgn} {\operatorname{sgn}}

\newcommand\Tr{\operatorname{Tr}}

\newcommand {\T} {\mathfrak T}

\newcommand {\Z} {\mathfrak Z}

\newcommand{\hooklongrightarrow}{\lhook\joinrel\longrightarrow}

\def\PP{\mathbb{P}}

\def\Z{\mathbb{Z}}
\def\N{\mathbb{N}}

\def\Q{\mathbb{Q}}

\newtheoremstyle{cited}%
  {3pt}
  {3pt}
  {\itshape}
  {}
  {\bfseries}
  {.}
  {.3em}
  {\thmname{#1} \thmnumber{#2}\thmnote{\normalfont#3}}

\theoremstyle{cited}
\newtheorem{citedthm}{Theorem}

\begin{document}

\title[Fock-Goncharov dual cluster varieties and Gross--Siebert mirrors]{Fock-Goncharov dual cluster varieties and Gross--Siebert mirrors}

\author[H.\,Arg\"uz]{H\"ulya Arg\"uz}
\address{University of Georgia, Department of Mathematics, Athens, GA 30605}
\email{Hulya.Arguz@uga.edu}

\author[P.\,Bousseau]{Pierrick Bousseau}
\address{University of Georgia, Department of Mathematics, Athens, GA 30605}
\email{Pierrick.Bousseau@uga.edu}

\date{}

\begin{abstract}
Cluster varieties come in pairs: for any $\mathcal{X}$ cluster variety there is an associated Fock--Goncharov dual $\mathcal{A}$ cluster variety. On the other hand, in the context of mirror symmetry, associated with any log Calabi--Yau variety is its mirror dual, which can be constructed using the enumerative geometry of rational curves in the framework of the Gross--Siebert program. In this paper we bridge the theory of cluster varieties with the algebro-geometric framework of Gross--Siebert mirror symmetry. Particularly, we show that the mirror to the $\mathcal{X}$ cluster variety is a degeneration of the Fock-Goncharov dual $\mathcal{A}$ cluster variety
and vice versa. To do this, we investigate how the cluster scattering diagram of Gross-Hacking-Keel-Kontsevich compares with the canonical scattering diagram defined by Gross-Siebert to construct mirror duals in arbitrary dimensions. Consequently, we derive an enumerative interpretation of the cluster scattering diagram. Along the way, we prove the Frobenius structure conjecture for a class of log Calabi--Yau varieties obtained as blow-ups of toric varieties. 
\end{abstract}

\maketitle

\tableofcontents
\setcounter{section}{-1}
\section{Introduction}

\subsection{Overview and context}
\label{section_context}
Cluster varieties, introduced by Fock-Goncharov \cite{FG2, FG}, are algebraic varieties constructed by gluing together algebraic tori via distinguished birational maps
referred to as cluster transformations \cite{FZ}. The theory of cluster varieties has deep connections with several areas of mathematics, particularly in algebraic geometry and representation theory: in algebraic geometry they play a significant role in the study of the moduli space of local systems on topological surfaces \cite{FG2, GonSquantum}, and in representation theory they are the geometric counterparts of the cluster algebras introduced by Fomin--Zelevinsky \cite{FZ}. In this paper, we establish new relationships between cluster varieties and mirror symmetry from an algebro-geometric point of view \cite{gross2021canonical}. 
 
A remarkable feature of cluster varieties is that they come in pairs. There are two distinct families of cluster transformations leading to two distinct cluster varieties: the $\mathcal{A}$ cluster variety, called the cluster $K_2$-variety, and the $\mathcal{X}$ cluster variety, called the cluster Poisson variety. Fock and Goncharov conjectured that $\mathcal{A}$ and $\mathcal{X}$ are \emph{dual cluster varieties} in the sense that there exists a canonical basis of 
regular functions on $\mathcal{A}$ parametrized by the integral points of the tropicalization of $\mathcal{X}$, and vice versa\footnote{In this paper we focus attention on skew-symmetric cluster varieties as reviewed in \S\ref{sec_intro_A_X}, and hence the duality simply exchanges $\mathcal{A}$ and $\mathcal{X}$. For the more general situation where one considers skew-symmetrizable cluster varieties see \cite[Appendix A]{GHKK}.}.

Gross-Hacking-Keel \cite{GHKbirational} showed that the original Fock--Goncharov conjecture for cluster varieties cannot hold without additional positivity assumptions. Moreover, without positivity assumptions, they conjectured a
``formal version" of the Fock-Goncharov conjecture, concerning formal families of dual cluster varieties.
In their seminal work Gross-Hacking-Keel-Kontsevich \cite{GHKK} proved the formal Fock-Goncharov conjecture of \cite{GHKbirational}, as well as the original Fock-Goncharov conjecture with the necessary positivity assumptions. Their proof relies on combinatorial methods, and uses the concept of \emph{cluster scattering diagrams} to construct canonical bases for cluster algebras.

On the other hand, the concept of a scattering diagram had emerged in the earlier works of Kontsevich-Soibelman \cite{KS2000}, and Gross-Siebert \cite{GSannals} in a more general context, aiming to understand mirror symmetry from an algebro-geometric point of view. 
Mirror symmetry is a phenomenon emerging from string theory, proposing that Calabi--Yau varieties arise in mirror dual pairs, exhibiting dual properties in their complex and symplectic geometries. 
Given a Calabi--Yau variety\footnote{More precisely, the input of the construction of \cite{gross2019intrinsic, gross2021canonical} is either a maximal log Calabi-Yau variety or a maximal log smooth degeneration of a Calabi-Yau variety.} Gross and Siebert propose that its mirror dual, as a family of complex varieties, can be constructed using logarithmic enumerative geometry \cite{gross2021canonical}.
More precisely, their construction is based on specific curve counts called \emph{punctured log Gromov--Witten invariants}, defined by Abramovich--Chen--Gross--Siebert \cite{ACGS}, which are encoded in a \emph{canonical scattering diagram}.

The cluster duality conjecture of Fock--Goncharov was placed in the framework of mirror symmetry in \cite{GHKbirational, GHKK}, where it is proposed that Fock--Goncharov dual cluster varieties shall arise as mirror duals (see also \cite{GonSmirror}).  
In this paper we verify this expectation, and bridge the theory of cluster varieties with the algebro-geometric framework of mirror symmetry of Gross and Siebert \cite{gross2021canonical}. To do this, we compare the a priori two distinct constructions of scattering diagrams: the combinatorially constructed cluster scattering diagram of Gross--Hacking--Keel--Kontsevich and the canonical scattering diagram of Gross--Siebert defined using the data of punctured log Gromov--Witten invariants. The comparison we obtain allows us to establish a precise relationship between Fock--Goncharov duals and mirror dual varieties in the sense of Gross--Siebert. Along the way, we deduce combinatorial descriptions of punctured log Gromov--Witten invariants of cluster varieties.

After a brief review of mirror symmetry in the framework of the Gross--Siebert program in \S \ref{sec_intro_background}, we state our main results in 
\S \ref{sec_main_results}.

\subsection{Background}
\label{sec_intro_background}
A log Calabi-Yau pair $(X, D)$ is a smooth projective 
variety $X$ over an algebraically closed field $\mathbf{k}$ of characteristic zero together with a reduced simple normal crossing divisor $D$ in $X$ with $K_X+D=0$ \footnote{The definition of a log Calabi-Yau pair in \cite{gross2019intrinsic, gross2021canonical} only requires that $K_X+D=\sum_i a_i D_i$ with $a_i \geq 0$, where $D_i$ are the irreducible components of $D$. The stronger assumption $K_X+D=0$ will always be satisfied for the particular pairs considered in this paper, as in \cite{HDTV}}. In particular, the complement $X \setminus D$
is a Calabi-Yau variety.
Gross and Siebert provide a general construction for the mirror to such a pair $(X,D)$
in any dimension \cite{gross2021canonical}. The mirror produced from this construction is
a family 
\begin{equation} \label{eq_intro_mirror}
\check{\mathfrak{X}} \longrightarrow \mathrm{Spf}\, \mathbf{k} \lfor NE(X) \rfor \, ,\end{equation} 
where $NE(X)$ stands for the monoid spanned by effective curve classes in the group $N_1(X)$ of curve classes in $X$ modulo numerical equivalence. 
The algebra $R(X,D)$ of functions on $\check{\mathfrak{X}}$ is shown to admit a canonical topological basis given by so-called \emph{theta functions}, and hence is referred to as the \emph{algebra of theta functions}. 

The main ingredient to construct the algebra of theta functions is a combinatorial gadget called the \emph{canonical scattering diagram} 
associated to $(X,D)$ and denoted by $\foD_{(X,D)}$. We use the notation $R(\foD_{(X,D)})$, or simply $R(X,D)$, to denote the algebra of theta functions defining the mirror family. We review in  \S\ref{Sec: scattering diagrams and theta functions} how to construct the algebra of theta functions $R(\foD)$ defined by a scattering diagram $\foD$.

As reviewed in \S\ref{Sec: canonical scattering}, the canonical scattering diagram 
encodes the data of counts of rational stable maps to $(X,D)$ with a single marked point where the contact order with $D$ is fixed. The counts of these curves give rise to well defined invariants of $(X,D)$, giving certain  punctured log Gromov--Witten invariants as defined by Abramovich--Chen--Gross--Siebert \cite{ACGS}. While log Gromov--Witten invariants, introduced by Abramovich--Chen \cite{logGWbyAC} and Gross--Siebert \cite{logGW}, are counts of curves in $X$ with prescribed tangencies along $D$, punctured log Gromov--Witten invariants are counts of more general curves which can admit particular types of marked points with negative tangencies along $D$. 

Computing these invariants is technically challenging in general and there are only a handful of cases so far where we know how to describe concretely the canonical scattering diagram and 
deduce from that an explicit description of the mirror to a log Calabi--Yau $(X,D)$ \cite{arguz2021heart,HDTV}: these are situations when $X$ is a projective variety obtained by a blow-up
\begin{equation}
\label{Eq: blow up}
    X \longrightarrow X_{\Sigma}
\end{equation}
of a toric variety $X_{\Sigma}$ associated to a complete fan $\Sigma$ in $\RR^n$, and $D\subset X$ is given by the strict transform of the toric boundary divisor $D_\Sigma$ of $X_\Sigma$. We assume that the center of the blow-up is a union of disjoint connected general smooth hypersurfaces of $D_\Sigma$,
\begin{equation}
    \label{Eq: hypersurfaces H}
    H=\bigcup_{i \in I} H_i \,,
\end{equation}
indexed by a finite set $I$. 
Following \cite{GHK}, we refer to the data of a blow-up map as in \eqref{Eq: blow up}, or equivalently to the data of the pair $(X_{\Sigma},H)$, as a \emph{toric model} for $(X,D)$.

Mark Gross and the first author, in their paper \emph{The Higher Dimensional Tropical Vertex} (HDTV) \cite{HDTV}, showed that for log Calabi--Yau pairs $(X,D)$ obtained from blow-ups of toric varieties as in 
\eqref{Eq: blow up}-\eqref{Eq: hypersurfaces H}, referred to as \emph{HDTV log Calabi--Yau pairs} in the present paper, there is an explicit algorithm computing the canonical scattering diagram for $(X,D)$. 
More precisely, they gave a combinatorial construction of
a scattering diagram $\foD_{(X_{\Sigma},H)}$, referred to as the \emph{HDTV scattering diagram} and reviewed in \S\ref{Sec: HDTV review}, and they proved that the canonical scattering diagram $\foD_{(X,D)}$, which encodes the enumerative geometric information of punctured log Gromov--Witten invariants, can be reconstructed from $\foD_{(X_{\Sigma},H)}$ \cite[Theorem 1.2]{HDTV}. 
Using the HDTV scattering diagram, the first author provided the first explicit equations for mirrors to higher dimensional log Calabi--Yau varieties \cite{arguz2021heart}. We investigate HDTV log Calabi--Yau pairs further in this paper.

\subsection{Main results}
\label{sec_main_results}
We first show in \S\ref{Sec:Extension} that the mirror to a HDTV log Calabi--Yau pair extends to a larger base $S_{(X_\Sigma,H)}$: the mirror family \eqref{eq_intro_mirror}
constructed over $\mathrm{Spf}\, \mathbf{k} \lfor NE(X) \rfor$ as in \cite{gross2021canonical} is recovered as a base change from the extended family over $S_{(X_\Sigma,H)}$. To define the extended base $S_{(X_\Sigma,H)}$,
we first introduce the monoid
\begin{equation}
    \mathcal{M}:= NE(X) \cap (NE(X_\Sigma) \oplus \NN^I) \subset N_1(X) \,,
    \end{equation}
where $NE(X_\Sigma) \oplus \NN^I$ is the submonoid  $N_1(X)$ consisting of all curve classes of the form 
$\overline{\beta}-\sum_{i\in I}a_i E_i$, where 
$\overline{\beta} \in NE(X_\Sigma)$, $(a_i)_{i\in I} \in \NN^I$ and the $E_i$'s are the curve classes of the $\PP^1$-fibers of the exceptional divisors over the $H_i$'s.
Then, we define $S_{(X_\Sigma,H)}$ as the formal completion of 
$\Spec \mathbf{k}[\mathcal{M}]$ along its subscheme $\Spec \mathbf{k}[NE(X_\Sigma)]$. On the other hand, the inclusion of 
$\mathcal{M}$
in $NE(X)$ induces a morphism 
$\mathrm{Spf}\, \mathbf{k}[\![NE(X)]\!] \rightarrow S_{(X_\Sigma,H)}$. We prove in Theorems \ref{thm_intro_extension_1} and \ref{thm_toric_mirror}:
\begin{citedthm}
\label{thm_intro_1}
Let $(X,D)$ be a HDTV log Calabi-Yau pair and $(X_\Sigma,H)$ a toric model of $(X,D)$. Then, there exists a canonical extension 
of the mirror family $\check{\mathfrak{X}}$ to a formal scheme
$\check{\mathfrak{X}}_{\mathrm{ext}}$ over $S_{(X_\Sigma,H)}$ fitting into a fiber diagram
\[\begin{tikzcd}
\check{\mathfrak{X}} 
\arrow[r]
\arrow[ d]
&
\check{\mathfrak{X}}_{\mathrm{ext}}
\arrow[d]\\
\mathrm{Spf}\, \mathbf{k}[\![NE(X)]\!]
\arrow[r]& S_{(X_\Sigma,H)}\,.
\end{tikzcd}\]
Moreover, the restriction of $\check{\mathfrak{X}}_{\mathrm{ext}}$
to $\Spec \mathbf{k}[NE(X_\Sigma)] \subset S_{(X_\Sigma,H)}$ is the mirror family of the toric pair $(X_\Sigma,D_\Sigma)$.
\end{citedthm}

In Theorem \ref{thm_intro_1}, the extension is ``canonical"
in the sense that the product structure on the algebra of functions for $\check{\mathfrak{X}}_{\mathrm{ext}}$ is uniquely determined by the product structure of the theta functions on $\check{\mathfrak{X}}$ -- see  Theorem \ref{thm_intro_extension_1} for details.

In \S\ref{Sec: applications} we investigate two significant applications of Theorem \ref{thm_intro_1}: the first one is the proof of the \emph{Frobenius structure conjecture} for HDTV log Calabi--Yau pairs \cite[Conjecture 0.8, arxiv version 1]{GHK}. This conjecture, roughly, says that there exists a unique algebra structure on the topologically free $\mathbf{k}\lfor NE(X) \rfor$-module generated by a set $\{\vartheta_m\}_m$ of elements indexed by integral points $m$ in the tropicalization of $(X,D)$, such that the $\vartheta_0$-components of products of $\vartheta_m$'s are given in terms of specific log Gromov--Witten invariants of $(X,D)$ with $\psi$ class insertions -- see Conjecture \ref{conj_frobenius} for a precise statement. The existence part of this conjecture follows from the recent work of Johnston \cite[Theorem 1.4]{johnston2022comparison} showing that the algebra of theta functions $R(X,D)$ satisfies the conditions of the conjecture. In Theorem \ref{thm_frob_unique} we prove the uniqueness part for HDTV log Calabi-Yau pairs, by showing that the product structure on the algebra of theta functions is uniquely determined by the trace of products of two or three theta functions. Consequently, in Theorem \ref{thm_frobenius} we obtain:
\begin{citedthm}
The Frobenius structure conjecture holds for HDTV log Calabi--Yau pairs.
\end{citedthm}

As a second application of Theorem \ref{thm_intro_1}, in \S\ref{Sec: the HDTV mirror family} we define the HDTV mirror family to a HDTV log Calabi--Yau pair $(X,D)$ as the base change of the extended family $\check{\mathfrak{X}}_{\mathrm{ext}} \to S_{(X_\Sigma,H)}$ along the morphism 
$\mathrm{Spf}\, \mathbf{k}[\![\NN^I]\!] \rightarrow S_{(X_\Sigma,H)}$
obtained by restriction from $\Spec \mathbf{k}[NE(X_\Sigma)]$ to the unit $1$ of its big torus orbit 
$\Spec \mathbf{k}[N_1(X_\Sigma)]$ 
-- see Definition \ref{def_intro_rest_mirror}. Theorem \ref{thm_restriction} then proves the following.
\begin{citedthm}
Let $(X,D)$ be a HDTV log Calabi-Yau pair and $(X_\Sigma,H)$ a toric model of $(X,D)$. 
Then, the algebra of functions on the HDTV mirror to $(X,D)$ is isomorphic to the algebra of theta functions $R(\foD_{(X_{\Sigma},H)})$ defined by the HDTV scattering diagram $\foD_{(X_{\Sigma},H)}$.
\end{citedthm}

In \S\ref{sec_C_scattering} after a brief review of cluster varieties, and their interpretation as blow-ups of toric varieties following \cite{GHKbirational}, we show how to compute the HDTV mirror to log Calabi-Yau compactifications $(X,D)$ of $\mathcal{X}$ and $\mathcal{A}$ cluster varieties, which are examples of HDTV log Calabi-Yau pairs. To do this we provide a concrete description of the HDTV scattering diagram $\foD_{(X_{\Sigma},H)}$ associated to such cluster varieties in \S\ref{Sec: HDTV for cluster}. We then provide a comparison of $\foD_{(X_{\Sigma},H)}$ with the cluster scattering diagram associated to 
$\mathcal{A}_{\mathrm{prin}}$,
the $\mathcal{A}$ cluster variety with principal coefficients.
Following the notation of \cite{GHKK} 
we denote by $\foD^{\mathcal{A}_{\mathrm{prin}}}_{\mathbf{s}}$ this cluster scattering diagram.

The cluster scattering diagram constructed in \cite{GHKK} and reviewed in \S\ref{Sec: C scattering diagrams}
does not satisfy all the requirements of a scattering diagram in the sense of
\S \ref{Sec: scattering diagrams and theta functions}.
Indeed, a scattering diagram $\foD$ as in \S \ref{Sec: scattering diagrams and theta functions}
is defined over a monoid, and asked to satisfy the finiteness condition in Definition \ref{def_scattering_general}, requiring that $\foD$ can be described order-by-order by finitely many walls. This enables one to build the algebra of theta functions $R(\foD)$ from $\foD$ order by order. However, for the cluster scattering diagram there is no data of a monoid, and one imposes different technical assumptions to be able to construct the algebra of cluster theta functions -- see \S\ref{Sec: C scattering diagram} for details. 

To clarify the distinction between scattering diagrams arising in the context of mirror symmetry \cite{gross2021canonical} and scattering diagrams arising in the context of cluster varieties \cite{GHKK}, we use
following \cite{KY} the terminology ``C-scattering diagram"  for a ``scattering diagram" in the sense of \cite{GHKK}. We review $C$-scattering diagrams, and show that cluster scattering diagrams are particular types of $C$-scattering diagrams in \S\ref{Sec: C scattering diagrams}. We then investigate the cluster scattering diagram for the $\mathcal{A}_{\mathrm{prin}}$ cluster variety, denoted by $\foD^{\mathcal{A}_{\mathrm{prin}}}_{\mathbf{s}}$, as a specific $C_{\mathrm{prin}}$-scattering diagram -- see \S\ref{Sec Cprin and Aprin} for the definition of a $ C_{\mathrm{prin}}$ scattering diagram. Theorem \ref{Thm: HDTV and cluster} enables us to compare scattering diagrams appearing in the framework of mirror symmetry, as in \S\ref{Sec: scattering diagrams}, with cluster scattering diagrams, and shows:
\begin{citedthm}
\label{Thm_intro_comparision}
There is a well-defined map
\[ \Psi \colon C_{\mathrm{prin}}-\mathrm{Scatt} \longrightarrow \mathrm{Scatt} \, , \]
from the set of $ C_{\mathrm{prin}}$-scattering diagrams to the set of scattering diagrams, which maps the $\mathcal{A}_{\mathrm{prin}}$ cluster scattering diagram to the HDTV scattering diagram for the $\mathcal{X}$ cluster variety:
\[ \Psi(\foD^{\mathcal{A}_{\mathrm{prin}}}_{\mathbf{s}})
=\foD_{(X_\Sigma,H)}\,.\]
\end{citedthm}

Theorem \ref{Thm_intro_comparision} provides us with one of the main ingredients to prove an isomorphism between the mirror to the $\mathcal{X}$ (resp.\ $\mathcal{A}$) cluster variety 
and a degeneration of the Fock--Goncharov dual $\mathcal{A}$ (resp.\ 
$\mathcal{X}$) cluster variety:

\begin{citedthm}
\label{Thm_intro_main}
Let $(X,D)$ be a log Calabi-Yau compactification of the $\mathcal{X}$
(resp.\ $\mathcal{A}$) cluster variety.
Then,
the HDTV mirror family $ \check{\mathfrak{X}}_{\mathrm{HDTV}}$ of $(X,D)$ 
is isomorphic to a degeneration of the
$\mathcal{A}$ (resp.\ $\mathcal{X}$) cluster variety.
\end{citedthm}

In Theorem \ref{Thm_intro_main}, the degenerations of the cluster varieties are defined as formal completions of cluster varieties with 
principal coefficients $\mathcal{A}_{\mathrm{prin}}$ and 
$\mathcal{X}_{\mathrm{prin}}$, studied in \cite{GHKK} in the $\mathcal{A}$ case and defined by \cite{Xcoeff} in the $\mathcal{X}$ case-- see Theorems \ref{thm_main} and \ref{thm_main_A} for details.
We also prove a version of Theorem \ref{Thm_intro_main} for the symplectic leaves of the $\mathcal{X}$ cluster variety and quotients of the 
$\mathcal{A}$ cluster variety-- see Theorems \ref{thm_mirror_X_fiber} and 
\ref{thm_mirror_A_quotient}.
Theorem \ref{Thm_intro_main} shows how the mirror symmetry heuristic motivating the work Gross-Hacking-Keel-Kontsevich \cite{GHKK} (see also \cite{HKicm} for a more expository presentation) is concretely realized by the general Gross-Siebert mirror construction. As an application, we obtain in Theorem \ref{Thm: structure constants in cluster} an enumerative interpretation for the structure constants in the algebra of cluster theta functions defined by the cluster scattering diagram:

\begin{citedthm}
The structure constants of the algebra of $\mathcal{A}_{\mathrm{prin}}$
(resp.\ $\mathcal{X}_{\mathrm{prin}}$)
cluster theta function are expressed in terms of punctured log Gromov--Witten invariants of log Calabi-Yau compactifications of the $\mathcal{X}$ (resp.\ $\mathcal{A}$) cluster variety.
\end{citedthm}

\subsection{Related work}
\label{Sec: related works}
The Frobenius structure conjecture was proven previously in two cases: for cluster varieties \cite{mandel2019theta} by Mandel -- see \cite[Theorem 1.5]{mandel2019theta} and \cite[Theorem 2.16]{mandelscattering}, and for log Calabi-Yau pairs $(X,D)$ with $X \setminus D$ affine and containing a torus by Keel--Yu \cite{KY} in the context of their non-archimedean mirror construction, which is shown to be equivalent to the Gross--Siebert mirror construction by Johnston \cite[Corollary 1.2]{johnston2022comparison}.
In the cluster case, under the assumption that $X \setminus D$ is affine, 
a comparison between the non-archimedean construction of \cite{KY} and the work of Gross--Hacking--Keel--Kontsevich \cite{GHKK} is given in \cite[Theorem 1.19]{KY}.
Note that for HDTV log Calabi-Yau pairs, while the complement 
$X\setminus D$ always contains a torus, it is not affine in general. 
The relationship between the HDTV scattering diagram of \cite{HDTV} and the cluster scattering diagram of \cite{GHKK} is also discussed in the work of Mou \cite[Lemma 6.18]{mou2021scattering} in the context of generalized cluster algebras, and by Bardwell-Evans--Cheung--Hong--Lin  \cite[\S 6.2]{bardwell2021scattering} in the special case of rank two cluster varieties.
Different manifestations of mirror symmetry from the symplectic point of view, such as homological mirror symmetry, were also investigated in the context of cluster varieties by Gammage and Le \cite{gammage2021mirror}.

\subsection{Acknowledgments} We thank Mark Gross and Tom Coates for many useful discussions related to the extensions of Gross--Siebert mirror families.
The research of Hülya Argüz was partially supported by the NSF grant DMS-2302116. The research of Pierrick Bousseau was partially supported by the NSF grant DMS-2302117.

\vspace{0.5cm}

\emph{Conventions:} Throughout the paper $\mathbf{k}$ denotes an algebraically closed field of characteristic zero. All monoids we consider in this paper are commutative. Given a monoid $Q$, we denote by $Q^{\mathrm{gp}}$ the associated Grothendieck group and $Q^{\mathrm{gp}}_{\RR} = Q^{\mathrm{gp}} \otimes_{\Z} \RR$.

\section{Scattering diagrams in mirror symmetry}
\label{Sec: scattering diagrams}

\subsection{Scattering diagrams and theta functions}
\label{Sec: scattering diagrams and theta functions}
The theory of scattering diagrams and theta functions is presented in great generality in \cite{GHS}. In this section, we briefly review the 
key notions that we will use, making simplifying assumptions that will always hold in this paper.
\subsubsection{Scattering diagrams}
\label{Sec:scatt}
\begin{definition}
\label{def_pseudomanifold}
A \emph{conical affine pseudomanifold} $(B, \P)$ of dimension $n$ consists of:
\begin{itemize}
    \item[(i)] a set $\P$ of at most $n$-dimensional rational polyhedral cones such that $\tau \in \P$ if $\sigma \in \P$ and $\tau$ is a face of $\sigma$. 
    \item[(ii)] A topological manifold $B$ obtained by gluing together
     without self-intersection the cones $\sigma \in \P$ along their faces by integral linear maps: 
    $B=\cup_{\sigma \in \P} \sigma$, such that in $B$, every cone is contained in an $n$-dimensional cone, every $(n-1)$-dimensional cone is the common face of two $n$-dimensional cones, and every point $x$ in the interior of a at most $(n-2)$-dimensional cone $\tau$ admits a basis of open neighborhoods $V$ such that $V \setminus (V \cap \tau)$ is connected.
    \item[(iii)] An integral affine structure on the complement $B_0=B\setminus \Delta$ of the union $\Delta$ of at most $(n-2)$-dimensional cones of $\P$, which restricts to the standard integral affine structure on each $n$-dimensional cone of $\P$.
\end{itemize}

We call an $n$-dimensional cone of $\P$ a \emph{maximal cone} of $\P$, and a
$(n-1)$-dimensional cone of $\P$ a \emph{codimension one cone} of $\P$.
\end{definition}

In the examples considered in this paper, $B$ will be always diffeomorphic to a real vector space, and $\P$ will be either isomorphic or piecewise-linearly isomorphic to a fan in a real vector space. 
From now on, we fix a conical affine pseudomanifold $(B,\P)$ of dimension $n$. Let $Q$ be a commutative monoid such that the associated abelian group $Q^{\mathrm{gp}}$ is free of finite rank, the natural map $Q \rightarrow Q^{\mathrm{gp}}$ is injective, and $Q$ is contained in a stricly convex cone in the real vector space $Q^{\mathrm{gp}}\otimes \RR$. In particular, there exists a linear form $\ell : Q^{\mathrm{gp}}\otimes \RR \rightarrow \RR$ which is positive on $Q \setminus \{0\}$. Then, for every $k \in \Z_{\geq 0}$, $I_k :=\{ \beta \in Q\,|\, \ell(\beta) \geq k\}$ is a monoid ideal such that $Q \setminus I_k$ is finite. Note that $I_1=Q\setminus \{0\}$ is the maximal monoid ideal of $Q$. 
For every ring $R$, we denote by $ R[\![Q]\!]$ the $R$-algebra
\[ R[\![Q]\!]=\varprojlim_k R[Q]/I_k\,.\]
Concretely elements of $R[\![Q]\!]$ are formal power series $\sum_{\beta \in Q} c_\beta t^\beta$ with $c_\beta \in R$. If one uses a different linear form $\ell'$ with the same properties as $\ell$, defining ideals $I_k'$, then for every $k$, we have $I_m' \subset I_k$ and $I_m \subset I_k'$ for $m$ large enough, and so 
the algebra $R[\![Q]\!]$ is actually independent of the choice of $\ell$.\footnote{In \cite{GHS, gross2019intrinsic, gross2021canonical}, the monoid $Q$ is assumed to be finitely generated, but it is in fact enough to have an exhausting increasing sequence of monoid ideals $I_k$ such that $Q\setminus I_k$ is finite for every $k$. In this paper, we will take $Q$ equal to monoids of effective curves, which are not finitely generated in general. One could instead as in \cite{gross2019intrinsic, gross2021canonical} choose a finitely generated monoid containing the monoid of effective curves, but this would make our results less transparent.}

For every maximal cone $\sigma \in \P$, we denote by $\Lambda_\sigma$ the rank $n$ lattice of integral tangent vectors to $B_0$ at a given point $x \in B_0 \cap \sigma$. 
For every $(n-1)$-dimensional rational polyhedral cone $\fod \subset B$ 
contained in an $n$ or $(n-1)$-dimensional cone of $\P$, we denote by $\Lambda_\fod$ 
the rank $(n-1)$ lattice of integral tangent vectors to $\fod$ at a given point $x \in B_0 \cap \fod$.
For every $m \in \Lambda_\sigma$ or $\Lambda_\fod$, we denote by $z^m$ the corresponding monomial in the monoid algebra $\mathbf{k}[\Lambda_\sigma]$ or $\mathbf{k}[\Lambda_\fod]$, and by $\mathbf{k}[z^m]$ the algebra of polynomials in the variable $z^m$.

We also fix a multivalued piecewise-linear (MVPL) function $\varphi$ on $(B,\P)$ as in \cite[Def.\ 1.8]{GHS}. By definition, 
$\varphi$ is the data, for every codimension one cone $\rho$ of $\P$, of a 
$Q_\RR^{\mathrm{gp}}$-valued piecewise-linear (PL) function $\varphi_\rho$ on $\sigma^+ \cup \sigma^-$, where $\sigma^+$ and $\sigma^-$ are the maximal cones of $\P$ having $\rho$ as a common face. 
For every $\rho$, the function $\varphi_\rho$ is determined, up to a linear function, by specifying its kink $\kappa_{\rho}
\in Q^{\gp}$, defined as follows \cite[Def.\ 1.6, Prop.\ 1.9]{GHS}. 
\begin{definition}
\label{def:kinks}
Let $\rho$ be a codimension one cone of $\P$ and let $\sigma^+,\sigma^-$ be the two maximal cones of $\P$
having  $\rho$ as a commom face, and let $\varphi_{\rho}$ be a piecewise linear
function on $\sigma^+\cup \sigma^-$.
For every point $x$ in $\rho \setminus \Delta$, let $\Lambda_x$ be the lattice of integral tangent vectors to $B_0$ at $x$, which is canonically identified by the integral affine structure on $\sigma^+ \cup \sigma^-$ to both $\Lambda_{\sigma^+}$ and
$\Lambda_{\sigma^-}$. Let $\delta:\Lambda_x\to\Z$ be the quotient by
$\Lambda_\rho \subset \Lambda_x$. We fix signs by requiring that $\delta$
is non-negative on tangent vectors pointing from $\rho$ into
$\sigma^{-}$. Let $n^+,n^-\in \Hom(\Lambda_x, \Z) \otimes Q^\gp$ be the slopes of
$\varphi_{\rho}|_{\sigma^+}$, $\varphi_{\rho}|_{\sigma^-}$, respectively. Then
$(n^--n^+)(\Lambda_\rho)=0$ and hence there exists $\kappa_{\rho}\in Q^\gp$
with
\begin{equation}
\nonumber
n^- -n^+ =\delta \cdot\kappa_{\rho}.
\end{equation}
We refer to $\kappa_{\rho}$ as the \emph{kink} of $\varphi_{\rho}$ along $\rho$.
\end{definition}
In what follows, we assume that $\varphi$ is \emph{$Q$-convex}, that is, $\kappa_\rho
\in Q\setminus \{0\}$ for all $\rho$.

\begin{definition}
\label{def_wall}
A \emph{wall} in $(B,\P)$ over $(Q,\varphi)$ is a pair
$(\fod,f_{\fod})$, consisting of an $(n-1)$-dimensional rational polyhedral cone $\fod \subset B$ contained in a maximal cone 
of $\P$, together with an attached function 
$f_{\fod} \in 
\mathbf{k}[z^{-m_0}][\![Q]\!] \subset \mathbf{k}[\Lambda_\fod][\![Q]\!]$
for some nonzero primitive $m_0 \in \Lambda_\fod$, and
such that $f_\fod=1 \mod I_1$.
Explicitly, $f_\fod$ is of the form
\[ f_\fod =\sum_{\substack{\beta \in Q \\ k \in \NN }}c_{\beta,k} t^{\beta} z^{-k m_0} \]
with $c_{\beta,m} \in \mathbf{k}$, $c_{0,0} =1$, and $c_{0,k}=0$
for all $k\geq 1$.
We say that a wall $(\fod,f_{\fod})$
is \emph{incoming} if 
$\fod=\fod-\RR_{\geq 0}m_0$, and \emph{outgoing} elsewise. We call $m_0$ the \emph{direction} of the wall.
\end{definition}

\begin{definition}
\label{def_scattering_general}
A \emph{scattering diagram} $\foD$ in $(B,\P)$ over $(Q,\varphi)$ is a set of walls $(\fod,f_{\fod})$, such that the \emph{finiteness condition} holds, that is: for every $k \geq 1$, there are only finitely many walls $(\fod,f_{\fod}) \in \foD$ with $f_{\fod} \neq 1 \mod I_k$.
\end{definition}

The \emph{support} of a scattering diagram $\foD$, denoted by $\mathrm{Supp}(\foD)$, is the union of all cones $\fod$ supporting a wall $(\fod,f_{\fod})$ of $\foD$. The \emph{singular locus} of $\foD$ is given by
\begin{equation}
\nonumber
\Sing(\foD)  := \Delta \cup \bigcup_{\fod\in\foD} \partial\fod
\cup \bigcup_{\fod,\fod'\in\foD} (\fod\cap\fod'),
\end{equation}
where the last union is over all pairs of walls $\fod,\fod'$ with
$\fod\cap\fod'$ codimension at least two, and $\Delta$ is as in Definition \ref{def_pseudomanifold}.

If $x\in B\setminus \Sing(\foD)$, we define
\begin{equation}
\label{eq: fx def}
f_x:=\prod_{x\in\fod} f_{\fod} \, , 
\end{equation}
where the product is over all the walls $(\fod,f_{\fod})$ of $\foD$ containing $x$. We say that two scattering diagram $\foD$, $\foD'$ are
\emph{equivalent} if $f_x=f'_x$ for all 
\[ x\in B\setminus (\Sing(\foD)
\cup\Sing(\foD'))\,.\]

\subsubsection{Broken lines}
In this section we overview how to write broken lines. Broken lines are used in the following \S \ref{Sec:theta_functions}
to define the algebra of theta functions.

\begin{definition}
\label{Def: broken line}
Let $\foD$ be a scattering diagram in $(B,\P)$ over $(Q,\varphi)$. 
A \emph{broken line} for $\foD$ is the data of 
\begin{itemize}
\item[(i)] a proper continuous map
\begin{equation}
    \nonumber
    \beta \colon (-\infty,0] \lra B_0 \setminus (B_0 \cap \Sing(\foD))
\end{equation}
along with a sequence $-\infty = t_0 <
t_1 < \cdots< t_r = 0$ such that $r\geq 1$ and,
for all $1 \leq i \leq r$, $\beta((t_{i-1},t_{i}))$ is contained in a maximal cone $\sigma_i \in \P$, 
\item[(ii)]for all $1 \leq i \leq r$,
monomials 
\[ a_i z^{m_i} \in \mathbf{k}[\Lambda_{\sigma_i}][\![Q]\!] \] 
with $a_i \in \mathbf{k}[\![Q]\!]$ and $m_i \in \Lambda_{\sigma_i}$, 
\end{itemize}

subject to the following conditions:
\begin{itemize}
    \item[(i)] for every $1 \leq i \leq r$, 
    $\beta|_{(t_{i-1},t_i)}$ is a non-constant affine map with \[\beta'(t)=-m_i\] for all $t \in (t_{i-1},t_i)$. Moreover, the intersections of the image of $\beta$ with $\Supp(\foD)$ and
    with codimension one cones of $\P$ are transverse. 
    Finally, we have $\beta(0) \notin \Supp(\foD)$, and for every $1 \leq i \leq r-1$, we have either $\beta(t_i) \in \Supp(\foD)$, or 
    $\beta(t_i)$ is contained in a codimension one cone of $\P$.
    \item[(ii)] $a_1=1$, and for every $1 \leq i \leq r-1$, we have the following relations between the monomials $a_i z^{m_i}$ and $a_{i+1} z^{m_{i+1}}$. If $\beta(t_i)$ is not on a codimension one cone of $\P$, then $\sigma_i=\sigma_{i+1}$, 
    $\beta(t_i) \in \Supp(\foD)$, and one requires $a_{i+1}z^{m_{i+1}}$
    to be a monomial distinct from $a_i z^{m_i}$ and contained in the expansion of 
    \begin{equation} \label{eq_bending_1}
    a_i z^{m_i} f_{\beta(t_i)}^{\langle n,m_i \rangle} \,,\end{equation}
    where $n\in \Hom(\Lambda_{\sigma_i},\Z)$ is the unique primitive normal vector to $\Supp(\foD)$ such that $\langle n, m_i \rangle >0$, and $f_{\beta(t_i)}$ is as in \eqref{eq: fx def} with $x=\beta(t_i)$.
    
    If $\beta(t_i)$ is on a codimension one cone $\rho$ of $\P$, then one requires $a_{i+1}z^{m_{i+1}}$
    to be a monomial contained in the expansion of 
    \begin{equation} \label{eq_bending_2}
    a_i z^{m_i}(t^{\kappa_\rho} f_{\beta(t_i)})^{\langle n,m_i \rangle} \,,\end{equation}
    where $n\in \Hom(\Lambda_{\sigma_i},\Z)$ is the unique primitive normal vector to $\rho$ such that $\langle n, m_i \rangle >0$, 
    and $f_{\beta(t_i)}$ is as in \eqref{eq: fx def} with $x=\beta(t_i)$, 
    $\kappa_\rho \in Q$ is the kink of $\varphi$ across $\rho$, and we identify $\Lambda_{\sigma_{i+1}}$ with $\Lambda_{\sigma_i}$ by parallel transport across $\rho$.
\end{itemize}

We call $m_1$ the \emph{asymptotic monomial of $\beta$}. We denote by $a_\beta z^{m_\beta}:= a_r z^{m_r}$ the \emph{final monomial} carried by $\beta$, and we call $\beta(0)$ the \emph{endpoint} of $\beta$.
\end{definition}

\begin{remark}
\label{Rem: crossing walls}
Definition \ref{Def: broken line}, roughly put, says that a broken line $\beta$ starts its life coming from infinity as a line with asymptotic direction $m_1$ decorated with the  monomial $z^{m_1}$. Each time $\beta$ crosses a codimension one cone of $\P$ or a wall $\fod$, 
it can either go straight without bending, or can bend in the direction of the wall, 
until it finally reaches an endpoint point $\beta(0)$ in $B$ and stops (see Figure \ref{Fig:figcan} for an illustration of broken lines). 
\end{remark}

\begin{definition}
\label{def_asymptotic_direction}
Let $B(\Z)$ be the set of integral points of $B$, defined as the union of the sets of integral points of the cones $\sigma \in \P$. For every $m \in B(\Z)$, we say that a broken line $\beta$ has \emph{asymptotic direction} $m$ if the cone $\sigma_1 \in \P$ containing the asymptotic part of $\beta$ contains $m$, and if $m=m_1$ after identification of the integral points of $\sigma_1$ with tangent directions. 
\end{definition}

\subsubsection{Theta functions} \label{Sec:theta_functions}
Now we are ready to define theta functions from broken lines following \cite[\S~3.3]{GHS}. 
Broken lines can be used to define an algebra of theta functions under the assumption that the scattering diagram $\foD$ is \emph{consistent}.
As we will not need the details of the definition in general, we refer to \cite[Definition 3.9]{GHS} for the notion of a consistent scattering diagram. We will describe and use a special case of this notion later in \S \ref{Sec:scattering_M_R}. 
Given a consistent scattering diagram $\foD$ in 
$(B,\P)$ over $(Q,\varphi)$, \cite[Theorem 3.19]{GHS}
produces a $\mathbf{k}[\![Q]\!]$-algebra structure on the $\mathbf{k}[\![Q]\!]$-module
\[ R(\foD)
:= \varprojlim_k
\bigoplus_{m \in B(\Z)} \big( \mathbf{k}[Q]/I_k \big) \vartheta_m \]
where $\vartheta_m$ are basis elements indexed by the integral points 
$m \in B(\Z)$. We refer to the basis elements 
$\vartheta_m$ as \emph{theta functions} and to $R(\foD)$ with this algebra structure 
as the \emph{algebra of theta functions} defined by $\foD$.

The $\mathbf{k}[\![Q]\!]$-algebra structure on $R(\foD)$ is determined by a set of structure constants 
\[C_{m_1m_2}^m \in \mathbf{k}[\![Q]\!]\,,\]
indexed by
$m_1,m_2,m\in B(\Z)$, 
such that for all $m_1,m_2 \in B(\Z)$,
\begin{equation}
    \label{Eq: Cpqr}
     \vartheta_{m_1} \vartheta_{m_2}=\sum_{m \in B(\Z)} C_{m_1m_2}^m \vartheta_m \,,
\end{equation}
and for all $k \geq 0$ there are only finitely many $m \in B(\Z)$ such that 
$C_{m_1m_2}^m \neq 0 \mod I_k$. 
For every $m_1,m_2,m \in B(\Z)$ and a general point $p \in B$ sufficiently close to $m$, it is shown in \cite[Theorem 3.24]{GHS} that
\begin{equation}
    \label{Eq: formula for structure constants}
    C_{m_1m_2}^m = \sum_{\beta_1,\beta_2} a_{\beta_1} a_{\beta_2}
\end{equation}
where the sum is over all pairs $(\beta_1,\beta_2)$ of broken lines for $\foD$ with asymptotic directions $m_1,m_2$, end point $p$, and such that $m_{\beta_1}+ m_{\beta_2} = m$ where $a_{\beta_i} z^{m_{\beta_i}}$, for $i\in \{1,2\}$, are the final monomials carried by the broken lines $\beta_i$'s, as in Definition \ref{Def: broken line}.

Finally, we recall from \cite[Theorem 3.19]{GHS} that theta functions have natural power series expansions. For every general point $p \in B_0$, contained in a maximal cone $\sigma \in\P$, there exists a morphism of $\mathbf{k}[\![Q]\!]$-algebras 
\begin{align*} R(\foD) &\longrightarrow \mathbf{k}[\Lambda_\sigma][\![Q]\!]\\
\vartheta_m &\longmapsto \vartheta_m(p)
\end{align*}
such that
 \begin{equation}
     \label{Eq: theta function defined by broken line}
     \nonumber
     \vartheta_m(p) := \sum_{\beta} a_\beta z^{m_\beta}\,,
 \end{equation}
where the sum runs over all broken lines $\beta$ with asymptotic direction $m$ and endpoint $p$, and where $a_{\beta} z^{m_\beta}$ 
is the final monomial of $\beta$, as in Definition \ref{Def: broken line}.

\subsection{Canonical scattering diagram and Gross-Siebert mirror families}
\label{Sec: canonical scattering}
Let $(X,D)$ be a log Calabi-Yau pair, consisting of an $n$-dimensional smooth projective variety $X$ and a reduced simple normal crossing anticanonical divisor $D$ in $X$\footnote{For a more general notion of log Calabi-Yau pair, which is not necessary in the context of the current paper, see \cite{gross2021canonical}.
}. Under the assumption that $D$ contains a $0$-dimensional stratum, Gross and Siebert construct in \cite{gross2021canonical} 
the \emph{canonical scattering diagram} $\foD_{(X,D)}$ using the enumerative geometry of rational curves in $(X,D)$. In this situation, the algebra of functions on the mirror family to $(X,D)$ is proposed to be the algebra of theta functions defined by $\foD_{(X,D)}$. In this section, we briefly review the definition of the canonical scattering diagram and the construction of the mirror family.

\subsubsection{The tropicalization $(B,\P)$ of $(X,D)$}
\label{Sec: tropicalization}
We first review how to define from $(X,D)$ a conical affine pseudomanifold $(B,\P)$ as in Definition \ref{def_pseudomanifold} called the \emph{tropicalization of $(X,D)$}.

Let $D_1, \dots, D_m$ be the irreducible components of $D$. We assume that for every $I \subset \{1,\cdots,m\}$, the stratum $\bigcap_{i \in I} D_i$ is connected. Then, $(B, \P)$ is the dual intersection complex of $(X,D)$: for every $I \subset \{1,\cdots,m\}$ with $\bigcap_{i \in I} D_i$ non-empty, $\P$ contains an $|I|$-dimensional simplicial cone,
where $|I|$ is the cardinality of $I$, and these cones are glued together to form $B$
according to the intersection pattern of the strata of $D$. 
Moreover, one can extend the integral affine structure across codimension one cones using a recipe modeled on toric geometry:
given a codimension one cone $\rho$ with generators $m_1, \dots, m_{n-1}$, contained in two maximal cones $\sigma^+$ and $\sigma^-$
with additional generators $m_n^+$ and $m_n^-$, an integral affine structure is defined on $\sigma^+ \cup \sigma^-$ by embedding 
$\sigma^+ \cup \sigma^-$ in $\RR^n$ in such a way that 
\[ m_n^++m_n^-= -\sum_{i=1}^{n-1}(D_\rho \cdot D_{m_i}) m_i \,,\]
where $D_\rho$ is the curve stratum corresponding to the codimension one cone $\rho$, and 
$D_{m_i}$ are the divisor strata corresponding to the rays $\RR_{\geq 0}m_i$. 
Assuming that $D$ contains a $0$-dimensional stratum, it follows from 
\cite[Propositions 1.3-1.6]{gross2021canonical} that $(B,\P)$
is an $n$-dimensional conical affine pseudomanifold as in Definition \ref{def_pseudomanifold}. For the HDTV log Calabi-Yau pairs that we will consider later in
\S \ref{Sec: HDTV review}, this will be clear as $(B,\P)$
will be piecewise-linearly isomorphic to a fan in a vector space in these cases 
(see \eqref{eq_PL_isom}).

\subsubsection{The canonical scattering diagram}

We review below the definition of the canonical scattering diagram of $(X,D)$ constructed in \cite{gross2021canonical}. It is a scattering diagram in the sense of Definition \ref{def_scattering_general}, where:
\begin{itemize}
    \item[(i)] $(B,\P)$ is the tropicalization of $(X,D)$ described in 
    \S \ref{Sec: tropicalization}.
    \item[(ii)] $Q=NE(X)$ is the monoid spanned by effective curve classes in the abelian group $N_1(X)$ of curves classes in $X$ modulo numerical equivalence. Note that $NE(X)$ satisfies the conditions on $Q$ listed in \S 
    \ref{Sec:scatt}. Indeed, $NE(X)$ is contained in a strictly convex cone of $N_1(X)\otimes \RR$ because $X$ is projective (if $L$ is an ample divisor, then $L \cdot C >0$ for every $C$ in the closure of $NE(X)$ by Kleiman's criterion).
    \item[(iii)] $\varphi$ is a $NE(X)$-convex MVPL function on $(B,\varphi)$ with kink
    \begin{equation}
    \label{Eq: kink}
    \kappa_{\rho}=D_{\rho} \in NE(X),
\end{equation}
across every codimension one cone $\rho \in \P$, where $D_\rho$ is the 
class of the curve stratum corresponding to $\rho$. 
\end{itemize}
 
The definition of the canonical scattering diagram is based on the enumerative geometry of maps $f \colon (C,x) \rightarrow X $ from rational curves $C$ to $X$ with a prescribed tangency condition along $D$ at a given marked point $x$ on $C$. Moduli spaces of rational curves in $X$ meeting $D$ at a single point have expected dimension $n-2$ \cite[Lemma 3.9]{gross2021canonical}. The enumerative invariants entering the definition of the canonical scattering diagram roughly count $0$-dimensional families of maximally degenerated configurations of such curves. The precise definition \cite[\S 3.2]{gross2021canonical} is based on logarithmic geometry and the theory of \emph{punctured log Gromov--Witten invariants} \cite{ACGS}. 

These invariants, denoted
$N_{\tau,\beta}$, are indexed by an $(n-2)$-dimensional family $\tau$ of tropical curves in $(B,\P)$ called a \emph{wall type} \cite[Definition 3.6]{gross2021canonical} and a class $\beta \in NE(X)$. The tropical curves in $\tau$ are dual intersection complexes of the curve $(C,x)$ mapping to the dual intersection complex $(B,\P)$ of $(X,D)$. In particular, they have a single unbounded edge corresponding to the marked point $x \in C$. The $(n-2)$-dimensional family of these unbounded edges traces out an $(n-1)$-dimensional cone $\fod_\tau$ in $B$ contained in an $n$ or $(n-1)$-dimensional cone of $\P$. The direction $u_\tau \in \Lambda_{\fod_\tau}$ of these unbounded edges prescribes the contact order of $C$ with $D$ at the marked point $x \in C$. Finally, a multiplicity $k_\tau$ is attached to $\tau$, equal to the lattice index in $\Lambda_{\fod_\tau}$ of the sublattice
given by the image of the integral tangent vectors to 
the $(n-1)$-dimensional cone formed by the unbounded edges in the family of source tropical curves \cite[(3.10)]{gross2021canonical}.

\begin{definition}
\label{Def: canonical scattering}
The \emph{canonical scattering diagram associated to $(X,D)$}, denoted by $\foD_{(X,D)}$, is the scattering diagram in the tropicalization $(B,\P)$ of $(X,D)$ over $(NE(X),\varphi)$, given by the union of walls $\{(\fod_\tau,f_{\tau,\beta}) \}$ with
\begin{equation}
\nonumber
f_{\tau,\beta} := \exp(k_{\tau}N_{\tau, \beta}t^{\beta} z^{-u_\tau}) \in \textbf{k} [ \Lambda_{\fod_\tau}] \lfor NE(X) \rfor  \, ,
\end{equation}
indexed by wall types $\tau$ and curve classes $\beta \in NE(X)$, where 
$N_{\tau,\beta}$ is the punctured log Gromov--Witten invariant counting rational curves in $(X,D)$ of type $\tau$ and class $\beta$.
\end{definition}

\subsubsection{The Gross-Siebert mirror family}
\label{Sec:mirror family}
According to \cite[Theorem B]{gross2021canonical}, the canonical scattering diagram of Definition \ref{Def: canonical scattering} is consistent and so one can apply the constructions reviewed in
\S \ref{Sec:theta_functions}. In particular, we have the 
$\mathbf{k}[\![NE(X)]\!]$-algebra 
$R(\foD_{(X,D)})$ of theta functions defined by the canonical scattering diagram $\foD_{(X,D)}$. 
We denote 
\[ R(X,D):=R(\foD_{(X,D)})\,,\]
and the \emph{mirror family} of $(X,D)$ defined in \cite{gross2021canonical} is 
\[ \check{\mathfrak{X}} \longrightarrow \mathrm{Spf}\, \mathbf{k}[\![NE(X)]\!] \,, \]
where 
\[\check{\mathfrak{X}}:=\mathrm{Spf} \, R(X,D)\] is the formal scheme defined as the formal spectrum of the ring of theta functions $R(X,D)$.

\begin{example}
\label{Mirror to the blowup of P2 at a point}
Let $X$ be the blow-up of $X_\Sigma=\PP^2$ at a non-toric point of the toric boundary divisor $D_{\Sigma} \subset \PP^2$, and $D$ be the strict transform of $D_{\Sigma}$. In this situation the complement $X\setminus D$ is a cluster variety, known as the $A_1$ cluster variety with principal coefficients. We illustrate the canonical scattering diagram associated to $(X,D)$ in Figure \ref{Fig:figcan} -- for details on how to obtain the tropicalization $(B,\P)$ see for instance \cite[Example 3.2]{arguz2021heart}.
\begin{figure}
\center{\scalebox{.4}{\input{figcan.pspdftex}}}
\caption{On the left hand figure is the canonical scattering diagram for the blow up or $\PP^2$ at a non-toric point, and on the right hand figure we illustrate the broken lines defining the theta functions generating the coordinate ring of its mirror}
\label{Fig:figcan}
\end{figure}
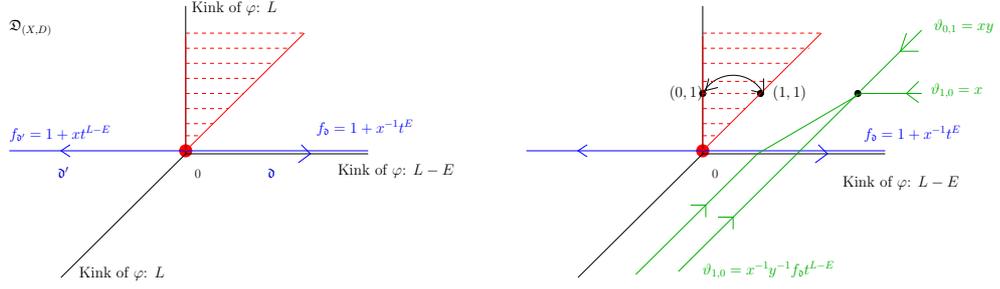
The three theta functions generating the ring of theta functions $R(X,D)$ are given by 
\begin{equation}
\nonumber
   \vartheta_{(1,0)} = x, \,\ \,\  \vartheta_{(0,1)} = xy, \,\ \,\ \mathrm{and} \,\ \,\ \vartheta_{(-1,-1)} = x^{-1}y^{-1}(1+x^{-1}t^{E})t^{L-E} \, , 
\end{equation}
where $L$ is the pullback of the class of a line in $\PP^2$, and by $E$ is the exceptional divisor. Hence, in this case, the mirror to $X\setminus D$ is given by
\begin{equation}
    \label{Eq: mirror equation for BlP2}
    \Spec \mathbf{k} [NE(X)] [ \vartheta_{(1,0)} , \vartheta_{(0,1)},  \vartheta_{(-1,-1)}   ] /  ( \vartheta_{(1,0)} \vartheta_{(0,1)} \vartheta_{(-1,-1)}  =  t^{L} + \vartheta_{(1,0)} t^{L-E})\,.  
\end{equation}
The equation $ \vartheta_{(1,0)} \vartheta_{(0,1)} \vartheta_{(-1,-1)}  =  t^{L} + \vartheta_{(1,0)} t^{L-E}$ involves only the variables $t^L$ and $t^{L-E}$ which span the monoid $\mathcal{M} \subset NE(X)$ discussed in \S\ref{Sec:Extension}.
\end{example}

\subsection{HDTV log Calabi-Yau pairs and scattering diagrams}
\label{Sec: HDTV review}

\subsubsection{HDTV log Calabi-Yau pairs}
\label{Sec:HDTV_log_CY}    
We first review the set-up of \cite{HDTV}.
Let $M \simeq \Z^n$ be a rank $n$ lattice and let $X_\Sigma$ be an $n$-dimensional smooth projective toric variety defined by a complete fan $\Sigma$ in $M_\RR := M \otimes \RR$. Let 
\[ H=\bigcup_{i \in I} H_i\]
be a union of disjoint connected smooth projective hypersurfaces of the toric boundary $D_\Sigma$ of $X_\Sigma$ indexed by a finite set $I$. We assume that for every $i\in I$, there is a unique irreducible component $D_i$ of $D_\Sigma$ such that $H_i \subset D_i$, that $H_i$ intersects transversally the toric boundary of $D_i$, and that for every $i \neq j$, either $D_i=D_j$ or $D_i \cap D_j=\emptyset$.

Let $X$ be the smooth projective variety obtained by blowing-up $X_\Sigma$ along $H$ and let $D \subset X$ be the strict transform of $D_\Sigma$. Then, $(X,D)$ is a \emph{log Calabi-Yau pair}, that is, $X$ is a smooth projective variety and $D$ is an anticanonical reduced simple normal crossings divisor on $X$. Because these log Calabi-Yau pairs are exactly those studied in \cite{HDTV}, we refer to such $(X,D)$ as a \emph{HDTV log Calabi-Yau pair} and we call
$(X_\Sigma,H)$ a \emph{toric model} for $(X,D)$.

The general construction of the mirror family based on the canonical scattering diagram $\foD_{(X,\Sigma)}$ and reviewed in \S \ref{Sec:mirror family} can be applied in particular to HDTV log Calabi-Yau pairs.
We showed in \cite{HDTV,arguz2021heart} how to recover the canonical scattering 
$\foD_{(X,\Sigma)}$ of a log Calabi-Yau pair from a combinatorially defined scattering diagram $\foD_{(X_\Sigma,H)}$, that we call the \emph{HTDV scattering diagram}. In the following sections, we first review the notion of a scattering diagram in $M_\RR$, then the definition of $\foD_{(X_\Sigma,H)}$ which is an example of scattering diagram in $M_\RR$, and finally the main result of \cite{HDTV} comparing $\foD_{(X,D)}$ and $\foD_{(X_\Sigma,H)}$. 

\subsubsection{Scattering diagrams in $M_\RR$}
\label{Sec:scattering_M_R}

Scattering diagrams in $M_\RR$ are particular examples of general scattering diagrams as in Definition \ref{def_scattering_general}.

\begin{definition} \label{def_scattering}
Let $Q$ be a monoid as in \S\ref{Sec:scatt}. A  \emph{scattering diagram in $M_\RR$ over $Q$} is by definition a scattering diagram in $(B,\P)$
over $(Q,\varphi)$ as in Definition \ref{def_scattering_general} where $B=M_\RR$, 
$\P$ is the trivial conical decomposition of $M_\RR$ consisting of the single cone $M_\RR$, and $\varphi=0$. 

In other words, a scattering diagram in $M_\RR$ over $Q$ is a set of walls $(\fod,f_{\fod})$ in $M_\RR$,
that is pairs 
$(\fod,f_{\fod})$, consisting of a codimension one rational polyhedral cone $\fod \subset M_{\RR}$, together with an attached functions $f_\fod \in \mathbf{k}[z^{-m_0}][\![Q]\!] \subset
\mathbf{k}[M][\![Q]\!]$
for some nonzero primitive vector $m_0 \in M$ tangent to $\fod$, and with $f_\fod=1 \mod I_1$.
Moreover, for every $k \geq 1$, there are only a finite number of walls $(\fod,f_{\fod}) \in \foD$ with $f_{\fod} \neq 1 \mod I_k$.
A wall is called \emph{incoming} if 
$-m_0 \in \fod$, and \emph{outgoing} elsewise. We call $m_0$ the \emph{direction} of the wall.
\end{definition}

For the construction of the HDTV scattering diagram and the comparison with cluster scattering diagram, the notion of a consistent scattering diagram in $M_\RR$ is particularly important and so we  review the definition, following \cite{HDTV,gross2021canonical}.
Given a scattering diagram $\foD$ in $M_\RR$ over $Q$ and a path \begin{align*} \gamma: [0,1] &\longrightarrow M_{\RR} \setminus \mathrm{Sing}(\foD)\\ 
 t &\longmapsto \gamma(t)\,\end{align*}
 transversal to the walls of $\foD$, 
the associated path ordered product $\mathfrak{\theta}_{\gamma, \foD}$  is the $\mathbf{k}[\![Q]\!]$-algebra automorphism of $\mathbf{k}[M][\![Q]\!]$ obtained as the ordered product of the $\mathbf{k}[\![Q]\!]$-algebra automorphisms 
$\theta_{\gamma, \fod_i}$ attached to the sequence of walls $(\fod_i, f_{\fod_i})$ crossed by $\gamma$ for $t=t_i$, where 
\begin{align} \label{eq_autom}
\theta_{\gamma, \fod_i} \colon \mathbf{k}[M][\![Q]\!]&\longrightarrow \mathbf{k}[M][\![Q]\!] \\
 \nonumber
 z^m &\longmapsto f_{\fod_i}^{\langle n_{\fod,\gamma},m\rangle} z^m \,,
\end{align}
where $n_{\fod,\gamma} \in N:=\Hom(M,\Z)$ is the primitive normal vector of $\fod$ which is negative on 
$\gamma((t_0-\epsilon, t_0))$ for all small $\epsilon >0$.
Finally, we call a scattering diagram \emph{consistent} if for any path $\gamma$ with $\gamma(0)=\gamma(1)$ 
the associated path ordered product $\theta_{\gamma,\foD}$ is the identity automorphism.

According to \cite[Theorem 5.6]{HDTV}\footnote{Similar reconstruction results of consistent scattering diagrams from initial walls were also obtained in \cite{GSannals,KS2000}.}, one can construct 
consistent scattering diagrams in $M_\RR$ from the data of particular collections of initial walls called \emph{widgets} in \cite{HDTV} -- for details see \cite[\S5.2.1]{HDTV}. We review below the notion of widget, which is based on the notion of tropical hypersurface.

\begin{definition}
\label{Def: tropical hypersurface}
A \emph{tropical hypersurface} in $M_{\RR}$
is a fan $\T$ in $M_{\RR}$ whose support $|\T|$ is pure dimension $\dim M_{\RR}-1$, along
with a positive integer weight attached to each cone of $\T$ of
dimension $\dim M_{\RR}-1$, which satisfies the following
balancing condition.
For every $\omega\in\T$ of dimension $\dim M_{\RR}-2$, let $\gamma$
be a loop in $M_{\RR}\setminus \omega$ around an interior point
of $\omega$, intersecting top-dimensional 
cones $\sigma_1,\ldots,\sigma_p$ of $\T$ of weights $w_1,\ldots,w_p$.
Let $n_i\in N:= \Hom(M,\Z)$ be the primitive normal vector of $\sigma_i$ which is negative on 
$\gamma((t_0-\epsilon, t_0))$ for all small $\epsilon >0$.
Then
\begin{equation}
\label{Eq:balancing}
\nonumber
\sum_{i=1}^p w_in_i=0\,.
\end{equation}
\end{definition}

\begin{definition}
\label{def:widget_general}
Let $\Sigma$ be a complete toric fan in $M_\RR$ and $v \in M \setminus \{0\}$ be a primitive integral vector such that $\RR_{\geq 0}v$
is a ray of $\Sigma$.
Denote by $\pi \colon M_\RR \rightarrow M_\RR/\RR v$ the quotient map
and by $\Sigma(v):=\{ \pi(\sigma)\,|\,\sigma \in \Sigma\,,\, v \in \sigma \}$ the quotient fan of $\Sigma$ in $M_\RR/\RR v$.
Let $\T$ be a tropical hypersurface in $M_\RR/\RR v$ with support contained in the union of the codimension one cones of $\Sigma(v)$. Finally, let $f \in \mathbf{k}[z^{v}][\![Q]\!]$ such that 
$f=1 \mod I_1$.
Then, the \emph{widget} associated to $\T$ and $f$ is the scattering diagram in $M_\RR$ over $Q$ given by
\begin{equation}
\nonumber
\foD :=\{(\fod_{\sigma}, f^{w_\sigma} ) \,|\,\sigma ~\mathrm{is  ~ a ~maximal ~ dimensional ~ cone ~ of ~ } \T\}\,.
\end{equation}
Here, $\fod_{\sigma}$ is the unique codimension one cone of $\Sigma$ containing $\RR_{\geq 0}v$ such that $\pi(\fod_{\sigma})=\sigma$,
and $w_\sigma$ is the weight of $\sigma$ in $\T$.
\end{definition}

\begin{remark}
In Definition \ref{def:widget_general}, the walls of the widget have directions $-v$ and are incoming.
\end{remark}

Let $\Sigma$ be a complete toric fan in $M_\RR$, and $(v_i)_{i\in I}$ be
a finite collection of primitive integral vectors in $M \setminus \{0\}$, indexed by a finite set $I$, such that $\RR_{\geq 0}v_i$ is a ray of $\Sigma$ for all $i \in I$. For every $i \in I$, let $\T_i$ be a tropical hypersurface in $M_\RR/\RR v_i$ with support contained in the union of codimension one cones of $\Sigma(v_i)$, and let 
$f_i \in \mathbf{k}[z^{-v_i}][\![Q]\!]$ such that 
$f_i=1 \mod I_1$.
Finally, let 
$\foD_{\mathrm{in}}= \bigcup_{i \in I} \foD_i$ be the \emph{initial scattering diagram} obtained by taking the union of the widgets $\foD_i$ associated to $\T_i$ and $f_i$ as in Definition \ref{def:widget_general}.
According to \cite[Theorem 5.6]{HDTV}, there exists a consistent scattering diagram $S(\foD_{\mathrm{in}})$ containing $\foD_{\mathrm{in}}$, such that all walls in $S(\foD_{\mathrm{in}}) \setminus \foD_{\mathrm{in}}$ are outgoing. Moreover,
$S(\foD_{\mathrm{in}})$ is unique up to equivalence with these properties and we refer to $S(\foD_{\mathrm{in}})$ as the \emph{consistent completion of $\foD_{\mathrm{in}}$}.

\subsubsection{The HDTV scattering diagram}
\label{Sec: HDTV scattering}
In the remaining part of this subsection we review the construction of the combinatorially constructed HDTV scattering diagram $\foD_{(X_\Sigma,H)}$ introduced in \cite{HDTV}.

Let $(X,D)$ be a HDTV log Calabi-Yau pair obtained as in 
\S \ref{Sec:HDTV_log_CY} from a toric model $(X_\Sigma, H)$, where 
$H = \cup_{i \in I} H_i$. 
The HDTV scattering diagram $\foD_{(X_\Sigma,H)}$ is a scattering diagram in $M_\RR$, as in Definition \ref{def_scattering}, over the monoid $\NN^I$ containing one copy of $\NN$ for each hypersurface $H_i$.

In \cite{HDTV}, the HDTV scattering diagram $\foD_{(X_{\Sigma},H)}$ is defined by first describing an initial scattering diagram $\foD_{(X_{\Sigma},H),\mathrm{in}}$.
The construction of $\foD_{(X_{\Sigma},H),\mathrm{in}}$ is based on the tropical hypersurfaces defined as follows.

\begin{definition} \label{Def: tropical hypersurface associated}
For every $i\in I$, let $v_i \in M$ be the primitive generator of the ray of $\Sigma$ corresponding to the toric divisor $D_i$ containing $H_i$.  
The \emph{tropical hypersurface associated to $H_i$} is the tropical hypersurface  $\scrH_{i}\subseteq M_\RR/\RR v_i$
which is supported on the union of codimension one cones of the toric fan 
$\Sigma(v_i)$ in $M_\RR/\RR v_i$
of the divisor $D_{i}$, 
with the weight on a codimension one cone $(\sigma+\RR v_i)/\RR v_i$ being
\begin{equation} \label{eq_weight}
w_{\sigma}=D_{\sigma}\cdot H_i\,,\end{equation}
where the intersection number
is calculated on $D_{i}$. Here $\sigma$ denotes a codimension one cone of $\Sigma$ containing the ray $\RR_{\geq 0}v_i$ and $D_{\sigma}$ is the toric curve in $D_{i}$ corresponding to $\sigma$. 
\end{definition}

We now define the initial HDTV scattering diagram.

\begin{definition}
\label{Def HDTV scattering initial}
The \emph{initial HDTV scattering diagram} for $(X_\Sigma,H)$ is the scattering diagram in $M_{\RR}$ over $\NN^I$, as in Definition 
\ref{def_scattering},
whose set of walls is given by  
\[ 
\foD_{{(X_\Sigma,H)},\mathrm{in}} := \bigcup_{i \in I} \foD_i 
\,,\]
where $\foD_i$ is the widget associated as in Definition 
\ref{def:widget_general}
to the tropical hypersurface $\scrH_i$ in $M_\RR/\RR v_i$ and to the function 
\begin{equation}
\nonumber
    f_i :=1+t_i z^{v_i} \in \mathbf{k}[M][\![\NN^I]\!]\,,
\end{equation}
where $t_i$ is the variable in the monoid algebra $\mathbf{k}[\NN^I]$
corresponding to the generator of the copy of $\NN$ labeled by $i\in I$.
\end{definition}

As reviewed at the end of \S \ref{Sec:scattering_M_R},
every initial scattering diagram $\foD_{\mathrm{in}}$ defined as a union of widgets has a consistent completion $S(\foD_{\mathrm{in}})$.

\begin{definition}
\label{Def HDTV scattering}
The \emph{HDTV scattering diagram} for $(X_\Sigma, H)$ 
\[\foD_{(X_\Sigma,H)} :=  S(\foD_{(X_{\Sigma},H),\mathrm{in}})   \, ,   \] 
is the scattering diagram in $M_\RR$ over $\NN^I$ obtained as the consistent completion of $\foD_{(X_{\Sigma},H),\mathrm{in}}$.
\end{definition}

\subsubsection{Comparison of the HDTV and canonical scattering diagrams}

Let $(X,D)$ be a HDTV log Calabi-Yau pair obtained as in 
\S \ref{Sec:HDTV_log_CY} from a toric model $(X_\Sigma, H)$.
We review the main result of \cite{HDTV} comparing the canonical scattering diagram $\foD_{(X,D)}$ and the HDTV scattering diagram 
$\foD_{(X_\Sigma,H)}$.

We first introduce some notations. Denote 
$\mathrm{Bl}_H \colon X \rightarrow X_\Sigma$ the blow-up morphism.
For every $i\in I$, the exceptional divisor 
$\mathcal{E}_i$ in $X$ over $H_i$ is a $\PP^1$-bundle over $H_i$, and we denote
by $E_i$ the class in $NE(X)$ of a $\PP^1$-fiber. In particular, for every $i,j \in I$, we have $\mathcal{E}_i \cdot E_j=-1$ if $i=j$, and $\mathcal{E}_i \cdot E_j=0$ if $i \neq j$. 
The map 
\begin{align}\nonumber
    \iota : N_1(X_\Sigma)\oplus \Z^I & \longrightarrow N_1(X) \\
 \nonumber   (C, (a_i)_{i\in I}) &\longmapsto \mathrm{Bl}_H^{*}C -\sum_{i\in I}a_i E_i \,,
\end{align}
is an isomorphism of abelian groups, whose inverse is 
\begin{align}\nonumber
      N_1(X) &\longrightarrow N_1(X_\Sigma)\oplus \Z^I \\
\nonumber    C &\longmapsto (\mathrm{Bl}_{H,*}C, (C \cdot \mathcal{E}_i)_{i\in I}) \,.
\end{align}
From now one, we implicitly use $\iota$ to identify 
$N_1(X)$ with $N_1(X_\Sigma) \oplus \Z^I$.

To compare the canonical scattering diagram $\foD_{(X,D)}$
with $\foD_{(X_{\Sigma},H)}$ first note that it follows from the definition of the tropicalization of
$(X,D)$ reviewed in \S \ref{Sec: tropicalization}
that there is a natural piecewise-linear isomorphism
\begin{equation} \label{eq_PL_isom}
\Upsilon:(M_{\RR},\Sigma)\rightarrow (B,\P)\,.
\end{equation}
In particular, $\Upsilon$ induces a bijection between $M$ and the set 
$B(\Z)$ of integral points as in Definition \ref{def_asymptotic_direction}, and from now on we identify $M$ and $B(\Z)$ using this bijection.
We use $\Upsilon$ to define from the HDTV scattering diagram 
$\foD_{(X_\Sigma,H)}$
in $M_\RR$ a new scattering diagram $\Upsilon_{*}(\foD_{(X_\Sigma,H)})$ in 
$(B,\P)$. The construction treats differently incoming and outgoing walls. Up to refining the walls, we may assume that every wall 
$(\fod, f_\fod)$ of 
$\foD_{(X_\Sigma,H)}$ is contained in some cone $\sigma \in \Sigma$.

If the wall $(\fod, f_\fod)$ is incoming, then by construction of $\foD_{(X_{\Sigma},H)}$ (see Definitions \ref{def:widget_general} and \ref{Def HDTV scattering initial})
it is of the form
\[ (\fod,(1+t_i z^{v_i})^w) \] 
for some positive integer
$w$. As $v_i$ is tangent to the cone of $\Sigma$ containing
$\fod$ and $\Upsilon$ is piecewise linear with respect to $\Sigma$,
$\Upsilon_*(v_i)$ makes sense as a tangent vector to $B$.
We then define 
\begin{equation}
    \label{Eq: Upsilon}
    \Upsilon_*(\fod,(1+t_i z^{v_i})^{w})
=(\Upsilon(\fod),(1+t^{E_i}z^{-\Upsilon_*(v_i)})^w)\,.
\end{equation}
If the wall $(\fod, f_\fod)$ is outgoing, with $\fod\subseteq\sigma\in\Sigma$, then
$f_\fod$ is a sum of monomials of the form 
\[ c \prod_{i \in I} (t_i z^{v_i})^{a_i}\]
with $c \in \mathbf{k}$, and $a_i$ non-negative integers for all $i \in I$.
The data of the list of integers $\mathbf{A}=\{a_{i}\}_{i\in I}$ and $\sigma$ determine a curve class 
$\bar\beta_{\mathbf{A},\sigma}\in N_1(X_{\Sigma})$ -- see \cite[\S$6$]{HDTV} for the precise description of this curve
class. Under the inclusion $N_1(X_{\Sigma}) \hookrightarrow N_1(X)$ 
given by the above mentioned splitting, we may view
$\bar\beta_{\mathbf{A},\sigma}$ as a curve class in $N_1(X)$, which we also denote by $\bar\beta_{\mathbf{A},\sigma}$. We then obtain a curve class
\begin{equation}
    \label{Eq: curve class}
    \beta_{\mathbf{A},\sigma}=\bar\beta_{\mathbf{A},\sigma} - \sum_{i \in I} a_{i} E_i \,.
\end{equation}
Further, as $v_{\out}:=-\sum_{i \in I} a_i v_i$ is tangent to $\sigma$, as before $\Upsilon_*(v_{\out})$ makes sense
as a tangent vector to $B$. We may thus define 
\begin{equation} \label{Eq: upsilon star}
\Upsilon_* (c \prod_{i\in I}(t_i z^{v_i})^{a_i})
:= c \,t^{\beta_{\mathbf{A},\sigma}}z^{-\Upsilon_*(v_{\out})}\,,
\end{equation}
then $\Upsilon_*(f_\fod)$ by linearity, and finally
the wall
\begin{equation}
\label{Eq: wall upsilon canonical}
\Upsilon_*(\fod,f_{\fod})=(\Upsilon(\fod),\Upsilon_*(f_\fod))\,.
\end{equation}
We then define  
\[
\Upsilon_*(\foD_{(X_{\Sigma},H)}) :=
\{\Upsilon_*(\fod,f_{\fod})\,|\, (\fod,f_{\fod})\in \foD_{(X_{\Sigma},H)}\}\,.
\]

The key theorem, \cite[Thm $6.1$]{HDTV}, then states:

\begin{theorem} \label{thm_canonical_hdtv}
Let $(X,D)$ be a HDTV log Calabi-Yau pair obtained as in 
\S \ref{Sec:HDTV_log_CY} from a toric model $(X_\Sigma, H)$.
Then, the scattering diagram
$\Upsilon_*(\foD_{(X_{\Sigma},H)})$ obtained from the HDTV scattering diagram 
$\foD_{(X_\Sigma,H)}$ by applying $\Upsilon_{*}$
is equivalent to the canonical scattering diagram $\foD_{(X,D)}$.
\end{theorem}

\section{The extension of the Gross--Siebert mirror family}
\label{Sec:Extension}

\subsection{Construction of the extended mirror family}
Let $(X,D)$ be a HDTV log Calabi-Yau pair obtained as in 
\S \ref{Sec:HDTV_log_CY} from a toric model $(X_\Sigma, H)$.
In this section, we prove that the mirror family
$\check{\mathfrak{X}} \rightarrow \mathrm{Spf}\, \mathbf{k}[\![NE(X)]\!]$
of $(X,D)$ naturally extends over a bigger base. 
The idea is to replace the monoid $NE(X)$ by the smaller monoid
\begin{equation}
    \label{M}
    \mathcal{M}:= NE(X) \cap (NE(X_\Sigma) \oplus \NN^I) \subset N_1(X) 
    \end{equation}
of effective curve classes of the form $\mathrm{Bl}_H^{*}\overline{\beta}-\sum_{i\in I}a_i E_i$, where $\mathrm{Bl}_H: X \rightarrow X_\Sigma$ is the blow-up morphism, $\overline{\beta}\in NE(X_\Sigma)$,  and $a_i \in \NN$ for all $i\in I$.

The projection $\mathcal{M} \subset NE(X_\Sigma) \oplus \NN^I \rightarrow NE(X_\Sigma)$ induces a closed embedding 
\[ \Spec \mathbf{k}[NE(X_\Sigma)] \hooklongrightarrow \Spec \mathbf{k}[\mathcal{M}]\,, \]
defined by the monoid ideal 
\begin{equation} \label{eq_J}
J:= \{ \beta \in \mathcal{M}\,|\, \beta \notin NE(X_\Sigma) \}\,.
\end{equation}
Let $S_{(X_\Sigma,H)}$ be the formal completion of $\Spec \mathbf{k}[\mathcal{M}]$ along  $\Spec \mathbf{k}[NE(X_\Sigma)]$, that is,
\[ S_{(X_\Sigma,H)}=\mathrm{Spf}\,\widehat{\mathbf{k}[\mathcal{M}]}\] where 
\[ \widehat{\mathbf{k}[\mathcal{M}]} := \varprojlim_k \mathbf{k}[\mathcal{M}]/J^k \,.\]
The inclusion $\mathcal{M} \subset NE(X)$ induces a morphism
\[ \mathrm{Spf}\, \mathbf{k}[\![NE(X)]\!] \longrightarrow S_{(X_\Sigma,H)}\,.\]

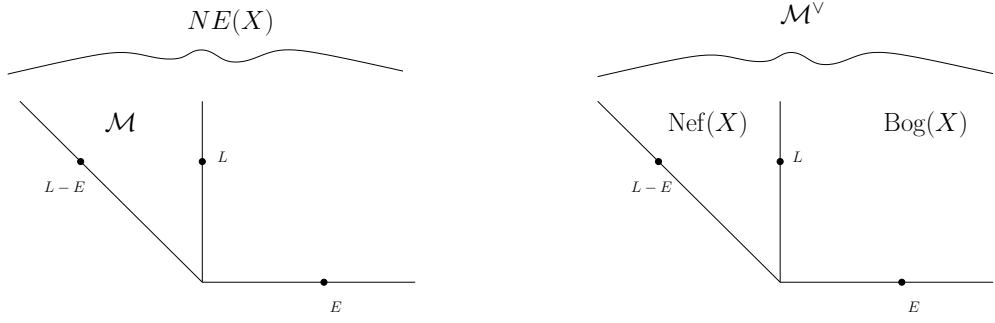
\begin{figure}
\center{\scalebox{.4}{\input{cones.pspdftex}}}
\caption{Cones of curves on the left and cones of divisors on the right
in Examples \ref{example_M} and \ref{example_cones}}
\label{Fig:cones}
\end{figure}

\begin{example}
\label{example_M}
Let $(X,D)$ be the log Calabi-Yau surface obtained by blowing-up a point in 
$\PP^2$, as in Example \ref{Mirror to the blowup of P2 at a point}. Recall that we denote by $L$ the pullback of the class of a line in $\PP^2$, and by $E$ the exceptional divisor. Then, 
$NE(X_\Sigma)=\NN L$ and $NE(X)=\NN E \oplus \NN (L-E)$, and so $\mathcal{M}= \NN L \oplus \NN (L-E)$, as illustrated in Figure \ref{Fig:cones}. 
\end{example}

Our first main result below shows that the mirror family $\check{\mathfrak{X}} \rightarrow \mathrm{Spf}\, \mathbf{k}[\![NE(X)]\!]$ naturally extends over $S_{(X_\Sigma,H)}$, and defines an \emph{extended} mirror family 
\[\check{\mathfrak{X}}_{\mathrm{ext}} \longrightarrow S_{(X_\Sigma,H)}\]
of $(X,D)$ with respect to the toric model $(X_\Sigma,H)$. 

\begin{theorem} \label{thm_intro_extension_1}
Let $(X,D)$ be a HDTV log Calabi-Yau pair obtained as in 
\S \ref{Sec:HDTV_log_CY} from a toric model $(X_\Sigma, H)$.
Then, the product of theta functions defines
a structure of topological $\widehat{\mathbf{k}[\mathcal{M}]}$-algebra on 
\[ R(X,D)_{\mathrm{ext}}:= \varprojlim_k \bigoplus_{m \in M} (\mathbf{k}[\mathcal{M}]/J^k) \, \vartheta_m\,,\]
such that, denoting $\check{\mathfrak{X}}_{\mathrm{ext}}:= \mathrm{Spf} R(X,D)_{\mathrm{ext}}$, we have a fiber diagram
\[\begin{tikzcd}
\check{\mathfrak{X}} 
\arrow[r]
\arrow[ d]
&
\check{\mathfrak{X}}_{\mathrm{ext}}
\arrow[d]\\
\mathrm{Spf}\, \mathbf{k}[\![NE(X)]\!]
\arrow[r]& S_{(X_\Sigma, H)}\,.
\end{tikzcd}\]
\end{theorem}

\begin{proof}
Recall from \eqref{Eq: Cpqr} that the algebra structure on the ring of theta functions is determined by a set of structure constants $C_{m_1m_2}^m \in \mathbf{k}[\![NE(X)]\!]$, for all
$m_1,m_2,m\in B(\Z)$,
such that, for all $k \geq 0$, there are only finitely many $m \in B(\Z)$ such that 
$C_{m_1m_2}^m \neq 0 \mod I_k$. 
By the definition of $\mathcal{M}$ in \eqref{M} and the definition of $S_{(X_\Sigma,H)}$ as the formal completion of $\Spec \mathbf{k}[\mathcal{M}]$ along  $\Spec \mathbf{k}[NE(X_\Sigma)]$, to prove Theorem \ref{thm_intro_extension_1}, it is enough to show that for every $m_1,m_2,m\in B(\Z)$, 
\[ C_{m_1m_2}^m \in \mathbf{k}[NE(X_\Sigma)][\![\NN^I]\!]\,,\]
and that for all $k \geq 0$ and for all $m_1, m_2 \in B(\Z)$, there are only finitely many $m \in B(\Z)$ such that $C_{m_1m_2}^m \neq 0 \mod \mathfrak{m}_{\NN^I}^k$, where $\mathfrak{m}_{\NN^I}:=\NN^I \setminus \{0\}$ is the maximal monoid ideal of $\NN^I$. 

To show $C_{m_1m_2}^m \in \mathbf{k}[NE(X_\Sigma)][\![\NN^I]\!]$, first recall that $C_{m_1m_2}^m$ can be expressed by \eqref{Eq: formula for structure constants}
in terms of final monomials carried by broken lines in the canonical scattering diagram $\foD_{(X,D)}$. Hence, it suffices to show that all the coefficients
$a_{\beta}$
appearing in \eqref{Eq: formula for structure constants} are in $\mathbf{k}[NE(X_\Sigma)][\![\NN^I]\!]$ and that for every $k \geq 0$, only finitely many of the coefficients $a_\beta$ are nonzero modulo 
$\mathfrak{m}_{\NN^I}^k$. We start proving the first claim by showing the stronger statement that for every broken line $\beta$ as in Definition \ref{Def: broken line}, we have  $a_i \in \mathbf{k}[NE(X_\Sigma)][\![\NN^I]\!]$
for every monomial $a_i z^{m_i}$ carried by any line
segment $\beta((t_{i-1},t_i))$ of $\beta$. 

We prove this by induction on $i$. By Definition \ref{Def: broken line}, the initial monomial is given by $z^{m_1}$ and so $a_1=1 \in  \mathbf{k} [NE(X_\Sigma)][\![\NN^I]\!]$. Now assume we have a monomial $a_i z^{m_i}$ carried by a broken line such that $a_i \in \mathbf{k}[NE(X_\Sigma)][\![\NN^I]\!]$. Then, the next monomial 
$a_{i+1} z^{m_{i+1}}$ is a monomial in the expansion of 
\eqref{eq_bending_1} or \eqref{eq_bending_2}.

If $\beta(t_i)$ does not belong to a codimension one cone of $\P$, then $a_{i+1}z^{m_{i+1}}$ is a monomial in the expansion of 
\eqref{eq_bending_1}. Moreover, by Theorem \ref{thm_canonical_hdtv}, all the walls of $\foD_{(X,D)}$
contributing to the bending are of the form 
$\Upsilon_* ((\fod, f_\fod))$, where 
$(\fod, f_\fod)$ is an outgoing wall of $\foD_{(X_{\Sigma},H)}$. 
In such case, $\Upsilon_* ((\fod, f_\fod))$ is defined by 
\eqref{Eq: upsilon star}-\eqref{Eq: wall upsilon canonical}, where the curve class $\beta_{\mathbf{A},\sigma}$ is defined by \eqref{Eq: curve class} and so is contained in $NE(X_{\Sigma}) \oplus \NN^I$.
Hence, the function 
$f$ in \eqref{eq_bending_1} belongs to $\mathbf{k}[\Lambda_{\sigma_i}][NE(X_{\Sigma}][\![\NN^I ]\!]$,
and so $a_{i+1} \in \mathbf{k}[NE(X_\Sigma)][\![\NN^I]\!]$.

If $\beta(t_i)$ does belong to a codimension one cone $\rho$ of $\P$, then $a_{i+1}z^{m_{i+1}}$ is a monomial in the expansion of \eqref{eq_bending_2}:
\[ a_i z^{m_i}(t^{\kappa_\rho} f)^{\langle n,m_i \rangle}\]
where $\langle n,m_i \rangle>0$, and, by \eqref{Eq: kink}, $\kappa_\rho$ is the class of the curve $D_\rho$ in $X$ corresponding to 
$\rho$.

If there are no $i \in I$ such that $v_i \in \rho$, then $D_\rho$ is simply the pullback of the class of the toric curve $\overline{D}_\rho$ in 
$X_\Sigma$ corresponding to $\rho$, that is $D_\rho \in NE(X_\Sigma)$, and as 
$\langle n,m_i \rangle>0$, we have $t^{\langle n,m_i \rangle \kappa_\rho}
\in \mathbf{k}[NE(X_\Sigma)]$. Moreover, all the walls contributing to the bending are images by $\Upsilon_{*}$ of outgoing walls of $\foD_{(X_{\Sigma},H)}$, and so 
\[ f^{\langle n,m_i \rangle} \in \mathbf{k}[\Lambda_{\sigma_i}][NE(X_\Sigma)][\![\NN^I]\!]\] as before. Hence, we also have $a_{i,i+1} \in \mathbf{k}[NE(X_\Sigma)][\![\NN^I]\!]$.

If there exists $i \in I$ such that $v_i \in \rho$, then we can write 
$f=f' f''$, where $f''$ is the contribution of walls images by $\Upsilon_*$
of outgoing walls of $\foD_{(X_{\Sigma},H)}$, and $f'$ is the image by 
$\Upsilon_*$ of an incoming wall of $\foD_{(X_{\Sigma},H)}$ with attached function \[ (1+t_i z^{v_i})^{w_\rho}\,,\]
as in Definitions \ref{def:widget_general} and \ref{Def HDTV scattering initial}.
As before, we have \[ (f'')^{\langle n,m_i \rangle} \in \mathbf{k}[\Lambda_{\sigma_i}][NE(X_\Sigma)][\![\NN^I]\!]\,.\] On the other hand, $f'$ is given by \eqref{Eq: Upsilon}, that is 
\[ f'=(1+t^{E_i}z^{-\Upsilon_*(v_i)})^{w_\rho}\,,\]
which does not belong to $\mathbf{k}[\Lambda_{\sigma_i}][NE(X_\Sigma)][\![\NN^I]\!]$ in general. However, the combination 
$t^{\kappa_\rho} f'$ does belong to  $\mathbf{k}[\Lambda_{\sigma_i}][NE(X_\Sigma)][\![\NN^I]\!]$. 
Indeed, $\kappa_\rho=D_\rho$ is the class of the strict transform of the toric curve $\overline{D}_\rho$ in $X_\Sigma$ corresponding to $\rho$, so 
\[ D_\rho= \overline{D}_\rho -(\overline{D}_\rho \cdot H_i)E_i \,,\]
and, as $\overline{D}_\rho \cdot H_i =w_\rho$ by \eqref{eq_weight}, we have 
\[ t^{\kappa_\rho} f'=t^{\overline{D}_\rho -w_\rho E_i}(1+t^{E_i}z^{-\Upsilon_*(v_i)})^{w_\rho} \in \mathbf{k}[\Lambda_{\sigma_i}][NE(X_\Sigma)][\![\NN^I]\!] \,.  \]
As 
$\langle n,m_i \rangle>0$, we also have $(t^{\kappa_\rho} f')^{\langle n,m_i \rangle}\in \mathbf{k}[\Lambda_{\sigma_i}][NE(X_\Sigma)][\![\NN^I]\!] $, and finally $a_{i,i+1} \in \mathbf{k}[NE(X_\Sigma)][\![\NN^I]\!]$.

It remains to show that for every $k \geq 0$, only finitely many of these monomials are nonzero modulo 
$\mathfrak{m}_{\NN^I}^k$. By Theorem
\ref{thm_canonical_hdtv}, it is enough to prove that the scattering diagram $\foD_{(X_\Sigma,H)}$ is finite modulo $\mathfrak{m}_{\NN^I}^k$. This holds because 
$\foD_{(X_\Sigma,H)}$ is a scattering diagram over 
$\mathbf{k}[\![\NN^I]\!]$, see Definition \ref{def_scattering} . Finally, the claim that for given $m_1,m_2 \in B(\Z)$, there are only finitely many $m \in B(\Z)$ such that 
$C_{m_1m_2}^m \neq 0 \mod I_k$, follows for the same reason. 
\end{proof}

The extended mirror family has the nice property to interpolate between the mirror family of the log Calabi-Yau pair $(X,D)$
and the mirror family of the toric model $(X_\Sigma, D_\Sigma)$.

\begin{theorem}\label{thm_toric_mirror}
Let $(X,D)$ be a HDTV log Calabi-Yau pair obtained as in 
\S \ref{Sec:HDTV_log_CY} from a toric model $(X_\Sigma, H)$.
The restriction of the extended mirror family 
$\check{\mathfrak{X}}_{\mathrm{ext}} \rightarrow S_{(X_\Sigma,H)}$ to 
$\Spec \mathbf{k}[NE(X_\Sigma)]$ is the mirror family
of the toric variety $(X_\Sigma, D_\Sigma)$.
In particular, the restriction of $\check{\mathfrak{X}}_{\mathrm{ext}}$
to the torus $\Spec \mathbf{k}[N_1(X_\Sigma)]$ is a family of tori $\Spec \mathbf{k}[M]$, and the theta functions $\{ \vartheta_m \}_{m \in M}$ restrict to the 
monomials $\{ z^m\}_{m \in M}$ on the fiber $\Spec \mathbf{k}[M]$ over $1 \in \Spec \mathbf{k}[NE(X_\Sigma)]$.
\end{theorem}

\begin{proof}
By setting $t^{-E_i}=0$ for all $i \in I$, the canonical scattering diagram of $(X,D)$ becomes the canonical scattering diagram of 
$(X_\Sigma,D_\Sigma)$: all the walls become trivial and the kinks reduce to the toric kinks.
\end{proof}

\begin{remark} \label{rem_affine_1}
Let $(X,D)$ be a HDTV log Calabi-Yau pair such that the complement $X \setminus D$ is affine. Then, by \cite[Corollary 1.2]{johnston2022comparison} (based on the comparison with \cite{KY}), the mirror family 
$\check{\mathfrak{X}} \rightarrow \mathrm{Spf}\, \mathbf{k}[\![NE(X)]\!]$ canonically extends over $\Spec \mathbf{k}[NE(X)]$. In this case,
the extended mirror family $\check{\mathfrak{X}}_{\mathrm{ext}} \rightarrow S_{(X_\Sigma,H)}$ given by Theorem \ref{thm_intro_extension_1} canonically extends further over $\Spec \mathbf{k}[\mathcal{M}]$.
\end{remark}

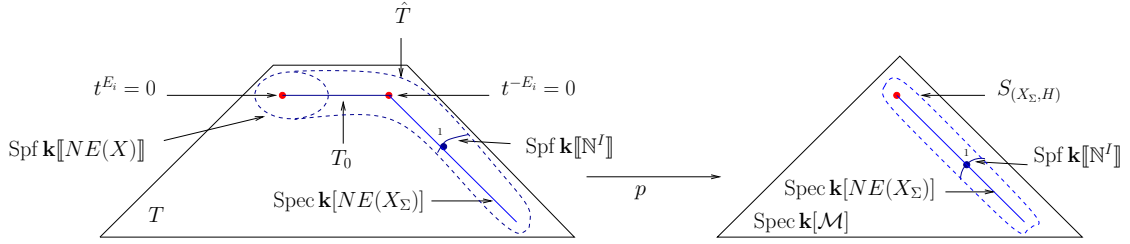
\begin{figure}
\center{\scalebox{.4}{\input{p.pspdftex}}}
\caption{The map $p: T \to \Spec \mathbf{k}[\mathcal{M}]$ contracting $T_0$ to the torus fixed point.}
\label{Fig:p}
\end{figure}

\subsection{Motivation from mirror symmetry for the extension}
In this section, we explain how the construction of the extended mirror family $ \check{\mathfrak{X}}_{\mathrm{ext}} \to S_{(X_\Sigma, H)}$ constructed in \S\ref{Sec:Extension} is motivated by mirror symmetry, which proposes that the (complexified stringy) moduli space of symplectic deformations of a Calabi--Yau variety should be isomorphic to the moduli space of complex deformations of its mirror \cite{morrison1993compactifications}. 

In the context of mirror symmetry for HDTV log Calabi--Yau pairs, obtained by a blow-up $(X,D) \to (X_{\Sigma},D_{\Sigma})$ from a toric variety, we obtain a natural family of symplectic structures by varying the volume of the exceptional divisors, hence a symplectic deformation space for $(X,D)$. In what follows we show that the complex deformation space of the extended mirror family corresponds to this symplectic deformation space of $(X,D)$. In a moment we describe a (formal) scheme $\hat{T}$, which we show is birational to $S_{(X_\Sigma, H)}$ and over which we have an extension of the mirror family to $(X,D)$, and for convenience we work with $\hat{T}$ in what follows.

First note that the cone of effective curves
$NE(X) \subset N_1(X)$ is actually contained in $NE(X_\Sigma) \oplus \Z^I$ because $C\in NE(X)$ implies $\mathrm{Bl}_{H,*}C \in NE(X_\Sigma)$.
More obviously, we also have 
\[ NE(X_\Sigma) \oplus \NN^I \subset NE(X_\Sigma) \oplus \Z^I\,.\]
Let $T$ be the scheme obtained by gluing together the two affine schemes 
\[ \Spec \mathbf{k}[NE(X)]\,\,\,\, \text{and} \,\,\,\,\Spec \mathbf{k}[NE(X_\Sigma) \oplus \NN^I]\]
along their common open subset $\Spec \mathbf{k}[NE(X_\Sigma) \oplus \Z^I]$.
By construction, 
$\Spec \mathbf{k}[\mathcal{M}]$ is the affinization of $T$ and we have a corresponding morphism 
\[ p \colon T \longrightarrow \Spec \mathbf{k}[\mathcal{M}]\,.\]

There are several interesting loci in $T$. 
First of all, the projections 
\[NE(X_\Sigma) \oplus \NN^I \rightarrow NE(X_\Sigma)\]
and \[NE(X_\Sigma) \oplus \NN^I \rightarrow \NN^I\]
induces closed embeddings 
\begin{equation}\nonumber
    \Spec \mathbf{k}[NE(X_\Sigma)] \hooklongrightarrow 
    \Spec \mathbf{k}[NE(X_\Sigma) \oplus \N^I] \hooklongrightarrow T
\end{equation}
and 
\begin{equation} \label{eq_subscheme_plus}
 \Spec \mathbf{k}[\N^I] \hooklongrightarrow 
    \Spec \mathbf{k}[NE(X_\Sigma) \oplus \N^I] \hooklongrightarrow T \,.  
\end{equation}
On the other hand, as the classes $E_i$ are effective for all 
$i\in I$, we have the inclusion 
\[ (-\NN)^I \simeq \bigoplus_{i\in I} \NN E_i \subset NE(X)\,,\] 
and so a closed embedding 
\begin{equation} \label{eq_subscheme_minus}
 \Spec \mathbf{k}[(-\N)^I] \hooklongrightarrow 
    \Spec \mathbf{k}[NE(X)] \hooklongrightarrow T  
\end{equation}
defined by the monoid ideal $\{ \beta \in NE(X) \,|\, \beta \notin (-\NN)^I\}$.
In $T$, the subschemes $\Spec \mathbf{k}[\NN^I]$ and 
$\Spec \mathbf{k}[(-\NN)^I]$ given by \eqref{eq_subscheme_plus} and 
\eqref{eq_subscheme_minus} are affine spaces with coordinates 
$(t^{-E_i})_{i\in I}$ and $(t^{E_i})_{i\in I}$ respectively, and are 
glued together along their common open subset 
$\Spec \mathbf{k}[\Z^I]$ via $t^{E_i} \mapsto t^{-E_i}$.
The resulting subscheme $T_0$ of $T$ is isomorphic to 
$(\PP^1)^I$ and is the fiber of 
$p$ over the torus fixed point of $\Spec \mathbf{k}[\mathcal{M}]$, defined by the monoid ideal 
$\mathcal{M} \setminus \{0\}$. 
Finally, let $\hat{T}$ be the formal completion of $T$ along $T_0 
\cup \Spec \mathbf{k}[NE(X_\Sigma)]$.
Then the restriction of $p$ to 
$\hat{T}$ defines a morphism 
\[ \hat{p} \colon \hat{T} \longrightarrow S_{(X_\Sigma,H)}\,,\] 
and one can consider the pullback $\check{\mathfrak{X}}_{\mathrm{ext}} \times_{S_{(X_{\Sigma},H)}}  \hat{T}$ of the mirror family to 
$\hat{T}$. From this point of view, the original mirror family was only defined on the formal neighborhood $\mathrm{Spf} \mathbf{k}[\![NE(X)]\!]$
of the point in 
$T_0$ with equations $t^{E_i}=0$ for all $i \in I$, whereas the extended mirror family is defined on a formal neighborhood of 
$T_0 
\cup \Spec \mathbf{k}[NE(X_\Sigma)]$. Moreover, by Theorem \ref{thm_toric_mirror}, the restriction of the extended mirror family to $\Spec \mathbf{k}[NE(X_\Sigma)]$ is the mirror family to the toric pair $(X_\Sigma,D_\Sigma)$.
We illustrate in Figure \ref{Fig:p} the map 
$p$ and the geometry of $T$ and  $\Spec \mathbf{k}[\mathcal{M}]$.

In the usual terminology of mirror symmetry, the point 
$t^{E_i}$=0 is the \emph{large complex structure limit} of the mirror family, corresponding to the \emph{large volume limit} of $X$. Indeed, the variable $t^{E_i}$ on the base of the space of complex deformations of the mirror should correspond by mirror symmetry to the function $e^{-\int_{E_i} \omega}$ on the space of classes of symplectic forms $\omega$ on $X$. In the large volume limit, we have
$t^{E_i} \rightarrow 0$, that is $\int_{E_i} \omega \rightarrow +\infty$, and the volume of the exceptional divisors with respect
to the symplectic form become indeed large.

The extension of the mirror family away from the large volume point $t^{E_i}=0$, until the point $t^{-E_i}=0$, correspond to 
moving away from the large volume point of $X$ in the space of symplectic forms, by decreasing the volume of the exceptional divisors. Such deformation is naturally provided by the birational map 
$\mathrm{Bl}_H \colon (X,D) \rightarrow (X_\Sigma, D_\Sigma)$:
symplectically, we have a continuous deformation from $(X,D)$ to 
$(X_\Sigma, D_\Sigma)$ obtained by continuously decreasing the volume of the exceptional divisors until they become of volume zero. 
While the classical symplectic geometry stops there, the symplectic form is complexified by the B-field in mirror symmetry, and the volume can be analytically continued until $-\infty$, that is, the point $t^{-E_i}=0$. 
To summarize, the extended mirror family gives a deformation of the mirror of $(X,D)$ to the mirror of $(X_\Sigma, D_\Sigma)$ which is mirror to the birational map 
$\mathrm{Bl}_H \colon (X,D) \rightarrow (X_\Sigma, D_\Sigma)$.

\begin{example}
\label{example_cones}
Let $(X,D)$ be the log Calabi-Yau surface obtained by blowing-up a point in 
$\PP^2$, as in Examples  \ref{Mirror to the blowup of P2 at a point} and \ref{example_M}. Recall that we denote by $L$ the pullback of the class of a line in $\PP^2$, by $E$ the exceptional divisor, and that we have
\[ NE(X_\Sigma)=\NN L\,,\,\,\, NE(X)=\NN E \oplus \NN (L-E)\,,\,\,\, \mathcal{M}= \NN L \oplus \NN (L-E)\,.\] 
As expected from Theorem \ref{thm_intro_extension_1}, the equation for the mirror to $(X,D)$ given in \eqref{Eq: mirror equation for BlP2} only involves 
$t^L$ and $t^{L-E}$, and so defines the extended mirror family over 
$\Spec \mathbf{k}[\mathcal{M}]$.
The dual cone of $NE(X)$ is the nef cone of $X$:
\[\mathrm{Nef}(X)=\NN L \oplus \NN (L-E)\,,\] whereas the dual of 
$NE(X_\Sigma) \oplus \NN$ is the ``bogus cone"
$\mathrm{Bog}(X)=\NN L \oplus \NN E$ in the terminology of 
\cite{hacking2020secondary}. The fan of $T$ is the union of the nef cone and of the bogus cone, and is in particular isomorphic to 
the blow-up of the affine plane at the origin. 
Erasing from this fan the middle ray $\NN L$, we obtain 
$\NN E \oplus \NN (L-E)$, which is the cone $\mathcal{M}^\vee$ dual to $\mathcal{M}$, and so the fan of
$\Spec \mathbf{k}[\mathcal{M}]$. In particular, $\Spec \mathbf{k}[\mathcal{M}]$ is isomorphic to an affine plane and 
$T \rightarrow \Spec \mathbf{k}[\mathcal{M}]$ is the blowup of the origin. The various cones of curves and divisors are illustrated 
in Figure \ref{Fig:cones}.
\end{example}

\section{Applications of the extended mirror family}
\label{Sec: applications}

\subsection{The Frobenius structure conjecture}

In this section, we use the extended mirror family to prove the Frobenius structure conjecture for HDTV log Calabi-Yau pairs.

We first review the statement of the Frobenius structure conjecture, proposed by Gross-Hacking-Keel as Conjecture 0.8 in the first arxiv version of \cite{GHK}.
Let $(X,D)$ be a log Calabi-Yau pair and 
$(B,\P)$ its tropicalization as in 
\S \ref{Sec: tropicalization}.
Let $A(X,D)$ be the topologically free 
$\mathbf{k}[\![NE(X)]\!]$-module over a set 
$\{ \vartheta_m\}_{m \in B(\Z)}$ indexed by the integral points of $B$:
\[ A(X,D):= \varprojlim_{k} \bigoplus_{m \in B(\Z)} 
\big( \mathbf{k}[NE(X)]/I_k\big)\, \vartheta_m  \,.\] 
We define the \emph{trace map} as being the projection on the coefficient of $\vartheta_0$:
\begin{align}\label{the trace map}
    \Tr : A(X,D) & \longrightarrow \mathbf{k}\lfor NE(X) \rfor \\
    \nonumber
    \sum_{m \in B(\Z)} a_m \vartheta_m & \longmapsto a_0 
\end{align}

We can now state the Frobenius structure conjecture:
\begin{conjecture}
\label{conj_frobenius}
Let $(X,D)$ be a log Calabi-Yau pair. Then, there exists a unique structure of 
$\mathbf{k}[\![NE(X)]\!]$-algebra on the $\mathbf{k}[\![NE(X)]\!]$-module $A(X,D)$ such that, for every $s\geq 2$ and $m_1, \dots, m_s \in B(\Z)$, 
\begin{equation} \label{eq_frobenius}
\Tr(\vartheta_{m_1} \cdots \vartheta_{m_s})
=\sum_{\beta \in NE(X)} N_\beta (m_1,\dots,m_s)
t^{\beta}\,, \end{equation}
where $N_\beta(m_1,\dots,m_s)$ is the log Gromov--Witten invariant of $(X,D)$
counting $(s+1)$-pointed rational curves 
$f: (C,x_0,x_1,\dots,x_s) \rightarrow X$, with contact orders $(0,m_1,\cdots,m_s)$ along $D$ at the marked points $(x_0,x_1,\cdots, x_s)$, and with insertion of $\psi_{x_0}^{s-2}$, where $\psi_{x_0}$
is the psi class attached to the marked point $x_0$.
\end{conjecture}

Recently, Johnston \cite[Theorem 1.4]{johnston2022comparison} proved that the algebra of theta functions $R(X,D)$, whose underlying $\mathbf{k}[\![NE(X)]\!]$-module is $A(X,D)$, satisfy \eqref{eq_frobenius}. This settles in particular the existence part of the conjecture (referred to as the weak Frobenius structure conjecture in 
\cite{johnston2022comparison}). The uniqueness part of the Conjecture \ref{conj_frobenius} has been proved by Keel-Yu \cite{KY} assuming that $X \setminus D$ is affine and contains a torus.
While the general case of the uniqueness part of the Conjecture \ref{conj_frobenius} is still open, we prove it below for HDTV log Calabi--Yau pairs $(X,D)$.

\begin{theorem} \label{thm_frob_unique}
Let $(X,D)$ be a HDTV log Calabi-Yau pair.
Then, the product structure of the algebra of theta functions 
$R(X,D)$ is uniquely determined by the values of the trace map 
\eqref{the trace map} on products of two and three theta functions.
\end{theorem}

\begin{proof}
We first introduce some notations.
For every $s \in \NN$, the \emph{s-point function} is the function
\begin{align}
    \Tr^s : R(X,D)^{\otimes s} & \longrightarrow \mathbf{k}\lfor NE(X) \rfor \\
    \nonumber
    x_1 \otimes \cdots \otimes x_s & \longmapsto \Tr (x_1 \cdots x_s) \,.
\end{align}
We similarly define the $s$-point function 
\[  \Tr^s_{\mathrm{ext}} : R(X,D)^{\otimes s }_{\mathrm{ext}}  \longrightarrow \widehat{\mathbf{k}[\mathcal{M}]}  \]
on the algebra $R(X,D)_{\mathrm{ext}}$ of functions on the extended mirror family given by Theorem \ref{thm_intro_extension_1}.
The structure constants $C_{m_1,m_2}^m$ defining the product of theta functions are the same in $R(X,D)$ and in $R(X,D)_{\mathrm{ext}}$: the only point of Theorem \ref{thm_intro_extension_1} is that the curve classes appearing in $C_{m_1,m_2}^m$, which are a priori in $NE(X)$, are actually contained in $\mathcal{M}$.
Therefore, for every $m_1, \cdots, m_s \in M$, we have 
\begin{equation}\label{eq_ext_trace}
    \Tr^s(\vartheta_{m_1},\cdots,\vartheta_{m_s})=\Tr^s_{\mathrm{ext}}(\vartheta_{m_1},\cdots,\vartheta_{m_s}) \,.
\end{equation}
In particular, it is enough to prove Theorem \ref{thm_frob_unique} for $R(X,D)_{\mathrm{ext}}$ to deduce it for $R(X,D)$: we have to show that the product structure on $R(X,D)_{\mathrm{ext}}$ can be recovered from the 2-point and 3-point functions 
$\Tr^2_{\mathrm{ext}}$ and $\Tr^3_{\mathrm{ext}}$.

We first prove that the trace pairing $\Tr^2 : R(X,D)_{\mathrm{ext}}^{\otimes 2}  \longrightarrow \widehat{\mathbf{k}[\mathcal{M}]} $ is non-degenerate, that is, the map
\begin{align}
    \label{nondegeneracy}
    R(X,D)_{\mathrm{ext}} & \longrightarrow \mathrm{Hom}(R(X,D)_{\mathrm{ext}} , \widehat{\mathbf{k}[\mathcal{M}]}) \\ \nonumber
x & \longmapsto (y \longmapsto \Tr^2_{\mathrm{ext}}(x,y))
\end{align}
is injective.

By Theorem \ref{thm_toric_mirror}, the restriction
of the extended mirror family to $\Spec \mathbf{k}[NE(X_\Sigma)]$ is the mirror family of the toric variety $(X_\Sigma, D_\Sigma)$.
It follows that for every $m \in M$, there exists $\kappa_m \in NE(X_\Sigma)$ such that, for every $n\in N$,
\[ \Tr^2_{\mathrm{ext}}(\vartheta_m, \vartheta_n)=t^{\kappa_m} \delta_{m+n,0} \mod J\,,\]
where $J$
is the ideal defining $\Spec \mathbf{k}[NE(X_\Sigma)]$ given by \eqref{eq_J}.

Now, to prove that the trace pairing is non-degenerate, it is enough to show that for every nonzero $x \in R(X,D)_{\mathrm{ext}}$, there exists $n \in M$ such that $\Tr^2 (x,\vartheta_n)\neq 0$. Let \[x=\sum_{m \in M} a_m \vartheta_m \in R(X,D)_{\mathrm{ext}} \setminus \{0\}\,,\] 
with $a_m \in \widehat{\mathbf{k}[\mathcal{M}]} $ for all $m \in M$. As $x \neq 0$, we can consider the smallest $k \geq 0$ such that 
$x \neq 0 \mod J^k$. Let $m_0 \in M$ be such that $a_{m_0} \neq 0 \mod J^k $. Then, 
\[ \Tr^2_{\mathrm{ext}}(x,\vartheta_{-m_0})=a_{m_0} t^{\kappa_m} \mod J^k\,,\]
and so in particular $\Tr^2_{\mathrm{ext}}(x,\vartheta_{-m_0}) \neq 0$.
This ends the proof that the trace pairing is non-degenerate.

We now conclude the proof of Theorem \ref{thm_frob_unique}
for $R(X,D)_{\mathrm{ext}}$ by showing that, for every $x_1,x_2 \in R(X,D)_{\mathrm{ext}}$ the product $x_1  x_2$ can be determined from the $2$-point and $3$-point functions. 
By the associativity of the ring of theta functions we can write $ \Tr^2_{\mathrm{ext}}(x_1 x_2, x_3)$ as a $3$-point function:
\[  \Tr^2_{\mathrm{ext}}(x_1 x_2, x_3) = \Tr^3_{\mathrm{ext}}( x_1,x_2,x_3) \, ,\]
Finally, by the non-degeneracy of the trace pairing on $R(X,D)_{\mathrm{ext}}$ the product $x_1 x_2$ is uniquely determined by the data of $ \Tr^2_{\mathrm{ext}}(x_1\ x_2, x_3)$ for all $x_3 \in R(X,D)_{\mathrm{ext}}$.
\end{proof}

As a corollary we obtain a proof of the Frobenius structure conjecture
for HDTV log Calabi-Yau pairs:

\begin{theorem} \label{thm_frobenius}
Let $(X,D)$ be a HDTV log Calabi-Yau pair. Then, the Frobenius structure conjecture (Conjecture \ref{conj_frobenius}) holds for $(X,D)$. More precisely, the algebra of theta functions $R(X,D)$ is the unique algebra satisfying \eqref{eq_frobenius}. 
\end{theorem}

\begin{proof}
By \cite[Theorem 1.4]{johnston2022comparison}, the algebra of theta functions $R(X,D)$ satisfies \eqref{eq_frobenius}.
In particular, the right-hand side of \eqref{eq_frobenius} is given by the trace of products of theta functions. By Theorem \ref{thm_frob_unique}, 
$R(X,D)$ is the unique algebra with the left-hand side 
of \eqref{eq_frobenius} equal to the trace of products of theta functions. Hence, it follows that $R(X,D)$ is the unique algebra satisfying \eqref{eq_frobenius}.
\end{proof}

\subsection{The HDTV mirror family}
\label{Sec: the HDTV mirror family}
We define in this section the HDTV mirror family to a HDTV log Calabi-Yau pair as a specific base change of the extended mirror family $\check{\mathfrak{X}}_{\mathrm{ext}} \to S_{(X_\Sigma,H)}$  -- see Definition \ref{def_intro_rest_mirror}. We then prove in Theorem \ref{thm_restriction} that the algebra of functions on the HDTV mirror family is isomorphic to the algebra of theta functions defined by the HDTV scattering diagram.

Let 
\[ \iota_0 \colon \mathrm{Spf} \,\mathbf{k}[\![\NN^I]\!] 
\hooklongrightarrow \mathrm{Spf}\, \mathbf{k}[NE(X_\Sigma)][\![\NN^I]\!] \]
be the morphism induced by restriction to the point 
\[ \Spec \mathbf{k} \hooklongrightarrow \Spec \mathbf{k}[NE(X_\Sigma)]  \]
defined by the equations $t^{\overline{\beta}}=1$ for all
$\overline{\beta} \in NE(X_\Sigma)$, i.e.\ the unit $1$ of the big torus orbit $\Spec \mathbf{k}[N_1(X_\Sigma)]$ in the toric variety 
$\Spec \mathbf{k}[NE(X_\Sigma)]$.
On the other hand, the inclusion
\[\mathcal{M} \subset NE(X_\Sigma)\oplus \NN^I\] induces a morphism 
\[ \mathrm{Spf} \,\mathbf{k}[NE(X_\Sigma)][\![\NN^I]\!]
\longrightarrow S_{(X_\Sigma,H)} \,.\]
Finally, we denote by 
\[ \iota \colon \mathrm{Spf}\, \mathbf{k}[\![\NN^I]\!] \longrightarrow S_{(X_\Sigma,H)} \]
the composition of this morphism with $\iota_0$.

\begin{definition} \label{def_intro_rest_mirror}
The \emph{HDTV} mirror family 
$\check{\mathfrak{X}}_{\mathrm{HDTV}}$ of $(X,D)$
with respect to the toric model $(X_\Sigma,H)$
is defined by the fiber diagram 
\[\begin{tikzcd}
\check{\mathfrak{X}}_{\mathrm{HDTV}} 
\arrow[r]
\arrow[ d]
&
\check{\mathfrak{X}}_{\mathrm{ext}}
\arrow[d]\\
\mathrm{Spf}\, \mathbf{k}[\![\NN^I]\!]
\arrow[r, "\iota"]& S_{(X_\Sigma,H)}\,.
\end{tikzcd}\]
\end{definition}

Restricting theta functions constructed from the canonical scattering diagram $\foD_{(X,D)}$ to $\check{\mathfrak{X}}_{\mathrm{HDTV}}$ we obtain
a topological basis $\{ \vartheta_m \}_{m \in M}$
of the algebra of functions on $\check{\mathfrak{X}}_{\mathrm{HDTV}}$ as a  
$\mathbf{k}[\![\NN^I]\!]$-module.
On the other hand, as reviewed in \S\ref{Sec:theta_functions}, one can construct using broken lines in the HDTV scattering diagram $\mathfrak{D}_{(X_\Sigma,H)}$ a topological $\mathbf{k}[\![\NN^I]\!]$-algebra $R(\mathfrak{D}_{(X_\Sigma,H)})$ 
with a basis of theta functions $\{ \vartheta_m^{HDTV} \}_{m \in M}$. In the following theorem, we compare the
theta functions $\vartheta_m$ and $\vartheta_m^{HDTV}$ using the main result of \cite{HDTV} comparing the canonical scattering diagram $\foD_{(X,D)}$ in $B$ and the HDTV canonical scattering diagram $\foD_{(X_\Sigma,H)}$ in $M_\RR$, along with a result of \cite{arguz2021heart} comparing the corresponding broken lines.

\begin{theorem} \label{thm_restriction}
Let $(X,D)$ be a HDTV log Calabi-Yau pair obtained as in 
\S \ref{Sec:HDTV_log_CY} from a toric model $(X_\Sigma, H)$.
Then, the algebra of theta functions $R(\mathfrak{D}_{(X_\Sigma,H)})$ defined by the HDTV scattering diagram $\mathfrak{D}_{(X_\Sigma,H)}$ is the algebra of functions of the HDTV mirror family, that is,
the map $\vartheta_m^{HDTV} \mapsto \vartheta_m$ for all $m \in M$ 
induces an isomorphism 
\[ \check{\mathfrak{X}}_{\mathrm{HDTV}} \simeq \mathrm{Spf}\, R(\mathfrak{D}_{(X_\Sigma,H)}) \,\]
over $\mathrm{Spf}\, \mathbf{k}[\![\NN^I]\!]$.
\end{theorem}

\begin{proof}
It is shown in \cite[Theorem 4.5]{arguz2021heart} that the mirror family $\check{\mathfrak{X}}$ to $(X,D)$ can be computed from a scattering diagram in $M_\RR$, called the heart of the canonical scattering diagram\footnote{More precisely, the heart of the canonical scattering diagram is not a scattering diagram as in Definition \ref{def_scattering_general}, but is a scattering diagram in the more general sense of \cite{GHS}. See \cite{arguz2021heart} for details.}. Briefly, to obtain the heart of the canonical scattering diagram, we first consider a degeneration to the normal cone of $(X,D)$ as described in \cite[\S3.1]{HDTV}. Considering the restriction of the canonical scattering diagram for the total space $(\tilde{X},\tilde{D})$ of this degeneration to the central fiber and then localizing around the origin we obtain a new scattering diagram, denoted by $T_0\foD^1(\tilde{X},\tilde{D})$, in \cite[Eq.5.1]{HDTV}. The group of curve classes of the central fiber of this degeneration is generated by the curve classes of $X_{\Sigma}$, the exceptional curves and some fiber classes of curves that arise as $\PP^1$ bundles. Setting these fiber classes to zero in $T_0\foD^1(\tilde{X},\tilde{D})$, we obtain the heart of the canonical scattering diagram -- for details see \cite[\S4]{arguz2021heart}. 

Now, the main point is that, by \cite[Theorem 6.2]{HDTV}, the heart of the canonical scattering diagram can be reconstructed from the HDTV scattering diagram $\foD_{(X_{\Sigma},H)}$ as follows.
For every incoming wall 
\begin{equation} \label{eq_incoming_hdtv}
(\fod,(1+t_i z^{v_i})^{w})
\end{equation} of $\foD_{(X_{\Sigma},H)}$ there is a corresponding wall
\begin{equation} \label{eq_incoming heart}
 (\fod,(1+t^{-E_i} z^{v_i})^{w}) \end{equation} 
of the heart. For every outgoing wall $(\fod, f_\fod)$, with $\fod \subset \sigma \in \Sigma$, and
\begin{equation} \label{eq_outgoing_hdtv}
f_\fod=  \sum c \prod_{i \in I} (t_i z^{v_i})^{a_i}\,,\end{equation}
there is a corresponding wall 
\begin{equation} \label{eq_outgoing_heart}
(\fod, \sum c \,t^{\beta_{\mathbf{A},\sigma}} \prod_{i \in I} z^{a_i v_i})\,,\end{equation}
of the heart, with $\beta_{\mathbf{A},\sigma}$ given by \eqref{Eq: curve class} and so in particular of the form
\[\beta_{\mathbf{A},\sigma} = \bar{\beta}_{\mathbf{A},\sigma} -\sum_{i \in I}a_i E_i \,,\]
with $\bar{\beta}_{\mathbf{A},\sigma} \in NE(X_\Sigma)$ and $a_i \geq 0$ for all $i \in I$.
In addition, the heart of the canonical scattering diagram
is defined using the PL function on $M_\RR$ whose kink across a codimension one cone $\rho$ of $\Sigma$ is given by the corresponding curve class $D_\rho$ in $NE(X_\Sigma)$.

We can now finish the proof of Theorem \ref{thm_restriction}.
By \cite[Theorem 4.5]{arguz2021heart}, the mirror family $\check{\mathfrak{X}}$ to $(X,D)$ is computed from the heart of the canonical scattering. Hence, the HDTV mirror family of Definition \ref{def_intro_rest_mirror} is obtained by setting to zero the curve classes of the heart of the canonical scattering diagrams contained in $NE(X_\Sigma)$.
The kinks of the PL function become equal to zero, the incoming walls \eqref{eq_incoming heart} do not change, and the outgoing walls \eqref{eq_outgoing_heart} become 
\[ (\fod, \sum c  \prod_{i \in I} (t^{-E_i}z^{v_i})^{a_i})\,.\]
In other words, writing $t_i=t^{-E_i}$, we recover exactly the HDTV scattering diagram $\foD_{(X_\Sigma,H)}$
whose walls are given by \eqref{eq_incoming_hdtv}-\eqref{eq_outgoing_hdtv}, and this concludes the proof that the HDTV mirror family is computed by the HDTV scattering diagram.
\end{proof}

\section{Cluster varieties and their HDTV mirrors}
\label{sec_C_scattering}
In this section we first briefly review some background on cluster varieties
in \S \ref{sec_cluster}. We then study their HDTV mirrors in
\S \ref{Sec: HDTV for cluster}.

\subsection{Cluster varieties}
\label{sec_cluster}

\subsubsection{$\mathcal{A}$ and $\mathcal{X}$ cluster varieties} \label{sec_intro_A_X}

Let $\mathbf{s}$ be a \emph{skew-symmetric seed}, given by the data of 
a rank $n$ lattice $N$, an integral skew-symmetric bilinear form 
\[\{ \cdot , \cdot  \} : N\times N  \longrightarrow \Z,\] 
an index set $\overline{I}$ of size $n$, a subset $I \subset \overline{I}$, and a basis
$(e_i)_{i \in \overline{I}}$ of $N$.
Basis elements $e_i$ are referred to as
\emph{unfrozen} if $i \in I$, and \emph{frozen} if 
$i \in \overline{I} \setminus I$. We denote by 
$M:= \Hom(N,\Z)$ the dual lattice of $N$, and by 
\[N_{\mathrm{uf}}=\Z^I=\bigoplus_{i\in I}\Z e_i\] the sublattice of $N$ spanned by the unfrozen basis elements.
We consider the map 
\begin{align}
 \nonumber
    p_1: N_{\mathrm{uf}} & \longrightarrow M \\
  \nonumber
   n & \longmapsto  \{ n , \cdot  \}
\end{align}
and for every $i\in I$ we denote 
$v_i :=p_1(e_i)=\{e_i,-\} \in M$. Finally, we assume that $v_i \neq 0$ for every 
$i \in I$.

From the data of a seed, there are two types of cluster varieties one can construct: the $\mathcal{X}$ and $\mathcal{A}$ cluster varieties.
The $\mathcal{X}$ (resp.\ $\mathcal{A}$) cluster variety 
is defined by gluing together copies of the torus 
$\Spec \mathbf{k}[N]$ (resp.\ $\Spec \mathbf{k}[M]$),
indexed by mutations of $\mathbf{s}$, using explicit birational maps referred to as cluster transformations, see \cite[\S 2]{GHKbirational} for details. However, throughout this paper we will use the interpretation of the $\mathcal{X}$ and $\mathcal{A}$ cluster varieties in terms of blow-ups of toric varieties, as shown in \cite{GHKbirational}. 
We review this below\footnote{In this paper we consider skew-symmetric seeds and cluster varieties constructed from such seeds, that is, skew-symmetric cluster varieties. However, we expect the results of this paper to hold in general for cluster varieties associated to skew-symmetrizable seeds. In that situation the blow-up construction of such varieties typically produce orbifolds/Deligne-Mumford (DM) stacks, rather than smooth algebraic varieties. Hence, treating such situations would require a generalization of the Gross--Siebert mirror construction \cite{HDTV,gross2021canonical} to DM stacks.}.

To construct the $\mathcal{X}$ cluster variety, we choose a smooth projective toric fan $\Sigma$ in $M_{\RR}:=M \otimes \RR$ whose set of rays contains $\RR_{\geq 0}v_i$ for all $i\in I$, and such that the hyperplanes 
$e_i^{\perp}$ are union of codimension one cones of $\Sigma$. 
We also assume that no cone of $\Sigma$ contains two distinct rays of the form $\RR_{\geq 0}v_i$.
Such fans always exist
and we denote by $X_{\Sigma}$ the corresponding toric variety.
For each $i \in I$, let $D_i$ be the toric boundary component corresponding to the ray $\RR_{\geq 0} v_i$ of $\Sigma$ and let $H_i \subset D_i$ be the hypersurface defined as the closure of the locus in $D_i$ of equation \[ (1+z^{e_i})^{|v_i|}=0\,, \]
where $|v_i|$ is the divisibility of $v_i$ in $M$. Note that $H_i$ is connected because 
$e_i$ is primitive in $N$.
Let $X$ be the blow-up of $X_\Sigma$ along  \[ H:=\bigcup_{i\in I}H_i\,,\] and let $D \subset X$ be the strict transform of the toric boundary divisor $D_\Sigma$ of $X_\Sigma$. 
Then, according to \cite[Theorem 3.9]{GHKbirational}, the complement $U_{\mathcal{X}}:= X \setminus D$ is isomorphic to the $\mathcal{X}$ cluster variety outside codimension two. We call $(X,D)$ a \emph{log Calabi-Yau compactification} of the $\mathcal{X}$ cluster variety.
In order to apply the mirror construction to $(X,D)$, one needs to be ensure that $X$ is smooth. This can fail for two reasons. If $|v_i|>1$, then blowing-up $H_i$ produces an $(n-2)$-dimensional family of 
$A_{|v_i|-1}$ surface quotient singularities. On the other hand, if 
$D_i=D_j$ when $i \neq j$, then $H_i$ and $H_j$ intersect and blowing-up 
$H_i \cup H_j$ also produces singularities.
So, when applying the mirror construction to $(X,D)$, we will make the following \emph{$\mathcal{X}$ assumptions} on the seed $\mathbf{s}$:
\begin{itemize}
    \item[(i)]for every $i \in I$, $v_i$ is primitive in $M$, that is, $|v_i|=1$, and 
    \item[(ii)]for every $i,j \in I$, 
$\RR_{\geq 0}v_i \neq \RR_{\geq 0} v_j$ if $i \neq j$.
\end{itemize}
Under these assumptions, $(X,D)$ is a HDTV log Calab-Yau pair as in 
\S \ref{Sec:HDTV_log_CY}.
We expect that the mirror construction and the main results of 
\cite{HDTV} should generalize to the orbifold case, thus removing assumption (i). We also expect the main results of \cite{HDTV} to generalize to the case of intersecting hypersurfaces contained in a common boundary divisors. Applying this generalization to successive blow-ups of $X_\Sigma$ along the hypersurfaces $H_i$ would allow us to remove assumption (ii).

Similarly, to construct the $\mathcal{A}$ cluster variety we choose a smooth projective toric fan $\Sigma'$ in $N_\RR := N \otimes \RR$ whose set of rays contains $\RR_{\geq 0}e_i$  for $i \in I$, and such that the hyperplanes $v_i^{\perp}$ are union of codimension one cones of $\Sigma'$.
We also assume that no cone of $\Sigma'$ contains two distinct rays of the form $\RR_{\geq 0}e_i$.
Such fans always exist
and we denote by $X_{\Sigma'}$ the corresponding toric variety.
For each $i\in I$, let $D_i'$ be the toric boundary component corresponding to the ray $\RR_{\geq 0}e_i$ of $\Sigma'$ and let $H_i' \subset D_i'$ be the hypersurface defined as the closure of the locus in 
$D_i'$ of equation 
\[ 1+z^{v_i}=0\,.\] 
Note that $H_i'$ consists of 
$|v_i|$ connected components, where $|v_i|$ is the divisibility of $v_i$
in $M$.
Let $X'$ be the blow-up of $X_{\Sigma'}$ along 
\[ H':= \bigcup_{i \in I} H_i'\] 
and let $D' \subset X'$ be the strict transform of the toric boundary divisor $D_{\Sigma'}$ of 
$X_{\Sigma'}$. 
By the assumptions on $\Sigma$, $X'$ is smooth and so 
$(X',D')$ is a HDTV log Calabi-Yau pair as in \S \ref{Sec:HDTV_log_CY}.
According to \cite[Theorem 3.9]{GHKbirational}, assuming that the seed 
$\mathbf{s}$ is totally coprime in the sense of 
\cite[Definition 3.7]{GHKbirational}, the complement 
$U_{\mathcal{A}}:= X' \setminus D'$ is isomorphic to the $\mathcal{A}$ cluster variety outside codimension two.
We call $(X',D')$ a \emph{log Calabi-Yau compactification} of the $\mathcal{A}$ cluster variety.

As reviewed in \S \ref{Sec:mirror family}, mirror symmetry naturally involves families of varieties. In particular, understanding mirror symmetry for cluster varieties requires working with families of cluster varieties. The relevant families are given by \emph{cluster varieties with principal coefficients}, which we review in the next section.

\subsubsection{Cluster varieties with principal coefficients}
The $\mathcal{A}$ and $\mathcal{X}$ cluster transformations are special cases of more general cluster transformations with principal coefficients
parametrized by the torus $\Spec \mathbf{k}[N]$.
We refer to \cite{GHKbirational, GHKK} for the $\mathcal{A}$ case, to 
\cite{Xcoeff} for the $\mathcal{X}$ case and to \cite[\S 3-4]{mou2021scattering}
for a uniform exposition.

Gluing families of tori $\Spec \mathbf{k}[M]$ parametrized by 
$\Spec \mathbf{k}[N]$ using the $\mathcal{A}$ cluster transformations with principal coefficients produces the 
\emph{$\mathcal{A}$ cluster variety with principal coefficients} 
\begin{equation} 
\nonumber
\pi_{\mathcal{A}}: \mathcal{A}_{\mathrm{prin}} \longrightarrow \Spec \mathbf{k}[N] \,.\end{equation}
The fiber over the unit in the torus $\Spec \mathbf{k}[N]$
recovers the $\mathcal{A}$ cluster variety:
$\pi_{\mathcal{A}}^{-1}(1)=\mathcal{A}$.
Moreover, allowing the coefficients $t_i:=t^{e_i}$ to vanish defines a partial compactification 
\begin{equation} \nonumber
\overline{\pi}_{\mathcal{A}} \colon \overline{\mathcal{A}}_{\mathrm{prin}} \longrightarrow \Spec \mathbf{k}[N^{\oplus}] \,,\end{equation} 
where $N^{\oplus}=\bigoplus_{i \in \overline{I}} \NN e_i$.
The special fiber over $0$ is the torus $\pi_{\mathcal{A}}^{-1}(0)=\Spec \mathbf{k}[M]$.
In fact, we will only need to consider the restriction
defined by setting $t_i=1$ for  all $i \in \overline{I}\setminus I$:
\begin{equation} \label{eq_compact_A_prin}\overline{\pi}_{\mathcal{A},\mathrm{uf}} \colon \overline{\mathcal{A}}_{\mathrm{prin}, \mathrm{uf}} \longrightarrow \Spec \mathbf{k}[N_{\mathrm{uf}}^{\oplus}] \,,\end{equation}where $N_{\mathrm{uf}}^{\oplus}=\bigoplus_{i\in I} \NN e_i$. 
The formal completion
\begin{equation} \label{eq_formal_A_prin}\widehat{\overline{\pi}}_{\mathcal{A},\mathrm{uf}} \colon \widehat{\overline{\mathcal{A}}}_{\mathrm{prin}, \mathrm{uf}}\longrightarrow \mathrm{Spf}\, \mathbf{k}[\![N_{\mathrm{uf}}^{\oplus}]\!]\end{equation}
of $\overline{\pi}_{\mathcal{A},\mathrm{uf}}$ along 
$\overline{\pi}_{\mathcal{A},\mathrm{uf}}^{-1}(0)$ will play an essential role in our study of mirror symmetry for cluster varieties.
One of the main result of Gross-Hacking-Keel-Kontsevich,
\cite[Proposition 6.4 (4)]{GHKK}, is the construction of  a topological basis $\{ \vartheta_m^C \}_{m \in M}$ of \emph{cluster theta functions} for the algebra of functions on 
$\widehat{\overline{\mathcal{A}}}_{\mathrm{prin}, \mathrm{uf}}$
as a $\mathbf{k}[\![N_{\mathrm{uf}}^{\oplus}]\!]$-module. 
The cluster theta functions are the theta functions
defined by a combinatorially constructed \emph{cluster scattering diagram} $\mathfrak{D}_{\mathbf{s}}^{\mathcal{A}_{\mathrm{prin}}}$, reviewed in \S \ref{Sec: cluster theta functions}.

In \cite{GHKbirational, GHKK}, $\mathcal{A}_{\mathrm{prin}}$
is described as the $\mathcal{A}$ cluster variety associated to the seed 
with principal coefficients $\tilde{\mathbf{s}}$ defined as follows.
\begin{definition}
\label{Def: seed with principal coefficients}
The \emph{seed with principal coefficients}, denoted by $\tilde{\mathbf{s}}$, is the data of the lattice
$\widetilde{N}:= N \oplus M$, endowed with 
the integral skew-symmetric bilinear form 
\begin{equation} \label{eq_A_prin_form}
\{(n_1,m_1),(n_2,m_2)\}:=\{n_1,n_2\}+\langle n_1,m_2 \rangle -\langle n_2,m_1\rangle\,,\end{equation}
where $\langle -,-\rangle$ is the duality pairing between $M$ and $N$, and the basis 
elements $((e_i,0))_{i\in \overline{I}}$
and $((0,e_i^{\star}))_{i\in \overline{I}}$, where 
$(e_i^{\star})_{i\in \overline{I}}$ is the basis of $M$ dual to 
$(e_i)_{i\in \overline{I}}$. Finally, the unfrozen basis elements are $((e_i,0))_{i\in I}$. 
\end{definition}

The cluster variety $\mathcal{A}_{\mathrm{prin}}$ has a very natural description in terms of blow-ups of toric varieties. Let $X_{\Sigma'}$ be a toric variety as in the description of the $\mathcal{A}$ cluster variety in \S \ref{sec_intro_A_X}, with the toric divisors $D_i'$. Then, 
$\mathcal{A}_{\mathrm{prin}}$ is isomorphic outside codimension two to the complement of the strict transform of $D' \times \Spec \mathbf{k}[N]$
in the variety obtained from $X' \times \Spec \mathbf{k}[N]$ by blowing-up the hypersurfaces in $D_i' \times \Spec \mathbf{k}[N]$ of equation 
\[ 1+t_i z^{v_i}=0 \,.\]
From this point of view, it is clear that one can allow the $t_i$'s to vanish and that the fiber over the point $\{t_i=0\}_i$ is the torus 
$\Spec \mathbf{k}[M]$. This description is compatible with viewing 
$\mathcal{A}_{\mathrm{prin}}$ as $\mathcal{A}$ for the seed $\tilde{s}$
because the unfrozen basis elements of $\tilde{s}$ are $(e_i,0)$ for $i\in I$, and we have $\tilde{v_i}:=\{(e_i,0),-\}=(v_i,e_i)$ and so 
$1+z^{\tilde{v_i}}=1+t_i z^{v_i}$.

The $\mathcal{X}$ case is similar.
Gluing families of tori $\Spec \mathbf{k}[N]$ parametrized by 
$\Spec \mathbf{k}[N]$ using the $\mathcal{X}$ cluster transformations with principal coefficients defined in \cite{Xcoeff} produces the 
\emph{$\mathcal{X}$ cluster variety with principal coefficients} 
\begin{equation} 
\label{eq:X_coeff}
\pi_{\mathcal{X}}: \mathcal{X}_{\mathrm{prin}} \longrightarrow \Spec \mathbf{k}[N] \,.\end{equation}
The fiber over the unit in the torus $\Spec \mathbf{k}[N]$
recovers the $\mathcal{X}$ cluster variety:
$\pi_{\mathcal{A}}^{-1}(1)=\mathcal{X}$.
Moreover, allowing the coefficients $t_i:=t^{e_i}$ to vanish defines a partial compactification 
\begin{equation} \nonumber
\overline{\pi}_{\mathcal{X}} \colon \overline{\mathcal{X}}_{\mathrm{prin}} \longrightarrow \Spec \mathbf{k}[N^{\oplus}] \,,\end{equation} 
whose special fiber over $0$ is the torus $\pi_{\mathcal{X}}^{-1}(0)=\Spec \mathbf{k}[N]$. We will mainly consider the restriction
defined by setting $t_i=1$ for  all $i \in \overline{I}\setminus I$:
\begin{equation} \label{eq_compact_X_prin}\overline{\pi}_{\mathcal{X},\mathrm{uf}} \colon \overline{\mathcal{X}}_{\mathrm{prin}, \mathrm{uf}} \longrightarrow \Spec \mathbf{k}[N_{\mathrm{uf}}^{\oplus}] \,,\end{equation}
and the formal completion
\begin{equation} \label{eq_formal_X_prin}\widehat{\overline{\pi}}_{\mathcal{X},\mathrm{uf}} \colon \widehat{\overline{\mathcal{X}}}_{\mathrm{prin}, \mathrm{uf}}\longrightarrow \mathrm{Spf}\, \mathbf{k}[\![N_{\mathrm{uf}}^{\oplus}]\!]\end{equation}
of $\overline{\pi}_{\mathcal{X},\mathrm{uf}}$ along 
$\overline{\pi}_{\mathcal{X},\mathrm{uf}}^{-1}(0)$.
An important difference with the $\mathcal{A}$ case is that
$\mathcal{X}_{\mathrm{prin}}$ is not the $\mathcal{X}$ cluster variety
associated to $\tilde{\mathbf{s}}$. Indeed, $\mathcal{X}_{\mathrm{prin}}$ is obtained by gluing tori \[\Spec \mathbf{k}[N] \times \Spec\mathbf{k}[N] =\Spec \mathbf{k}[N \oplus N]\,,\] whereas the  $\mathcal{X}$ cluster variety
associated to $\tilde{\mathbf{s}}$ is obtained by gluing tori \[\Spec \mathbf{k}[\widetilde{N}]=\Spec \mathbf{k}[N \oplus M]\,.\] In particular, $\mathcal{X}_{\mathrm{prin}}$ in \eqref{eq:X_coeff} is different from the variety denoted by $\mathcal{X}_{\mathrm{prin}}$ in
 \cite{GHKK}. As stressed in \cite{Xcoeff}, $\mathcal{X}_{\mathrm{prin}}$ as in \cite{GHKK} is the cluster dual to $\mathcal{A}_{\mathrm{prin}}$ 
 as cluster varieties over $\Spec \mathbf{k}$, but $\mathcal{X}_{\mathrm{prin}}$ as in \eqref{eq:X_coeff} is the cluster dual to $\mathcal{A}_{\mathrm{prin}}$ 
 as cluster varieties with coefficients, that is, over 
 $\Spec \mathbf{k}[N]$.

Finally, $\mathcal{X}_{\mathrm{prin}}$ has also a very natural description in terms of blow-ups of toric varieties. Let $X_{\Sigma}$ be a toric variety as in the description of the $\mathcal{X}$ cluster variety in \S \ref{sec_intro_A_X}, with the toric divisors $D_i$. Then,
by \cite[Proposition 5.14]{mou2021scattering},
$\mathcal{X}_{\mathrm{prin}}$ is isomorphic outside codimension two to the complement of the strict transform of $D \times \Spec \mathbf{k}[N]$
in the variety obtained from $X \times \Spec \mathbf{k}[N]$ by blowing-up the hypersurfaces in $D_i \times \Spec \mathbf{k}[N]$ of equation 
\[ (1+t_i z^{e_i})^{|v_i|}=0 \,.\]
By contrast, the blow-up description of  the $\mathcal{X}$ cluster variety for $\tilde{\mathbf{s}}$
would involve a fan in $M \oplus N$ with rays $\RR_{\geq 0}\tilde{v}_i=\RR_{\geq 0}(v_i,e_i)$ and hypersurfaces of equation $(1+z^{(e_i,0)})$.

\subsection{The HDTV mirror for cluster varieties}
\label{Sec: HDTV for cluster}
In this section, we describe the HDTV scattering diagram of \S \ref{Sec: HDTV scattering} for the HDTV log Calabi-Yau pairs obtained as log Calabi-Yau compactifications of $\mathcal{X}$ and 
$\mathcal{A}$ cluster varieties. 
By Theorem \ref{thm_restriction}, the algebra of functions of the HDTV mirror of such log Calabi-Yau pair is given by the algebra of theta functions defined by this scattering diagram. 

\subsubsection{The HDTV mirror for $\mathcal{X}$ cluster varieties}
\label{Sec:HDTV X}
Let $\mathbf{s}$ be a skew-symmetric seed satisfying the $\mathcal{X}$ assumptions and 
let $(X,D)$ be a log Calabi-Yau compactification of the corresponding
$\mathcal{X}$ cluster variety, with a toric model
$(X_\Sigma,H)$, as in \S \ref{sec_intro_A_X}. As reviewed in \S \ref{sec_intro_A_X}, $(X,D)$ is an example of HDTV log Calabi-Yau pair in the sense of 
\S \ref{Sec:HDTV_log_CY}. In this section, we describe the HDTV scattering diagram 
$\foD_{(X_\Sigma,H)}$ of \S \ref{Sec: HDTV scattering} for $(X,D)$.
 
Recall from \S\ref{sec_intro_A_X} that for every $i\in I$, we have a ray $\RR_{\geq 0} v_i$ of $\Sigma$, corresponding to a toric divisor $D_i$ of $X_\Sigma$, and the hypersurface $H_i \subset D_i$ defined as the closure of the locus of equation 
$1+z^{e_i}=0$ in $D_i$. In \S \ref{Sec: HDTV scattering}, the HDTV scattering diagram $\foD_{(X_\Sigma,H)}$
is described as a scattering diagram in $M_\RR$ over $R=\mathbf{k}[\![\NN^I]\!]$.
From now on, we identify $\NN^I$ with $N_{\mathrm{uf}}^{\oplus}=\bigoplus_{i \in I} \NN e_i$, that is $t_i$ with $t^{e_i}$, and we view 
$\foD_{(X_\Sigma,H)}$ as a scattering diagram in $M_\RR$ over 
$R=\mathbf{k}[\![N_{\mathrm{uf}}^{\oplus}]\!]$.

We first defined  the initial scattering diagram $\foD_{(X_\Sigma,H), \mathrm{in}}$ in \S \ref{Sec: HDTV scattering} in terms of tropical hypersurfaces $\scrH_i$ associated to the hypersurfaces $H_i$. The following result describes these tropical hypersurfaces in the $\mathcal{X}$ cluster case.

\begin{lemma}
\label{Lemma standard toric}
For every $i \in I$, let $\sigma$ be a codimension one cone of $\Sigma$ containing the ray $\RR_{\geq 0}v_i$. Then the weight $w_\sigma$ of the tropical hypersurface 
$\scrH_i$ on the cone $(\sigma+\RR v_i)/\RR v_i$
is given by
$w_\sigma=1$ if $\sigma \subset e_i^{\perp}$,
and $w_\sigma=0$ otherwise.
\end{lemma}

\begin{proof}
The quotient map $M \to M / \Z v_i$ induces an injective map
\begin{equation}
\label{Eq quotient}
\mathrm{Hom}(M / \Z v_i, \Z) \longrightarrow  \mathrm{Hom}(M, \Z) \, ,    
\end{equation}
whose image consists of the linear forms on $M$ which vanish on $v_i$. Since, \[\langle e_i, v_i \rangle =
\langle e_i , p_1(e_i) \rangle =\{ e_i,e_i\}= 0\,,\] 
there is a unique element $\overline{e}_i \in \mathrm{Hom}(M / \Z v_i, \Z)$, whose image under the map in \eqref{Eq quotient} is $e_i$.

By definition, the tropicalization of the equation $1+z^{e_i}$ restricted to $D_i$,
\[  1+z^{e_i}|_{D_i} =  1+z^{\overline{e}_i} \, , \]
corresponds to the PL function $\varphi_i$ on  $M_\RR/\RR v_i$ given by $\mathrm{max}\{ \overline{e}_i, 0 \}$, which restricts to linear functions on the complement of $\overline{e}_i^{\perp}$ and has kink $1$ along $\overline{e}_i^{\perp}$. By standard toric geometry, analogously as in the proof of \cite[Prop.6.2]{arguz2021heart}, for every codimension one cone $\sigma \subset \Sigma$ containing the ray $\RR_{\geq 0} v_i$, the intersection number $w_{\sigma} = D_{\sigma} \cdot H_i$ 
from Definition \ref{Def: tropical hypersurface associated}
equals to the kink of $\varphi_i$ along the cone $(\sigma+\RR v_i)/\RR v_i$. Hence, $w_{\sigma} = 1$ if  $(\sigma+\RR v_i)/\RR v_i \subset \overline{e}_i^{\perp}$, that is, $\sigma \subset e_i^{\perp}$; and $w_{\sigma} = 0$ otherwise.
\end{proof}

We can now describe the widgets 
$\foD_i$ of Definition \ref{Def HDTV scattering initial} defining the initial scattering diagram 
$\foD_{(X_\Sigma,H),\mathrm{in}}$. 

\begin{lemma} \label{Lem: widget computation}
For every $i \in I$, the \emph{widget} $\foD_i$ associated to the hypersurface $H_i$ is equivalent to the scattering diagram 
\begin{equation}
\nonumber
\foD_i^C:=\{(\fod_{\sigma}, f_i ) \,|\,\sigma ~\mathrm{is  ~ a ~maximal ~ dimensional ~ cone ~ in ~ } \overline{e}_i^{\perp} ~ \},
\end{equation}
where $f_i=1+t_i z^{v_i}$.
\end{lemma}

\begin{proof}
The result follows from the Definition \ref{def:widget_general} of the widget $\foD_i$ and from Lemma \ref{Lemma standard toric} computing the weights $w_\sigma$.
\end{proof}

For the comparison with the $\mathcal{A}_{\mathrm{prin}}$ cluster scattering diagram in \S \ref{sec_comparing_scattering}, it is more convenient to work with the consistent completion of the widgets. 

\begin{lemma} \label{Lem: widget completion}
For every $i \in I$, the consistent completion of 
$\foD_i^C$ is equivalent to
\[ S(\foD_i^C)=\{ (e_i^\perp, f_i) \}\]
where $f_i=1+t_iz^{v_i}$.
\end{lemma}

\begin{proof}
The hyperplane $e_i^{\perp}$ divides the real vector space $M_\RR$ into the two half-spaces $\{ e_i>0\}$ and $\{e_i<0\}$. In the scattering diagram
$\{ (e_i^\perp, f_i) \}$, the wall-crossing transformation for crossing from $\{ e_i>0\}$ to $\{e_i<0\}$ is the inverse of the wall-crossing transformation for crossing from $\{ e_i<0\}$ to $\{e_i>0\}$. As a general oriented loop in $M_\RR$ goes from $\{ e_i>0\}$ to $\{e_i<0\}$ as many times that it goes from $\{e_i<0\}$ to $\{ e_i>0\}$, it follows that the scattering diagram
$\{ (e_i^\perp, f_i) \}$ is consistent.

Moreover, as $\foD_i^C$ consists of all codimension one cones of $\Sigma$ contained in $e_i^{\perp}$ and containing $v_i$,  the scattering diagram  $\{ (e_i^\perp, f_i) \}$ is obtained from $\foD_i^C$ by adding walls which are all outgoing.  As the consistent completion of a scattering diagram is unique up to equivalence, it follows that $\{ (e_i^\perp, f_i) \}$ is equivalent to the consistent completion of $\foD_i^C$.

\end{proof}

We arrive at the following description of the HDTV scattering diagram in the $\mathcal{X}$ cluster case.

\begin{theorem} \label{thm:hdtv cluster}
Let $\mathbf{s}$ be a skew-symmetric seed satisfying the $\mathcal{X}$ assumptions and 
let $(X,D)$ be a log Calabi-Yau compactification of the corresponding
$\mathcal{X}$ cluster variety, with a toric model
$(X_\Sigma,H)$, as in \S \ref{sec_intro_A_X}. Then, the HDTV scattering diagram 
$\foD_{(X_\Sigma, H)}$ is equivalent to the consistent completion of the set of initial walls
\[\{ (e_i^\perp, f_i)\}_{i\in I} \,,\]
where $f_i=1+t_i z^{v_i}$.
\end{theorem}

\begin{proof}
By Definition \ref{Def HDTV scattering}, 
$\foD_{(X_\Sigma,H)}$ is the consistent completion of 
\[\foD_{(X_\Sigma,H) \mathrm{in}}=\bigcup_{i\in I} \foD_i\,.\]
By Lemma \ref{Lem: widget computation}, $\foD_i$ is equivalent to 
$\foD_i^C$, and so $\foD_{(X_\Sigma,H)}$ is equivalent to the consistent completion of 
$\cup_{i\in I} \foD_i^C$. 
Hence, by uniqueness of the consistent completion, $\foD_{(X_\Sigma,H)}$ is also equivalent to the consistent completion of 
$\cup_{i\in I} S(\foD_i^C)$. The result then follows from the explicit description of $S(\foD_i^C)$ in Lemma 
\ref{Lem: widget completion}.
\end{proof}

\subsubsection{The HDTV mirror for $\mathcal{A}$ cluster varieties}

Let $\mathbf{s}$ be a skew-symmetric seed and 
let $(X',D')$ be a log Calabi-Yau compactification of the corresponding
$\mathcal{A}$ cluster variety, with a toric model
$(X_{\Sigma'},H')$, as in \S \ref{sec_intro_A_X}. As reviewed in \S \ref{sec_intro_A_X}, $(X',D')$ is an example of HDTV log Calabi-Yau pair in the sense of 
\S \ref{Sec:HDTV_log_CY}. In this section, we describe the HDTV scattering diagram 
$\foD_{(X_{\Sigma'},H')}$ of \S \ref{Sec: HDTV scattering} for $(X',D')$.
 
Recall from \S\ref{sec_intro_A_X} that for every $i\in I$, we have a ray $\RR_{\geq 0} e_i$ of $\Sigma'$, corresponding to a toric divisor $D_i'$ of $X_\Sigma$, and the hypersurface $H_i' \subset D_i$ defined as the closure of the locus of equation 
$1+z^{v_i}=0$ in $D_i'$. Unlike what happens in the $\mathcal{X}$ case, $H_i'$ is not connected in general and consists of $|v_i|$ disjoint connected components $H_{i,j}'$ with $1 \leq j \leq |v_i|$. 
Let 
\begin{equation} \label{eq_tilde_I}
\widetilde{I}:=\{ (i,j)\,|\, i\in I \,,\,1 \leq j \leq |v_i| \} \,.\end{equation}
The HDTV scattering diagram $\foD_{(X_{\Sigma'},H')}$
is then a scattering diagram in $N_\RR$ over $R=\mathbf{k}[\![\NN^{\widetilde{I}}]\!]$.
We denote by $t_{ij}$ the monomial corresponding to the $(i,j)$ copy of $\NN$.

Using arguments similar to those use for the
$\mathcal{X}$ case in \S \ref{Sec:HDTV X}, we obtain the following description of the HDTV scattering diagram in the 
$\mathcal{A}$ cluster case.

\begin{theorem} \label{thm:hdtv cluster A}
Let $\mathbf{s}$ be a skew-symmetric seed and 
let $(X',D')$ be a log Calabi-Yau compactification of the corresponding
$\mathcal{A}$ cluster variety, with a toric model
$(X_{\Sigma'},H')$, as in \S \ref{sec_intro_A_X}. Then, the HDTV scattering diagram 
$\foD_{(X_{\Sigma'}, H')}$ is equivalent to the consistent completion of the set of initial walls
\[\{ (v_i^\perp, f_i)\}_{i\in I} \,,\]
where 
\begin{equation} f_i=\prod_{j=1}^{|v_i|}(1+t_{ij} z^{e_i})\,.\end{equation}
\end{theorem}

\section{Mirror symmetry and cluster dual varieties}
\label{Sec:mirror cluster}

In \S \ref{Sec: C scattering diagrams} we first review the description of the cluster scattering diagram of \cite{GHKK}, as a particular type of a $C$-scattering diagram. 
While the description of a $C$-scattering diagram we provide is similar to the one of a scattering diagram discussed in \S\ref{Sec: scattering diagrams}, there are some technical differences between the two notions. We give a precise comparison in \S\ref{sec_comparing_scattering}. Particularly, in \eqref{Eq: the map Psi}, we construct a map from a set of $C$-scattering diagrams to the set of scattering diagrams. 
Using this, we compare in Theorem \ref{Thm: HDTV and cluster} the
$\mathcal{A}_{\mathrm{prin}}$ cluster scattering diagram with the HDTV scattering diagram for cluster varieties described in \S \ref{Sec: HDTV for cluster}.
In Theorem \ref{thm_main} we prove the main result of this section, showing that the HDTV mirror to the $\mathcal{X}$ cluster variety is a deformation of the dual $\mathcal{A}$ cluster variety. Consequently, in \S\ref{sec: enumerative interpretation} we obtain enumerative interpretations of the structure constants in the algebra of theta functions $R(\foD_{\mathbf{s}}^{\mathcal{A}_{\mathrm{prin}}})$ defined by the $\mathcal{A}_{\mathrm{prin}}$ cluster scattering diagram.

\subsection{Cluster scattering diagrams}
\label{Sec: C scattering diagrams}
In this section we first introduce the general notion of a $C$-scattering diagram, then describe the cluster scattering diagram of \cite{GHKK} as a particular type of a $C_{\mathrm{prin}}$ scattering diagram.

\subsubsection{$C$-scattering diagrams}
\label{Sec: C scattering diagram}
Let $\mathbf{s}$ be a skew-symmetric seed, as in \S \ref{sec_intro_A_X}.
In this section and in \S \ref{Sec: cluster theta functions}, we make the following injectivity assumption on $\mathbf{s}$, as in \cite[pg 17]{GHKK}: we assume that the map
\begin{align}
\nonumber
    p_1: N_{\mathrm{uf}} & \longrightarrow M \\
  \nonumber
   n & \longmapsto  \{ n , \cdot  \}
\end{align}
is injective. 
Later, we will only apply this section and \S \ref{Sec: cluster theta functions}
to the case of the seed $\tilde{\mathbf{s}}$ with principal coefficients as in Definition \ref{Def: seed with principal coefficients}, 
for which the injectivity assumption is always satisfied.

In what follows we set 
\[ N^{\oplus} :=\bigoplus_{i \in I} \NN e_i \,, \]
and $N^+ := N^{\oplus} \setminus \{ 0 \}$. 
We define the monoid
\begin{equation}
    \nonumber
    P := p_1(N^{\oplus}) = \{ p_1(n) ~ | ~ n \in N^{\oplus}\} \,,
\end{equation}
and denote by $\mathbf{k}[P]$ the associated monoid algebra.
By the injectivity assumption, $0$ is the only invertible element of $P$, so $\mathfrak{m}_P:=P\setminus \{0\}$ is a maximal monoid ideal of 
$P$, and we denote by $\mathbf{k}[\![P]\!]:=\varprojlim_k \mathbf{k}[P]/\mathfrak{m}_P^k$
the completion of $\mathbf{k}[P]$ with respect to $\mathfrak{m}_P$.

\begin{definition}
\label{def_C_wall}
A \emph{$C$-wall} in $M_\RR$ for $N^+$ is a pair
$(\fod,f_{\fod})$, consisting of a codimension one rational polyhedral cone $\fod \subset n_\fod^{\perp} \subset M_{\RR}$ for some primitive $n_\fod \in N^+$, together with an attached function
$f_\fod \in \mathbf{k}[\![P]\!]$ of the form
\[ f_{\fod}=1+\sum_{k \geq 1} c_k z^{kp_1(n_\fod)}\,. \]
We say that a $C$-wall is \emph{incoming} if 
$p_1(n_\fod) \in \fod$, and \emph{outgoing} elsewise. We call $-p_1(n_\fod)$ the \emph{direction} of the $C$-wall.
\end{definition}

\begin{definition}
\label{def_C_scattering}
A \emph{$C$-scattering diagram} in $M_\RR$ is a set of $C$-walls $(\fod,f_{\fod})$, such that for every $k \geq 1$, there are only a finite number of $(\fod,f_{\fod}) \in \foD$ with $f_{\fod} \neq 1 \mod \mathfrak{m}_P^k$.
\end{definition}

\begin{remark}
It follows from \cite[Lemma 1.3]{GHKK} that the notion of a $C$-scattering diagram in Definition \ref{def_C_scattering} is equivalent to the notion of a scattering diagram introduced in \cite[\S 1.1]{GHKK} in the context of cluster algebras. We introduce the terminology of ``$C$-scattering diagram" in order to reserve ``scattering diagram" for the more general notion appearing in the context of mirror symmetry as in Definition \ref{def_scattering}.
\end{remark}

The \emph{support} of a $C$-scattering diagram $\foD$, denoted by $\mathrm{Supp}(\foD)$, is the union of all cones $\fod \subset M_{\RR}$ supporting a wall $(\fod,f_{\fod})$ of $\foD$. We define the \emph{singular locus} of $\foD$ 
\begin{equation}
\nonumber
\Sing(\foD)  :=  \bigcup_{\fod\in\foD} \partial\fod
\cup \bigcup_{\fod,\fod'\in\foD} (\fod\cap\fod'),
\end{equation}
where the last union is over all pairs of walls $\fod,\fod'$ with
$\fod\cap\fod'$ codimension at least two.

In what follows we review the notion of a consistent $C$-scattering diagram following \cite[\S 1.1]{GHKK}.
Every $C$-wall $(\fod,f_\fod)$ defines an automorphism $\mathfrak{p}_\fod$ of $\mathbf{k}[\![P]\!]$
given by 
\begin{align} \label{eq_p_d}
\mathfrak{p}_\fod \colon \mathbf{k}[\![P]\!] &\longrightarrow k[\![P]\!] \\ \nonumber
 z^m &\longmapsto f_\fod^{\langle n_\fod,m\rangle} z^m \,.\end{align}
Given a $C$-scattering diagram $\foD$ and a path \begin{align*} \gamma: [0,1] &\longrightarrow M_{\RR} \setminus \mathrm{Sing}(\foD)\\ 
 t &\longmapsto \gamma(t)\,\end{align*}
 transversal to the walls of $\foD$,
the associated path ordered product $\mathfrak{p}_{\gamma, \foD}$  is the automorphism of $\mathbf{k}[\![P]\!]$ obtained as the ordered product of the automorphisms 
$\mathfrak{p}_{\fod_i}^{\epsilon_{\gamma,\fod_i}}$ attached to the sequence of walls $\fod_i$ crossed by $\gamma$ for $t=t_i$, where 
\begin{equation} \label{eq_C_sign}
\epsilon_{\gamma,\fod_i}:= -\sgn (\langle n_{\fod_i}, \gamma'(t_i) \rangle )\in \{\pm 1\}\,.\end{equation}
Two $C$-scattering diagrams are \emph{equivalent} if they have the same path ordered products.
Finally, we call a scattering diagram \emph{consistent} 
if for any path $\gamma$ with $\gamma(0)=\gamma(1)$ 
the associated path ordered product is the identity automorphism.

According to \cite[\S 1.2]{GHKK}, for every 
$C$-scattering diagram $\foD_{\mathrm{in}}$
consisting of incoming $C$-walls $(\fod, f_\fod)$ in the sense of Definition \ref{def_C_wall} such that $\fod = n_\fod^{\perp}$, one can always find a consistent $C$-scattering diagram $S(\foD_{\mathrm{in}})$ containing $\foD_{\mathrm{in}}$, such that all $C$-walls in $S(\foD_{\mathrm{in}}) \setminus \foD_{\mathrm{in}}$ are outgoing. Moreover, $S(\foD_{\mathrm{in}})$ is unique up to equivalence with these properties and we refer to $S(\foD_{\mathrm{in}})$ as the \emph{consistent completion of
$\foD_{\mathrm{in}}$}.

\begin{remark}
\label{Rem: consistent C scattering}
From a consistent $C$-scattering diagram $\foD$, one can construct
using broken lines, analogously as in \S\ref{Sec: scattering diagrams and theta functions} -- see \cite[\S 3]{GHKK} for details-- a $\mathbf{k}[\![P]\!]$-algebra structure on the $\mathbf{k}[\![P]\!]$-module
\[ R(\foD):=  \varprojlim_k \bigoplus_{m \in M} \big(
\mathbf{k}[P]/\mathfrak{m}_P^k \big) \,\vartheta_m^C  \]
with basis elements $\vartheta_m^C$ indexed by $m \in M$.
We call $R(\foD)$ the \emph{algebra of theta functions} defined by the $C$-scattering diagram $\foD$, and we refer to the basis elements 
$\vartheta_m^C$ as \emph{theta functions}.
\end{remark}

\subsubsection{The cluster scattering diagram as a $C$-scattering diagram}
\label{Sec: cluster theta functions}

We review the description of the cluster scattering diagram following \cite[Theorem 1.12]{GHKK}, as a particular $C$- scattering diagram.

\begin{definition}
\label{Def cluster scattering}
Let $\foD_{\mathrm{in},\mathbf{s}}$ be the $C$-scattering diagram formed by the set of incoming $C$-walls given by 
\begin{equation} \label{eq_initial_cluster}
(e_i^{\perp}, 1+z^{v_i})
\end{equation}
for $i \in I$.
Then, the \emph{cluster scattering diagram}, denoted by 
\[\foD_{\mathbf{s}}:= S(\foD_{\mathrm{in},\mathbf{s}})  \, , \]
is the $C$-scattering diagram in $M_\RR$ obtained as the consistent completion of $\foD_{\mathrm{in},\mathbf{s}}$.
\end{definition}

\begin{remark}
Note that $v_i=p_1(e_i) \in e_i^{\perp}$ and so the 
$C$-walls $(e_i^{\perp}, 1+z^{v_i})$ are indeed incoming.
\end{remark}

We refer to the algebra of theta functions $R(\foD_{\mathbf{s}})$ defined by the cluster scattering diagram (see Remark \ref{Rem: consistent C scattering}) as the \emph{algebra of cluster theta functions}, and to the theta functions $\vartheta_m^C$ as the \emph{cluster theta functions}.

\subsubsection{$C_{\mathrm{prin}}$-scattering diagrams and the $\mathcal{A}_{\mathrm{prin}}$ cluster scattering diagram}
\label{Sec Cprin and Aprin}
In this section we introduce particular types of $C$-scattering diagrams in $\widetilde{M}_\RR := M_\RR \oplus N_\RR$, which we call $C_{\mathrm{prin}}$-scattering diagrams. Finally, we describe  $\mathcal{A}_{\mathrm{prin}}$ cluster scattering diagrams, which are important examples of both $C_{\mathrm{prin}}$ and cluster scattering diagrams. These are of particular interest as they can be used to reconstruct $\mathcal{A}_{\mathrm{prin}}$ cluster varieties \cite{GHKK}.

\begin{definition} \label{def_C_prin_scattering} 
Let $\mathfrak{D}=\{ (\fod, f_\fod)\}$ be a $C$-scattering diagram in 
$\widetilde{M}_\RR := M_\RR \oplus N_\RR$.
We say that $\foD$ is a \emph{$C_{\mathrm{prin}}$-scattering diagram}
in $\widetilde{M}_\RR$ if
\begin{itemize}
    \item[(i)] every wall $\fod$ of 
$\foD$ is contained in a hyperplane of the form $(n,0)^{\perp}$, where \[(n,0)\in \widetilde{N}=N \oplus M \,,\] 
and is invariant under translation by $N_\RR$, that is $\fod +N_\mathbb{R} \subset \fod$.  
    \item[(ii)] for every wall $\fod$ of $\foD$, we have $f_\fod \in \mathbf{k}[M][\![N_{\mathrm{uf}}^{\oplus}]\!]$. 
\end{itemize}
\end{definition}

We are now ready to define $\mathcal{A}_{\mathrm{prin}}$ cluster scattering diagrams.

\begin{definition}
\label{Def: DAprin}
The \emph{$\mathcal{A}_{\mathrm{prin}}$ cluster scattering diagram}, denoted by \[\foD^{\mathcal{A}_{\mathrm{prin}}}_{\mathbf{s}} :=\foD_{\tilde{\mathbf{s}}}   \,, \] 
is the cluster scattering diagram in $\widetilde{M}_\RR = M_\RR \oplus N_\RR$ associated to the seed with principal coefficients $\tilde{\mathbf{s}}$, defined in Definition \ref{Def: seed with principal coefficients}.
\end{definition}

\begin{lemma} \label{lem_initial_A_prin}
The initial walls of the  $\mathcal{A}_{\mathrm{prin}}$ cluster scattering diagram $\foD^{\mathcal{A}_{\mathrm{prin}}}_{\mathbf{s}}$ are given by 
\begin{equation} 
\nonumber
((e_i,0)^\perp, 1+z^{(v_i,e_i)}) \end{equation}
for all $i\in I$.
\end{lemma}

\begin{proof}
 Using the definition \eqref{eq_A_prin_form} of the skew-symmetric form for 
 $\tilde{\mathbf{s}}$, we have
\[ \{ (e_i,0),-\}=\{e_i,-\}+\langle e_i, -\rangle = (v_i,e_i) \in M \oplus N \,,\]
and so the result follows from the description of the initial walls of the cluster scattering diagram in \eqref{eq_initial_cluster}.
\end{proof}

\begin{lemma}
The $\mathcal{A}_{\mathrm{prin}}$ cluster scattering diagram $\foD^{\mathcal{A}_{\mathrm{prin}}}_{\mathbf{s}}$
is a $C_{\mathrm{prin}}$-scattering diagram
in $\widetilde{M}_\RR$.
\end{lemma}

\begin{proof}
Condition (i) of Definition \ref{def_C_prin_scattering} holds for the initial scattering diagram 
$\foD^{\mathcal{A}_{\mathrm{prin}}}_{\mathrm{in},\mathbf{s}}$ by Lemma 
\ref{lem_initial_A_prin}. Since this is preserved by the consistent completion, the condition holds for 
$\foD^{\mathcal{A}_{\mathrm{prin}}}_{\mathbf{s}}$.
On the other hand, the functions attached to the walls of 
$\foD^{\mathcal{A}_{\mathrm{prin}}}_{\mathbf{s}}$ are power series in the variables $z^{(v_i,e_i)}$ for $i \in I$. Indeed, it is the case for the initial walls by Lemma \ref{lem_initial_A_prin} and this property is preserved by the consistent completion. It follows that
condition (ii) of Definition \ref{def_C_prin_scattering}
holds for $\foD^{\mathcal{A}_{\mathrm{prin}}}_{\mathrm{in},\mathbf{s}}$.
\end{proof}

For every $C_{\mathrm{prin}}$-scattering diagram $\foD$, one can consider the algebra of theta functions $R(\foD)$ defined by 
$\foD$ as in Remark \ref{Rem: consistent C scattering}. As $\foD$ is a scattering diagram in 
$\widetilde{M}_\RR=M_\RR \oplus N_\RR$, the theta functions 
$\vartheta_{(m,n)}^C$
are indexed by $(m,n)\in \widetilde{M}=M \oplus N$.

\begin{lemma} \label{Lem:module}
Let $\foD$ be a consistent $C_{\mathrm{prin}}$-scattering diagram and let 
\[ \{ \vartheta_{(m,n)}^C \}_{(m,n) \in M \oplus N}\] 
be the corresponding set of theta functions.
Then, for every $(m,n) \in M \oplus N$, we have 
\[ \vartheta_{(m,n)}^C=z^{(0,n)} \vartheta_{(m,0)}^C\,.\] 
In particular, the topological span of the theta functions (see \cite[\S 2.2.2]{DavisonMandel}) of the form $\vartheta_{(m,n)}^C$ with $n \in N_{\mathrm{uf}}^{\oplus}$ has a natural structure of $\mathbf{k}[\![N_{\mathrm{uf}}^{\oplus}]\!]$-module, for which
$\{ \vartheta_{(m,0)}^C \}_{m \in M}$ is a 
topological $\mathbf{k}[\![N_{\mathrm{uf}}^{\oplus}]\!]$-basis.

Moreover, the span of the theta functions of the form $\vartheta_{(m,n)}^C$ with $n \in N_{\mathrm{uf}}^{\oplus}$
is a subalgebra of the algebra of theta functions, and so has a natural structure of $\mathbf{k}[\![N_{\mathrm{uf}}^{\oplus}]\!]$-algebra.
\end{lemma}

\begin{proof}
As the walls of $\mathfrak{D}$ are of the form $(n,0)^\perp$ and invariant by translation by $N$, a broken line of asymptotic direction $(0,n)$ never intersects any wall, so never bends. Hence, the first part of the lemma follows.

The second part follows from (ii) in Definition \ref{def_C_prin_scattering} of a $C_{\mathrm{prin}}$-scattering diagram and the formula \eqref{Eq: formula for structure constants}
computing product of theta functions in terms of broken lines.
\end{proof}

\subsection{From cluster scattering diagrams to scattering diagrams}

\label{sec_comparing_scattering}

Let $\foD=\{(\fod,f_\fod)\}$ be a $C_{\mathrm{prin}}$-scattering diagram in $\widetilde{M}_\RR=M_\RR \oplus N_{\RR}$,
as in Definition \ref{def_C_prin_scattering}. 
By Definition \ref{def_C_prin_scattering}(i), every wall $\fod$ of $\foD$ is invariant by translation by $N_\RR$ and so one can consider the quotient 
\[\fod/N_\RR \subset \widetilde{M}_\RR/N_\RR =M_\RR\,.\]

\begin{lemma} \label{lem_psi}
Let $\foD=\{(\fod,f_\fod)\}$ be a $C_{\mathrm{prin}}$-scattering diagram in $\widetilde{M}_\RR$.
Then 
\begin{equation}
\nonumber
    \Psi(\foD) := \{(\fod/N_{\RR}, f_\fod)\}
\end{equation}
is a scattering diagram in $M_\RR$ over $R=\mathbf{k}[\![N^{\oplus}_{\mathrm{uf}}]\!]$. 

Moreover, a $C$-wall 
$(\fod, f_\fod)$ of $\foD$ is incoming if and only if
the wall $(\fod/N_\RR, f_\fod)$ of $\Psi(\foD)$ is incoming.
\end{lemma}

\begin{proof}
Let $(\fod,f_\fod)$ be a C-wall of $\mathfrak{D}$. 
By Definition \ref{def_C_wall}, $\fod$ is a codimension one rational polyhedral cone in $\widetilde{M}_\RR$, contained in a hyperplane 
$n_{\fod}^\perp$ for some primitive $n_\fod \in \widetilde{N}^+$, and 
$f_\fod$ is of the form 
\[ f_\fod = 1+\sum_{k \geq 1} c_k z^{kp_1(n_\fod)}\,.\]
By Definition \ref{def_C_prin_scattering}, $n_{\fod}=(n,0)$ with $n\in N^+$, and $\fod$ is invariant under translation by $N_\RR$. 
In particular, $\fod/N_\RR$ is a codimension one rational polyhedral cone contained in $n^{\perp} \subset M_\RR$.
Using \eqref{eq_A_prin_form}, we find
\[ p_1(n_\fod)=p_1((n,0))=\{ n,-\}+\langle n,-\rangle 
=(p_1(n), n) \in M \oplus N \]
and so 
\[ f_\fod = 1+\sum_{k \geq 1} c_k t^{kn} z^{kp_1(n)}\,.\]
As $\langle p_1(n),n \rangle =\{n ,n\}=0$, it follows that 
 $-p_1(n) \in M$ is tangent to $\fod \subset n^{\perp}$.
Finally, we have $f_\fod \in \mathbf{k}[M][\![N_{\mathrm{uf}}^\oplus]\!]$ by Definition 
\ref{def_C_prin_scattering}(ii), and so
$(\fod/N_\RR, f_\fod)$ is a wall of direction positively proportional to $-p_1(n)$ as in Definition \ref{def_scattering}. Moreover, it follows that 
$(\fod,f_\fod)$ is incoming if and only if 
$(\fod/N_\RR, f_\fod)$ is incoming.

As $\foD$ is a $C$-scattering diagram as in Definition \ref{def_C_scattering}, it contains finitely many $C$-walls  $(\fod, f_\fod)$ with $f_\fod \neq 1 \mod \mathfrak{m}_P^k$. It follows
from $f_\fod \in \mathbf{k}[M][\![N_{\mathrm{uf}}^{\oplus}]\!]$
that $\Psi(\foD)$ also contains finitely many walls $(\fod/N_\RR, f_\fod)$ with $f_\fod \neq 1 \mod \mathfrak{m}_R^k$, where $\mathfrak{m}_R = N_{\mathrm{uf}}^{\oplus
}\setminus \{0\}$, and so $\Psi(\foD)$ is a scattering diagram over $R$ as in Definition \ref{def_scattering}.
\end{proof}

By Lemma \ref{lem_psi}, we have a well-defined map 
\begin{align} 
\label{Eq: the map Psi}
\Psi \colon C_{\mathrm{prin}}-\mathrm{Scatt} & \longrightarrow \mathrm{Scatt}\\
\nonumber \foD & \longmapsto \Psi(\foD) \end{align}
from the set $C_{\mathrm{prin}}-\mathrm{Scatt}$ of $C_{\mathrm{prin}}$-scattering diagrams in $\widetilde{M}_\RR$ up to equivalence to the set $\mathrm{Scatt}$ of scattering diagrams in $M_\RR$ over 
$R=\mathbf{k}[\![N_{\mathrm{uf}}^{\oplus}]\!]$ up to equivalence. This map $\Psi$ is injective: indeed, if $\Psi(\foD)=\{ (\fod, f_{\fod}) \}$, then we necessarily have
$\foD=\{ (\pi_M^{-1}(\fod), f_{\fod})\}$, where $\pi_M$ is the projection $\widetilde{M}_\RR=M_\RR \oplus N_\RR \rightarrow M_\RR$, since the walls of a $C_{\mathrm{prin}}$-scattering diagram are of the restricted form given in Definition \ref{def_C_prin_scattering}.

\begin{lemma} \label{lem_psi_cons}
Let $\foD$ be a $C_{\mathrm{prin}}$-scattering diagram in $\widetilde{M}_\RR$. Then, $\foD$ is consistent if and only if $\Psi(\foD)$ is consistent.
\end{lemma}

\begin{proof}
As $\mathfrak{D}$ is invariant by translation by $N$, it is enough to consider paths
$\gamma$ in $M_\RR$ transversal to the walls of $\Psi(\foD)$ to test the consistency of $\Psi(\foD)$, and their lifts $(\gamma,0)\in \widetilde{M}_\RR =M_\RR \otimes N_\RR$ to test the consistency of $\foD$.
Let $\fod$ be a wall of $\foD$. Then $n_{\fod}=(n,0)$ for some $n\in N^+$
and $n_{\fod/N_\RR} =n$.
Comparing \eqref{eq_p_d}-\eqref{eq_C_sign} with \eqref{eq_autom}, it is enough to show that 
\[ \epsilon_{\gamma,\fod} n =  n_{\fod/N_\RR, \gamma} \,.\]
This equality holds because by \eqref{eq_C_sign}, $\epsilon_{\gamma, \fod} (n,0)$ is a primitive normal vector to $\fod$ which is negative
on $(\gamma((t_0-\epsilon,t_0)),0)$ for all small $\epsilon>0$, and so
$\epsilon_{\gamma, \fod} n$ satisfies  the defining property of $n_{\fod/N_\RR, \gamma}$.
\end{proof}

\begin{lemma} \label{lem_psi_comm}
The map $\Psi$ commutes with the operation of taking consistent completions. That is, given a $C_{\mathrm{prin}}$-scattering diagram $\foD_{\mathrm{in}}$ in $\widetilde{M}_\RR$ with support a union of hyperplanes, we have
\[   S(\Psi (\foD_{\mathrm{in}}) ) = \Psi (  S(\foD_{\mathrm{in}}) ) \, .    \]

\end{lemma}

\begin{proof}
This follows from Lemma \ref{lem_psi} and Lemma \ref{lem_psi_cons}.
\end{proof}

Finally, we describe how $\Psi$ behaves with respect to 
theta functions.

\begin{lemma}
\label{Lemma: psi is nice}
Let $\foD$ be a consistent $C_{\mathrm{prin}}$-scattering diagram. Let $\{ \vartheta_{(m,n)}^C \}_{(m,n) \in M \oplus N}$ be the theta functions defined by $\mathfrak{D}$
and $\{ \vartheta_m\}_{m \in M}$ the theta functions defined by $\Psi(\mathfrak{D})$.
Then, the map 
\[ \vartheta_{(m,0)}^C \longmapsto \vartheta_m \]
induces an isomorphism of 
$\mathbf{k}[\![N_{\mathrm{uf}}^{\oplus}]\!]$-algebras between the algebra of theta functions of the form $\vartheta_{(m,n)}^C$ with $n\in N_{\mathrm{uf}}^{\oplus}$ and the algebra of theta functions $\{\vartheta_m\}_{m \in M}$. 
\end{lemma}

\begin{proof}
First note that we have indeed a structure of $\mathbf{k}[\![N_{\mathrm{uf}}^{\oplus}]\!]$-algebra on the algebra of theta functions of the form $\vartheta_{(m,n)}^C$ with $n \in N_{\mathrm{uf}}^{\oplus}$ by Lemma \ref{Lem:module}. It remains to compare the theta functions $\vartheta_{(m,0)}^C$ and $\vartheta_m$ for $m \in M$. 
As described in the proof of Lemma \ref{lem_psi},
walls of $\foD$ are of the form $\fod \subset (n,0)^\perp$ and $f_\fod=1+\sum_{k \geq 1} z^{k(p_1(n),n)}$ for some $n \in N$, and so 
all monomials of a broken line for $\foD$ of asymptotic direction $(m,0)$ are of the form $a_i z^{(m+p_1(n_i),n_i)}$ for some $n_i \in N$.
Let $\beta$ be a broken line for $\foD$ of asymptotic direction $(m,0)$ and ending at a point $x \in M_\RR \oplus N_\RR$. Let $\pi_M(\beta)$ be the image in $M_\RR$ of $\beta$ by the projection $\pi_M \colon M_\RR \oplus N_\RR \rightarrow M_\RR$. For every domain of linearity $\beta_i$ of $\beta$, with a monomial of the form $a_i z^{(m+p_1(n_i),n_i)}$ for some $n_i\in N$, we attach the monomial $a_i t^{n_i} z^{m+p_1(n_i)}$ to the domain of linearity $\pi_M(\beta_i)$ of $\pi_M(\beta)$. Then, it follows from the descriptions of walls of $\Psi(\foD)$ given in the proof of Lemma \ref{lem_psi} that $\pi_M(\beta)$ is a broken line for $\Psi(\foD)$ of asymptotic direction $m$ and endpoint $\pi_M(x)$. Moreover, the map $\beta \mapsto \pi_M(\beta)$ is a bijection between broken lines for $\foD$ of asymptotic direction $(m,0)$ and ending at $x$, and broken lines for $\Psi(\foD)$ of asymptotic direction $m$ and ending at $\pi_M(x)$, and so Lemma \ref{Lemma: psi is nice} follows.
\end{proof}

\subsection{Fock--Goncharov duals and HDTV mirrors to cluster varieties}
\label{Sec: mirrors and duals}

\subsubsection{Mirrors to $\mathcal{X}$ cluster varieties}

Let $\mathbf{s}$ be a skew-symmetric seed satisfying the $\mathcal{X}$ assumptions and 
let $(X,D)$ be a log Calabi-Yau compactification of the corresponding
$\mathcal{X}$ cluster variety, with a toric model
$(X_\Sigma,H)$, as in \S \ref{sec_intro_A_X}.
We first use the map $\Psi$ defined in 
\S \ref{sec_comparing_scattering} to compare the HDTV scattering diagram $\foD_{(X_\Sigma,H)}$
described in \S \ref{Sec:HDTV X} with the 
$\mathcal{A}_{\mathrm{prin}}$ cluster scattering diagram reviewed in \S \ref{Sec Cprin and Aprin}.

\begin{theorem}
\label{Thm: HDTV and cluster}
Let $\mathbf{s}$ be a skew-symmetric seed satisfying the $\mathcal{X}$ assumptions and 
let $(X,D)$ be a log Calabi-Yau compactification of the corresponding
$\mathcal{X}$ cluster variety, with a toric model
$(X_\Sigma,H)$, as in \S \ref{sec_intro_A_X}.
Then, the HDTV scattering diagram $\mathfrak{D}_{(X_\Sigma,H)}$
is the image of the $\mathcal{A}_{\mathrm{prin}}$ cluster scattering diagram by the map $\Psi$:
\[ \Psi(\foD^{\mathcal{A}_{\mathrm{prin}}}_{\mathbf{s}})
=\foD_{(X_\Sigma,H)}\,.\]
\end{theorem}

\begin{proof}
By Lemma \ref{lem_initial_A_prin}, the initial walls of
$\foD^{\mathcal{A}_{\mathrm{prin}}}$
are 
 $((e_i,0)^{\perp}, 1+z^{(v_i,e_i)})$
for $i \in I$.
On the other hand, from Theorem \ref{thm:hdtv cluster}, we know that $\mathfrak{D}_{(X_\Sigma,H)}$
is the consistent completion of the walls 
$(e_i^\perp, 1+t_i z^{v_i})$ for $i \in I$. 
Using the identifications 
$\NN^I=N_{\mathrm{uf}}^{\oplus}$, and $t_i=t^{e_i}$,
this can be rewritten as $ (e_i^\perp, 1+z^{(v_i,e_i)})$.
As $(e_i,0)^\perp/N_\RR=e_i^\perp$, we have 
\[ \Psi(\{ ((e_i,0)^{\perp}, 1+z^{(v_i,e_i)})\}_{i\in I})=\{(e_i^\perp, 1+z^{(v_i,e_i)}) \}_{i\in I} \,,\]
and then $\Psi(\foD^{\mathcal{A}_{\mathrm{prin}}}_{\mathbf{s}})
=\foD_{(X_\Sigma,H)}$ follows by Lemma \ref{lem_psi_comm}.
\end{proof}

As described in Definition \ref{def_intro_rest_mirror}, $(X,D)$ has a HDTV mirror, which is a family
\[ \check{\mathfrak{X}}_{\mathrm{HDTV}} \longrightarrow \mathrm{Spf}\, \mathbf{k}[\![\NN^I]\!]=
\mathrm{Spf}\, \mathbf{k}[\![N_{\mathrm{uf}}^{\oplus}]\!]\,.\]
The main result of this section below shows that the HDTV mirror family to $(X,D)$ is a degeneration of the dual $\mathcal{A}$ cluster variety.

\begin{theorem} \label{thm_main}
Let $\mathbf{s}$ be a skew-symmetric seed satisfying the $\mathcal{X}$ assumptions and 
let $(X,D)$ be a log Calabi-Yau compactification of the corresponding
$\mathcal{X}$ cluster variety as in \S \ref{sec_intro_A_X}.
Then, there exists a unique isomorphism over 
$\mathrm{Spf}\, \mathbf{k}[\![N_{\mathrm{uf}}^{\oplus}]\!]$ between 
the HDTV mirror family
\[ \check{\mathfrak{X}}_{\mathrm{HDTV}} \longrightarrow 
\mathrm{Spf}\, \mathbf{k}[\![N_{\mathrm{uf}}^{\oplus}]\!]\,\]
of $(X,D)$ and the formal partially compactified 
$\mathcal{A}$ cluster variety with principal coefficients as in \eqref{eq_formal_A_prin},
\[ \widehat{\overline{\pi}}_{\mathrm{uf}} \colon \widehat{\overline{\mathcal{A}}}_{\mathrm{prin}, \mathrm{uf}}
\longrightarrow \mathrm{Spf}\, \mathbf{k}[\![N_{\mathrm{uf}}^{\oplus}]\!]\,,\]
matching the basis of theta functions 
$\{ \vartheta_m \}_{m \in M}$ on 
$\check{\mathfrak{X}}_{\mathrm{HDTV}}$ with the basis of cluster theta functions $\{ \vartheta_{(m,0)}^C \}_{m\in M}$
on 
$\widehat{\overline{\mathcal{A}}}_{\mathrm{prin}, \mathrm{uf}}$.
\end{theorem}

\begin{proof}
To compare the cluster variety $\widehat{\overline{\mathcal{A}}}_{\mathrm{prin}, \mathrm{uf}}$ with the HDTV mirror family $\check{\mathfrak{X}}_{\mathrm{HDTV}}$ to $(X,D)$, first note that by Theorem \ref{thm_restriction} the algebra of functions of HDTV mirror $\check{\mathfrak{X}}_{\mathrm{HDTV}}$ is the algebra of theta functions defined by the scattering diagram $\foD_{(X_{\Sigma},H)}$. On the other hand, by
\cite[Proposition 6.4 (4)]{GHKK}, the algebra of functions of the cluster variety $\widehat{\overline{\mathcal{A}}}_{\mathrm{prin}, \mathrm{uf}}$ is
the algebra of theta functions of the form $\vartheta_{(m,n)}^C$ with $n \in N_{\mathrm{uf}}^{\oplus}$
constructed from the $\mathcal{A}_{\mathrm{prin}}$ cluster scattering diagram $\foD_{\mathbf{s}}^{\mathcal{A}_{\mathrm{prin}}}$ defined in Definition \ref{Def: DAprin}. The result then follows immediately from Lemma \ref{Lemma: psi is nice} and Theorem \ref{Thm: HDTV and cluster}.
\end{proof}

\begin{remark} \label{rem_X}
Let $(X,D)$ be a log Calabi-Yau compactification of a $\mathcal{X}$ cluster variety as in Theorem \ref{thm_main}. By Remark \ref{rem_affine_1}, when the complement $X \setminus D$ is affine, the mirror family canonically extends over $\Spec \mathbf{k}[NE(X)]$ and the extended mirror family extends further over $\Spec \mathbf{k}[\mathcal{M}]$. In this case, the HDTV mirror family extends over $\Spec \, \mathbf{k}[N_{\mathrm{uf}}^{\oplus}]$ and it follows from the proof of Theorem \ref{thm_main} that the resulting family is isomorphic to the partially compactified $\mathcal{A}_{\mathrm{prin}}$ cluster variety 
$\overline{\pi}_{\mathcal{A},\mathrm{uf}} \colon \overline{\mathcal{A}}_{\mathrm{prin}, \mathrm{uf}} \longrightarrow \Spec \mathbf{k}[N_{\mathrm{uf}}^{\oplus}]$
as in \eqref{eq_compact_A_prin}. In particular, the fiber over $1 \in \Spec \mathbf{k}[N_{\mathrm{uf}}]$ of the HDTV mirror family extended over 
$\Spec \mathbf{k}[N_{\mathrm{uf}}^{\oplus}]$, and so the fiber over 
$1 \in \Spec \mathbf{k}[N_1(X)]$ of the mirror family extended over 
$\Spec \mathbf{k}[NE(X)]$, are then isomorphic to the 
$\mathcal{A}$ cluster variety. Moreover, the restriction of the theta functions to this fiber are the cluster theta functions for the $\mathcal{A}$ cluster varieties defined in \cite{GHKK}. In the context of the non-archimedean mirror construction of \cite{KY}, this result is essentially \cite[Theorem 1.19]{KY}.
By contrast, Theorem \ref{thm_main}, stated in terms of formal families, holds in complete generality without the assumption that $X \setminus D$ is affine.
\end{remark}

\subsubsection{Mirrors to $\mathcal{A}$ cluster varieties}

Let $\mathbf{s}$ be a skew-symmetric seed and 
let $(X',D')$ be a log Calabi-Yau compactification of the corresponding
$\mathcal{A}$ cluster variety, with a toric model
$(X_{\Sigma'},H')$, as in \S \ref{sec_intro_A_X}.
As described in Definition \ref{def_intro_rest_mirror}, $(X',D')$ has a HDTV mirror, which is a family
\[ \check{\mathfrak{X}}_{\mathrm{HDTV}}' \longrightarrow \mathrm{Spf}\, \mathbf{k}[\![\NN^{\widetilde{I}}]\!]
\]
where $\widetilde{I}$ is given by \eqref{eq_tilde_I}.
We will consider the restriction
\begin{equation} \label{eq_HDTV_A_restricted}
\check{\mathfrak{X}}_{\mathrm{HDTV},I}' \longrightarrow \mathrm{Spf}\, \mathbf{k}[\![\NN^{I}]\!]=
\mathrm{Spf}\, \mathbf{k}[\![N_{\mathrm{uf}}^{\oplus}]\!]
\end{equation}
defined by setting $t_{ij}=t_i$ for all $(i, j)\in \widetilde{I}$.

The main result of this section, Theorem
\ref{thm_main_A}
below, shows that the HDTV mirror family to $(X',D')$, 
restricted as in \eqref{eq_HDTV_A_restricted},
is a degeneration of the dual $\mathcal{X}$ cluster variety.
The corresponding statement of Theorem \ref{thm_main} for the mirror to 
$\mathcal{X}$ cluster varieties involves a comparison with the cluster theta functions defined by \cite{GHKK} on $\mathcal{A}_{\mathrm{prin}}$. 
For the mirror to $\mathcal{A}$ cluster varieties, one would like a comparison with cluster theta functions on $\mathcal{X}_{\mathrm{prin}}$. However, $\mathcal{X}_{\mathrm{prin}}$ was not studied in \cite{GHKK} and so no cluster theta functions were constructed. Therefore, we have to explain how to adapt the techniques of \cite{GHKK} to construct cluster theta functions on $\mathcal{X}_{\mathrm{prin}}$. 

We first review how the $\mathcal{X}$ cluster variety is studied in \cite[\S 7]{GHKK}. The map of lattices
\begin{align*} 
\tilde{p}: N \oplus M &\longrightarrow M \\ (n,m)& \longmapsto m -p_1^{*}(n) 
 \end{align*}
induces a map 
\begin{equation} \label{eq_p_AX}
\tilde{p}: \mathcal{A}_{\mathrm{prin}} \longrightarrow \mathcal{X} 
\end{equation}
which realizes $\mathcal{X}$ as the quotient of 
$\mathcal{A}_{\mathrm{prin}}$ by the torus 
$\Spec \mathbf{k}[M]$. By duality, $\tilde{p}$ induces an inclusion 
\begin{align} \label{eq_inclusion}
\iota: N_\RR &\hooklongrightarrow \widetilde{M}_\RR= M_\RR \oplus N_\RR \\ \nonumber
n &\longmapsto (p_1^{*}(n),n)
\end{align}
As all monomials in the functions attached to the walls of 
$\foD_{\mathbf{s}}^{\mathcal{A}_{\mathrm{prin}}}$
are of the form $z^{(p_1^{*}n,n)}$, it makes sense to intersect $\foD_{\mathbf{s}}^{\mathcal{A}_{\mathrm{prin}}}$ with $N_\RR$ to obtain a $C$-scattering diagram $\foD_{\mathbf{s}}^{\mathcal{X}}$ in $N_\RR$ \cite[Construction 7.11]{GHKK}. By Lemma \ref{lem_initial_A_prin}, 
$\foD_{\mathbf{s}}^{\mathcal{A}_{\mathrm{prin}}}$ is the consistent completion of the initial walls 
$((e_i,0)^{\perp}, 1+z^{(v_i,e_i)})$ for $i \in I$. As $\tilde{p}((e_i,0))=-v_i$, it follows that $\foD_{\mathbf{s}}^{\mathcal{X}}$ is the consistent completion of the initial walls $(v_i^{\perp}, (1+z^{e_i})^{|v_i|})$\footnote{Whereas $(e_i,0)$ is the primitive normal vector to $(e_i,0)^{\perp}$, the primitive normal vector to $v_i^{\perp}$ is $v_i/|v_i|$, and so one has to include a power of $|v_i|$ in the attached function for the wall-crossing automorphism to stay the same.} for $i \in I$. Cluster theta functions for $\mathcal{X}$ are then defined using the $C$-scattering diagram $\foD_{\mathbf{s}}^{\mathcal{X}}$ , or equivalently as $\Spec \mathbf{k}[M]$-invariant cluster theta functions for 
$\mathcal{A}_{\mathrm{prin}}$.

To study $\mathcal{X}_{\mathrm{prin}}$, we consider an analogue of \eqref{eq_p_AX} with coefficients.
Let $\mathcal{A}_{\mathrm{prin},c}$ be the cluster variety over 
$\Spec \mathbf{k}[N]$ obtained from $\mathcal{A}_{\mathrm{prin}}$, viewed as a cluster variety over $\Spec \mathbf{k}$, 
by adding coefficients for the basis elements $(e_i,0)$, $i \in \overline{I}$, of the seed $\tilde{s}$. Then, $\tilde{p}$ extends to a map 
\[ \mathcal{A}_{\mathrm{prin},c} \longrightarrow \mathcal{X}_{\mathrm{prin}}\]
which realizes $\mathcal{X}_{\mathrm{prin}}$ as a quotient of $\mathcal{A}_{\mathrm{prin},c}$ by $\Spec \mathbf{k}[M]$.
Let $\mathfrak{D}_{\mathbf{s}}^{\mathcal{A}_{\mathrm{prin}, c}}$ be the scattering diagram in $\widetilde{M}_\RR$ over 
$\mathbf{k}[\![ N_{\mathrm{uf}}^{\oplus}]\!]$ obtained as the consistent completion of the initial walls $((e_i,0)^{\perp}, 1+t_i z^{(v_i,e_i)})$ for $i \in I$. As $\mathfrak{D}_{\mathbf{s}}^{\mathcal{A}_{\mathrm{prin}, c}}$
has the same support as 
$\mathfrak{D}_{\mathbf{s}}^{\mathcal{A}_{\mathrm{prin}}}$, the arguments of \cite{GHKK}, based on the existence of the cluster complex in the complement of the walls of $\mathfrak{D}_{\mathbf{s}}^{\mathcal{A}_{\mathrm{prin}}}$, generalize to $\mathfrak{D}_{\mathbf{s}}^{\mathcal{A}_{\mathrm{prin}, c}}$ and imply that the cluster theta functions $\{\vartheta_{(m,n)}^C\}_{(m,n)\in M \oplus N}$ defined by $\mathfrak{D}_{\mathbf{s}}^{\mathcal{A}_{\mathrm{prin}, c}}$
form a topological basis of the 
$\mathbf{k}[\![N^{\oplus}]\!]$-algebra of functions on the formal completion
\[ \widehat{\overline{\pi}}_{\mathcal{A},c}: \widehat{\overline{\mathcal{A}}}_{\mathrm{prin},c}
\longrightarrow \mathrm{Spf}\, \mathbf{k}[\![N^{\oplus}]\!]\,.\]
The intersection of $\mathfrak{D}_{\mathbf{s}}^{\mathcal{A}_{\mathrm{prin}, c}}$ with $N_\RR$
embedded in $\widetilde{M}_\RR$ as in \eqref{eq_inclusion}
is the scattering diagram $\mathfrak{D}_{\mathbf{s}}^{\mathcal{X}_{\mathrm{prin}}}$ in $N_\RR$ over $\mathbf{k}[\![N_{\mathrm{uf}}^{\oplus}]\!]$ obtained as the consistent completion of the initial walls $(v_i^{\perp}, (1+t_i z^{e_i})^{|v_i|})$ for $i \in I$.
We conclude that the cluster theta functions 
$\{ \vartheta_n^C\}_{n \in N}$ defined by 
$\mathfrak{D}_{\mathbf{s}}^{\mathcal{X}_{\mathrm{prin}}}$ form a topological basis of the 
$\mathbf{k}[\![N_{\mathrm{uf}^\oplus}]\!]$-algebra of functions on the formal completion 
\[ \widehat{\overline{\pi}}_{\mathcal{X},\mathrm{uf}} \colon \widehat{\overline{\mathcal{X}}}_{\mathrm{prin}, \mathrm{uf}}
\longrightarrow \mathrm{Spf}\, \mathbf{k}[\![N_{\mathrm{uf}}^{\oplus}]\!]\,,\]
as in \eqref{eq_formal_X_prin}.

\begin{theorem} \label{thm_main_A}
Let $\mathbf{s}$ be a skew-symmetric seed and 
let $(X',D')$ be a log Calabi-Yau compactification of the corresponding
$\mathcal{A}$ cluster variety as in \S \ref{sec_intro_A_X}.
Then, there exists a unique isomorphism over 
$\mathrm{Spf}\, \mathbf{k}[\![N_{\mathrm{uf}}^{\oplus}]\!]$ between 
the HDTV mirror family of $(X',D')$, restricted as in 
\eqref{eq_HDTV_A_restricted}, 
\[ \check{\mathfrak{X}}_{\mathrm{HDTV},I}' \longrightarrow \mathrm{Spf}\, \mathbf{k}[\![N_{\mathrm{uf}}^{\oplus}]\!]\,,
\]
and the formal partially compactified 
$\mathcal{X}$ cluster variety with principal coefficients as in \eqref{eq_formal_X_prin},
\[ \widehat{\overline{\pi}}_{\mathcal{X},\mathrm{uf}} \colon \widehat{\overline{\mathcal{X}}}_{\mathrm{prin}, \mathrm{uf}}
\longrightarrow \mathrm{Spf}\, \mathbf{k}[\![N_{\mathrm{uf}}^{\oplus}]\!]\,,\]
matching the basis of theta functions 
$\{ \vartheta_n \}_{n \in N}$ on 
$\check{\mathfrak{X}}_{\mathrm{HDTV}}'$ with the basis of cluster theta functions 
$\{ \vartheta_n^C \}_{n\in N}$
on 
$\widehat{\overline{\mathcal{X}}}_{\mathrm{prin}, \mathrm{uf}}$.
\end{theorem}

\begin{proof}
By Theorem \ref{thm_restriction} the algebra of functions of the HDTV mirror $\check{\mathfrak{X}}_{\mathrm{HDTV}}'$ is the algebra of theta functions defined by the HDTV scattering diagram $\foD_{(X_{\Sigma'},H')}$, which is given for $\mathcal{A}$ cluster varieties by Theorem 
\ref{thm:hdtv cluster A}. Setting $t_{ij}=t_i$ in Theorem
\ref{thm:hdtv cluster A}, we obtain the scattering diagram $\mathfrak{D}_{\mathbf{s}}^{\mathcal{X}_{\mathrm{prin}}}$ producing the cluster theta functions on 
$\widehat{\overline{\mathcal{X}}}_{\mathrm{prin}, \mathrm{uf}}$ and this concludes the proof of the theorem.
\end{proof}

\begin{remark} \label{rem_A}
Let $(X',D')$ be a log Calabi-Yau compactification of a $\mathcal{A}$ cluster variety as in Theorem \ref{thm_main_A}. As in Remark \ref{rem_X}, when the complement $X' \setminus D'$ is affine, the HDTV mirror family extends over $\Spec \, \mathbf{k}[N_{\mathrm{uf}}^{\oplus}]$ and  the resulting family is isomorphic to the partially compactified $\mathcal{X}_{\mathrm{prin}}$ cluster variety 
$\overline{\pi}_{\mathcal{X},\mathrm{uf}} \colon \overline{\mathcal{X}}_{\mathrm{prin}, \mathrm{uf}} \longrightarrow \Spec \mathbf{k}[N_{\mathrm{uf}}^{\oplus}]$
as in \eqref{eq_compact_X_prin}. In particular, the fiber over $1 \in \Spec \mathbf{k}[N_{\mathrm{uf}}]$ of the HDTV mirror family is isomorphic to the 
$\mathcal{X}$ cluster variety. Moreover, the restriction of the theta functions to this fiber are the cluster theta functions for the $\mathcal{X}$ cluster varieties defined in \cite{GHKK}. 
\end{remark}

\subsubsection{Mirrors to $\mathcal{X}$ symplectic fibers and $\mathcal{A}$ torus quotients}

Let $\mathbf{s}$ be a skew-symmetric seed. We consider the map 
\begin{align*}
    p :N &\longrightarrow M\\
    n &\longmapsto \{n,-\}
\end{align*}
and we denote by $K \subset N$ the kernel of $p$. The inclusion $K \subset N$
induces a  map of tori \[\nu: \Spec \mathbf{k}[N] \longrightarrow \Spec \mathbf{k}[K]\,,\] and then a map \[\lambda \colon \mathcal{X} \longrightarrow \Spec \mathbf{k}[K]\,.\] 
The skew-symmetric form $\{ -,-\}$ on $N$ defines a Poisson structure on $\mathcal{X}$ and the fibers of $\lambda$ are the corresponding symplectic leaves. Dually, denoting $K^\vee :=\Hom(K,\Z)$, the projection 
$M \rightarrow K^\vee$ induces maps of tori $\Spec \mathbf{k}[K^\vee] \rightarrow \Spec \mathbf{k}[M]$, and then an action of the torus $\Spec \mathbf{k}[K^\vee]$ on $\mathcal{A}$. To state a version of our results for the $\mathcal{X}$ symplectic fibers 
$\lambda^{-1}(t)$ and the quotient $\mathcal{A}/\Spec \mathbf{k}[K^\vee]$, we remark that the map $\lambda$ extends to a map 
$\widehat{\overline{\lambda}}_{\mathrm{prin},\mathrm{uf}}: 
\widehat{\overline{\mathcal{X}}}_{\mathrm{prin},\mathrm{uf}} \rightarrow 
\Spec \mathbf{k}[K]$, and that the action of $\Spec \mathbf{k}[K^\vee]$ on $\mathcal{A}$ extends to an action on 
$\widehat{\overline{\mathcal{A}}}_{\mathrm{prin},\mathrm{uf}}$.

\begin{theorem} \label{thm_mirror_X_fiber}
Let $\mathbf{s}$ be a skew-symmetric seed satisfying the 
$\mathcal{X}$ assumptions and let $(X,D)$ be a log Calabi-Yau compactification of a general $\mathcal{X}$ symplectic fiber $\lambda^{-1}(t)$. 
Then the HDTV mirror family 
\[\check{\mathfrak{X}}_{HDTV}\longrightarrow \mathrm{Spf}\, \mathbf{k}[\![N_{\mathrm{uf}}^{\oplus}]\!]\]
of $(X,D)$ is isomorphic to the family of quotients 
\[\widehat{\overline{\mathcal{A}}}_{\mathrm{prin},\mathrm{uf}}/
\Spec \mathbf{k}[K^\vee] \longrightarrow \mathrm{Spf}\, \mathbf{k}[\![N_{\mathrm{uf}}^{\oplus}]\!]\,.\]
\end{theorem}

\begin{proof}
Let $\mu \colon M_\RR \rightarrow K^\vee_\RR$ be the projection dual to the 
inclusion $K \subset N$. By Theorem \ref{thm:hdtv cluster}, the HDTV scattering diagram of a log Calabi-Yau compactification of the $\mathcal{X}$ cluster variety is the consistent scattering diagram in $M_\RR$ with initial walls 
$(e_i^{\perp}, 1+t_i z^{v_i})_{i\in I}$. By Theorem \ref{thm_main}, the corresponding theta functions $\{ \vartheta_m \}_{m \in M}$ are the cluster theta functions for $\widehat{\overline{\mathcal{A}}}_{\mathrm{prin},\mathrm{uf}}$. 
Since $\{e_i,K\}=0$, it follows that $\mu(v_i)=0$ for all $i\in I$, and so $\mu(v)=0$ for every direction $v$ of a wall. Hence, the intersection of this scattering diagram with $\mu^{-1}(0)$ is naturally a scattering diagram in $\mu^{-1}(0)$, which is in fact the HDTV scattering diagram $\foD_{(X_\Sigma,H)}$
of a log Calabi-Yau compactification $(X,D)$ of a $\mathcal{X}$ symplectic fiber 
$\lambda^{-1}(t)$ for $t \in \Spec \mathbf{k}[K]$. 
Indeed, the toric model $(X_\Sigma, H)$ of $(X,D)$ has the following description.  The toric variety $X_\Sigma$ is a toric compactification of the torus $\nu^{-1}(t) \simeq \Spec \mathbf{k}[N/K]$ fiber of $\nu: \Spec \mathbf{k}[N] 
\rightarrow \Spec \mathbf{k}[K]$, with fan in $\mu^{-1}(0)$ containing the rays $\RR_{\geq 0}v_i$. Moreover, the hypersurfaces $H_i$ have equations given by the restriction to the divisors $D_i$ of the rational function on $X_\Sigma$ obtained by restricting $1+z^{e_i}$ from $\Spec \mathbf{k}[N]$ to $\nu^{-1}(t)$.

Therefore, the theta functions for the mirror of $(X,D)$ are 
$\{ \vartheta_m \}_{m \in M, \mu(m)=0}$. The result follows because for every $m \in M$, $\mu(m)\in K^\vee$ is the weight of the action of 
$\Spec \mathbf{k}[K^\vee]$ on $\vartheta_m$, and theta functions on the quotient
$\widehat{\overline{\mathcal{A}}}_{\mathrm{prin},\mathrm{uf}}/\Spec \mathbf{k}[K^\vee]$ are exactly the weight zero theta functions.
\end{proof}

\begin{remark}
    The conclusion of Theorem \ref{thm_mirror_X_fiber} is still valid if one does not assume part (ii) of the $\mathcal{X}$ assumptions of  \S\ref{sec_intro_A_X}. Indeed, for general $t$, the fibers over $t$ of the various hypersurfaces $H_i$ do not intersect, and so a log Calabi-Yau compactification of the $\mathcal{X}$ symplectic fiber $\lambda^{-1}(t)$ is automatically an HDTV log Calabi-Yau pair.
\end{remark}

\begin{theorem} \label{thm_mirror_A_quotient}
Let $\mathbf{s}$ be a skew-symmetric seed and let 
$(X,D)$ be a log Calabi-Yau compactification of the quotient 
$\mathcal{A}/\Spec \mathbf{k}[K^\vee]$. Then the restricted HDTV mirror family 
\[\check{\mathfrak{X}}_{HDTV,I}\longrightarrow \mathrm{Spf}\, \mathbf{k}[\![N_{\mathrm{uf}}^{\oplus}]\!]\]
of $(X,D)$ is isomorphic to the family of 
$\mathcal{X}$ symplectic fibers 
\[\widehat{\overline{\lambda}}_{\mathrm{prin},\mathrm{uf}}^{-1}(1) \longrightarrow \mathrm{Spf}\, \mathbf{k}[\![N_{\mathrm{uf}}^{\oplus}]\!]\,.\]
\end{theorem}

\begin{proof}
By Theorem \ref{thm:hdtv cluster A}, the HDTV scattering diagram of a log Calabi-Yau compactification of the $\mathcal{A}$ cluster variety after setting 
$t_{ij}=t_i$ is the consistent scattering diagram in $N_\RR$ with initial walls 
$(v_i^{\perp}, (1+t_i z^{e_i}))^{|v_i|})_{i\in I}$. By Theorem \ref{thm_main_A}, the corresponding theta functions $\{ \vartheta_n\}_{n \in N}$ are the cluster theta functions for $\widehat{\overline{\mathcal{X}}}_{\mathrm{prin},\mathrm{uf}}$. As $\{e_i,K\}=0$, all the walls are invariant under the action of 
$K_\RR$ by translation on $N_\RR$. 
Hence, this scattering diagram naturally defines a scattering diagram in the quotient $(N/K)_\RR$, which is in fact the HDTV scattering diagram $\foD_{(X_\Sigma,H)}$
of a log Calabi-Yau compactification $(X,D)$ of the quotient 
$\mathcal{A}/\Spec \mathbf{k}[K^\vee]$. 
Indeed, the toric model $(X_\Sigma, H)$ of $(X,D)$ has the following description.  
The toric variety $X_\Sigma$ is a toric compactification of the torus $\Spec \mathbf{k}[M/K^\vee]=(\Spec \mathbf{k}[M])/(\Spec \mathbf{k}[K^\vee])$, with fan in $(N/K)_\RR$ containing the images of the rays $\RR_{\geq 0}e_i$ by the projection $N \rightarrow N/K$. 
Moreover, as $v_i(K)=\{e_i,K\}=0$, the functions $1+z^{v_i}$ on $\Spec \mathbf{k}[M/K]$ are invariant under the action of $\Spec \mathbf{k}[K^\vee]$, so induce functions on the quotient \[\Spec \mathbf{k}[M/K^\vee]=(\Spec \mathbf{k}[M])/(\Spec \mathbf{k}[K^\vee])\,,\] 
whose restrictions to the divisors $D_i$ are the equations of the hypersurfaces $H_i$.

Therefore, the theta functions for the mirror of $(X,D)$ are obtained by setting $\vartheta_n=1$ for every $n \in K$, 
and so the mirror is the family of $\mathcal{X}$ symplectic fibers $\widehat{\overline{\lambda}}_{\mathrm{prin},\mathrm{uf}}^{-1}(1)$.
\end{proof}

\begin{remark}
The map $p \colon N \rightarrow M$ induces first a map 
$\mathcal{A} \rightarrow \mathcal{X}$, and then a generically finite map
$\mathcal{A}/\Spec \mathbf{k}[K^\vee] \rightarrow \lambda^{-1}(1)$,
which is an isomorphism when the skew-symmetric form induced by 
$\{-,-\}$ on $N/K$ is unimodular. In this case, Theorems 
\ref{thm_mirror_X_fiber} and \ref{thm_mirror_A_quotient} imply that the 
holomorphic symplectic variety $\mathcal{A}/\Spec \mathbf{k}[K^\vee] \simeq \lambda^{-1}(1)$ is essentially self-mirror. 
\end{remark}

\subsection{Enumerative geometry of cluster varieties}
\label{sec: enumerative interpretation}
Let $\mathbf{s}$ be a skew-symmetric seed satisfying the $\mathcal{X}$ assumptions and 
let $(X,D)$ be a log Calabi-Yau compactification of the corresponding
$\mathcal{X}$ cluster variety as in \S \ref{sec_intro_A_X}.
We consider the algebra of theta functions $R(\foD_{\mathbf{s}}^{\mathcal{A}_{\mathrm{prin}}})$ obtained from the $\mathcal{A}_{\mathrm{prin}}$ cluster scattering diagram 
$\foD_{\mathbf{s}}^{\mathcal{A}_{\mathrm{prin}}}$. A basis of this algebra is given the cluster theta functions of the form 
$\vartheta_{(m,0)}^C$, as discussed in \S\ref{Sec: cluster theta functions}. By Lemma \ref{Lem:module}, they naturally span a 
$\mathbf{k}[\![N_{\mathrm{uf}}^{\oplus}]\!]$-algebra, and so for every 
$m_1, m_2, m$, we have structure constants 
$C_{m_1 m_2}^m \in \mathbf{k}[\![N_{\mathrm{uf}}^{\oplus}]\!]$
such that 
\[ \vartheta_{(m_1,0)}^C \vartheta_{(m_2,0)}^C = \sum_{m \in B(\Z)} 
C_{m_1 m_2}^m \vartheta_{(m,0)}^C \,.\]
The structure constants are power series 
\[ C_{m_1 m_2}^m=\sum_{\mathbf{A}=(a_i)_{i \in I} \in \NN^I} C_{m_1m_2, \mathbf{A}}^m \prod_{i \in I} z^{(0,a_i e_i)}\,,\]
with coefficients $C_{m_1 m_2,\mathbf{A}}^m \in \mathbf{k}$.
Actually, it follows from \cite[Theorem 1.13]{GHKK}
that 
\[ C_{m_1 m_2,\mathbf{A}}^m \in \NN\]
for every $m_1,m_2,m\in M$ and $\mathbf{A} \in \NN^I$.

On the other hand, for every $m_1, m_2, m\in M$ and 
$\beta \in \NE(X_\Sigma)$, Gross and Siebert define in \cite{gross2019intrinsic} counts of curves $N_{m_1 m_2, \beta}^m \in \Q$
in $(X,D)$. Let $\sigma_1$, $\sigma_2$, $\sigma$ be the smallest cones of $\Sigma$ containing respectively $m_1$, $m_2$, $m$, and let $D_{\sigma_1}$, $D_{\sigma_2}$, $D_{\sigma_3}$ be the corresponding strata of $(X,D)$.
Then, roughly, $N_{m_1 m_2, \beta}^m$ is a count of $3$-pointed degree $\beta$ rational curves \[f: (C,x_1,x_2,x_3) \longrightarrow X\] 
such that 
\[f(x_1) \in D_{\sigma_1}\,,\,\, 
f(x_2)\in D_{\sigma_2}\,,\,\, f(x_3) \in D_{\sigma_3}\,,\] 
with prescribed tangency conditions determined by 
$m_1,m_2,-m$, and such that $f(x_3)$ coincides with a fixed given point on $D_{\sigma}$. The precise definition of 
$N_{m_1 m_2, \beta}^m$, given in \cite[Definition 3.21]{gross2019intrinsic} uses logarithmic geometry and 
$N_{m_1 m_2, \beta}^m$ is in general an example of 
punctured log Gromov--Witten invariant\footnote{If $m \neq 0$, then $-m \notin \sigma$, and so the tangency condition $-m$ involves negative contact orders.}.

\begin{theorem}
\label{Thm: structure constants in cluster}
Let $\mathbf{s}$ be a skew-symmetric seed satisfying the $\mathcal{X}$ assumptions and 
let $(X,D)$ be a log Calabi-Yau compactification of the corresponding
$\mathcal{X}$ cluster variety as in \S \ref{sec_intro_A_X}. Then, the structure constants of the algebra of theta functions defined by the cluster scattering diagram $\foD_{\mathbf{s}}^{\mathcal{A}_{\mathrm{prin}}}$
are expressed in terms of the punctured log Gromov--Witten invariants of $(X,D)$ as follows:
for every $m_1, m_2, m \in M$ and $\mathbf{A}=(a_i)_{i\in I}\in \NN^I$, we have 
\[ C_{m_1 m_2,\mathbf{A}}^m  = \sum_{\overline{\beta}\in NE(X_\Sigma)}
N_{m_1 m_2, ( \overline{\beta}-\sum_{i\in I}a_i E_i)}^m \,.\]
\end{theorem}

\begin{proof}
Combining Theorem \ref{thm_main} and Theorem \ref{thm_restriction}, the structure constants computed from the 
$\mathcal{A}_{\mathrm{prin}}$ scattering diagram are obtained by setting to zero the curves classes coming from 
$NE(X_\Sigma)$ with respect to the decomposition 
$NE(X_\Sigma) \oplus \NN^I$ in the structure constants computed from the canonical scattering diagram $\foD_{(X,D)}$.
On the other hand, by \cite[Theorem C]{gross2021canonical}
the structure constants computed from the canonical scattering diagram are the punctured log Gromov-Witten invariants $N_{m_1 m_2, \beta}^m$ introduced in \cite{gross2019intrinsic}.
\end{proof}

\begin{remark}
Using \cite[Theorem 4.5]{gross2021canonical}, one can similarly give an enumerative interpretation of 
the coefficients of the monomial expansions $\vartheta_{(m,0)}^C(p)$  
of the cluster theta functions as in \eqref{Eq: theta function defined by broken line}
in terms of the counts of log broken lines defined in \cite[Definition 3.21]{gross2021canonical}.
\end{remark}

\begin{remark} \label{rem_enum_A}
One can also exchange the roles of $\mathcal{A}$ and $\mathcal{X}$: given a skew-symmetric seed $\mathbf{s}$ and $(X',D')$ a log Calabi-Yau compactification of the corresponding
$\mathcal{A}$ cluster variety as in \S \ref{sec_intro_A_X}, we obtain an enumerative interpretation of the structure constants of the algebra of theta functions defined by the cluster scattering diagram $\foD_{\mathbf{s}}^{\mathcal{X}_{\mathrm{prin}}}$
in terms of the punctured log Gromov--Witten invariants of $(X',D')$.
\end{remark}

 \FloatBarrier

\bibliographystyle{plain}
\bibliography{bibliography}

\end{document}

%% file: figcan.pspdftex
\begin{picture}(0,0)%
\includegraphics{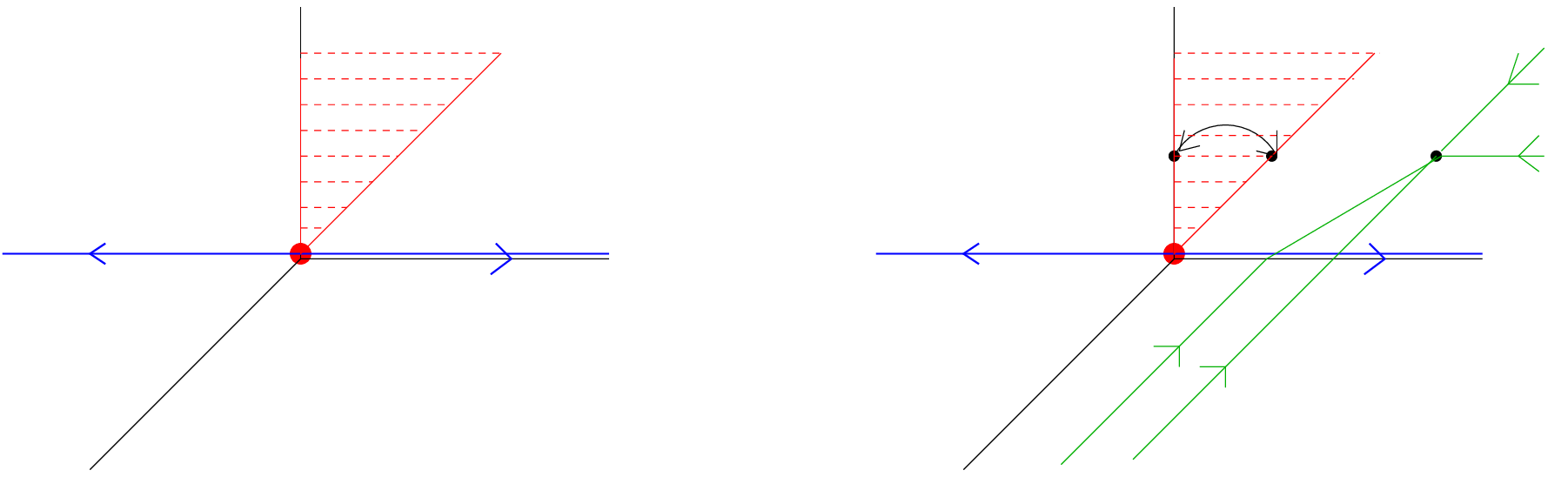}%
\end{picture}%
\setlength{\unitlength}{4144sp}%
\begingroup\makeatletter\ifx\SetFigFont\undefined%
\gdef\SetFigFont#1#2#3#4#5{%
  \reset@font\fontsize{#1}{#2pt}%
  \fontfamily{#3}\fontseries{#4}\fontshape{#5}%
  \selectfont}%
\fi\endgroup%
\begin{picture}(13717,4171)(-4926,-4220)
\put(5491,-2671){\makebox(0,0)[lb]{\smash{{\SetFigFont{12}{14.4}{\rmdefault}{\mddefault}{\updefault}{\color[rgb]{0,0,0}$0$}%
}}}}
\put(4861,-1456){\makebox(0,0)[lb]{\smash{{\SetFigFont{14}{16.8}{\rmdefault}{\mddefault}{\updefault}{\color[rgb]{0,0,0}$(0,1)$}%
}}}}
\put(6391,-1456){\makebox(0,0)[lb]{\smash{{\SetFigFont{14}{16.8}{\rmdefault}{\mddefault}{\updefault}{\color[rgb]{0,0,0}$(1,1)$}%
}}}}
\put(-2159,-2671){\makebox(0,0)[lb]{\smash{{\SetFigFont{12}{14.4}{\rmdefault}{\mddefault}{\updefault}{\color[rgb]{0,0,0}$0$}%
}}}}
\put(-4904,-466){\makebox(0,0)[lb]{\smash{{\SetFigFont{14}{16.8}{\rmdefault}{\mddefault}{\updefault}{\color[rgb]{0,0,0}$\foD_{(X,D)}$}%
}}}}
\put(-4904,-2086){\makebox(0,0)[lb]{\smash{{\SetFigFont{14}{16.8}{\rmdefault}{\mddefault}{\updefault}{\color[rgb]{0,0,1}$f_{\fod'}=1+xt^{L-E}$}%
}}}}
\put(8776,-466){\makebox(0,0)[lb]{\smash{{\SetFigFont{14}{16.8}{\rmdefault}{\mddefault}{\updefault}{\color[rgb]{0,.69,0}$\vartheta_{0,1}=xy$}%
}}}}
\put(8731,-1411){\makebox(0,0)[lb]{\smash{{\SetFigFont{14}{16.8}{\rmdefault}{\mddefault}{\updefault}{\color[rgb]{0,.69,0}$\vartheta_{1,0}=x$}%
}}}}
\put(5356,-4111){\makebox(0,0)[lb]{\smash{{\SetFigFont{14}{16.8}{\rmdefault}{\mddefault}{\updefault}{\color[rgb]{0,.69,0}$\vartheta_{1,0}=x^{-1}y^{-1}f_{\fod}t^{L-E}$}%
}}}}
\put(-2204,-196){\makebox(0,0)[lb]{\smash{{\SetFigFont{14}{16.8}{\rmdefault}{\mddefault}{\updefault}{\color[rgb]{0,0,0}Kink of $\varphi$: $L$}%
}}}}
\put(-359,-1996){\makebox(0,0)[lb]{\smash{{\SetFigFont{14}{16.8}{\rmdefault}{\mddefault}{\updefault}{\color[rgb]{0,0,1}$f_{\fod}=1+x^{-1}t^{E}$}%
}}}}
\put(7741,-2086){\makebox(0,0)[lb]{\smash{{\SetFigFont{14}{16.8}{\rmdefault}{\mddefault}{\updefault}{\color[rgb]{0,0,1}$f_{\fod}=1+x^{-1}t^{E}$}%
}}}}
\put(-3869,-4156){\makebox(0,0)[lb]{\smash{{\SetFigFont{14}{16.8}{\rmdefault}{\mddefault}{\updefault}{\color[rgb]{0,0,0}Kink of $\varphi$: $L$}%
}}}}
\put(7426,-2806){\makebox(0,0)[lb]{\smash{{\SetFigFont{14}{16.8}{\rmdefault}{\mddefault}{\updefault}{\color[rgb]{0,0,0}Kink of $\varphi$: $L-E$}%
}}}}
\put(-44,-2626){\makebox(0,0)[lb]{\smash{{\SetFigFont{14}{16.8}{\rmdefault}{\mddefault}{\updefault}{\color[rgb]{0,0,0}Kink of $\varphi$: $L-E$}%
}}}}
\put(-1079,-2671){\makebox(0,0)[lb]{\smash{{\SetFigFont{14}{16.8}{\rmdefault}{\mddefault}{\updefault}{\color[rgb]{0,0,.82}$\fod$}%
}}}}
\put(-4184,-2671){\makebox(0,0)[lb]{\smash{{\SetFigFont{14}{16.8}{\rmdefault}{\mddefault}{\updefault}{\color[rgb]{0,0,.82}$\fod'$}%
}}}}
\end{picture}%

%% file: cones.pspdftex
\begin{picture}(0,0)%
\includegraphics{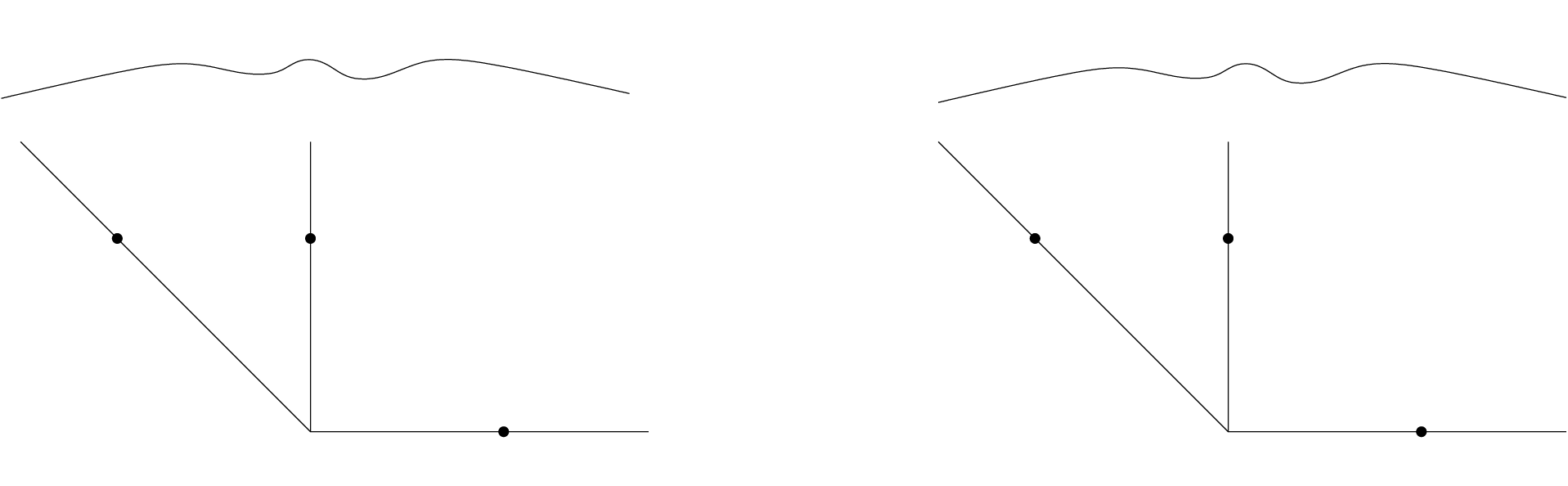}%
\end{picture}%
\setlength{\unitlength}{4144sp}%
\begingroup\makeatletter\ifx\SetFigFont\undefined%
\gdef\SetFigFont#1#2#3#4#5{%
  \reset@font\fontsize{#1}{#2pt}%
  \fontfamily{#3}\fontseries{#4}\fontshape{#5}%
  \selectfont}%
\fi\endgroup%
\begin{picture}(14604,4522)(2059,-4166)
\put(2611,-2311){\makebox(0,0)[lb]{\smash{{\SetFigFont{14}{16.8}{\rmdefault}{\mddefault}{\updefault}{\color[rgb]{0,0,0}$L-E$}%
}}}}
\put(5176,-1861){\makebox(0,0)[lb]{\smash{{\SetFigFont{14}{16.8}{\rmdefault}{\mddefault}{\updefault}{\color[rgb]{0,0,0}$L$}%
}}}}
\put(6841,-4111){\makebox(0,0)[lb]{\smash{{\SetFigFont{14}{16.8}{\rmdefault}{\mddefault}{\updefault}{\color[rgb]{0,0,0}$E$}%
}}}}
\put(11296,-2311){\makebox(0,0)[lb]{\smash{{\SetFigFont{14}{16.8}{\rmdefault}{\mddefault}{\updefault}{\color[rgb]{0,0,0}$L-E$}%
}}}}
\put(13681,-1861){\makebox(0,0)[lb]{\smash{{\SetFigFont{14}{16.8}{\rmdefault}{\mddefault}{\updefault}{\color[rgb]{0,0,0}$L$}%
}}}}
\put(15391,-4111){\makebox(0,0)[lb]{\smash{{\SetFigFont{14}{16.8}{\rmdefault}{\mddefault}{\updefault}{\color[rgb]{0,0,0}$E$}%
}}}}
\put(3511,-1411){\makebox(0,0)[lb]{\smash{{\SetFigFont{25}{30.0}{\rmdefault}{\mddefault}{\updefault}{\color[rgb]{0,0,0}$\mathcal{M}$}%
}}}}
\put(13501,209){\makebox(0,0)[lb]{\smash{{\SetFigFont{25}{30.0}{\rmdefault}{\mddefault}{\updefault}{\color[rgb]{0,0,0}$\mathcal{M}^{\vee}$}%
}}}}
\put(4726,119){\makebox(0,0)[lb]{\smash{{\SetFigFont{25}{30.0}{\rmdefault}{\mddefault}{\updefault}{\color[rgb]{0,0,0}$NE(X)$}%
}}}}
\put(15031,-1411){\makebox(0,0)[lb]{\smash{{\SetFigFont{25}{30.0}{\rmdefault}{\mddefault}{\updefault}{\color[rgb]{0,0,0}$\mathrm{Bog}(X)$}%
}}}}
\put(11836,-1411){\makebox(0,0)[lb]{\smash{{\SetFigFont{25}{30.0}{\rmdefault}{\mddefault}{\updefault}{\color[rgb]{0,0,0}$\mathrm{Nef}(X)$}%
}}}}
\end{picture}%

%% file: p.pspdftex
\begin{picture}(0,0)%
\includegraphics{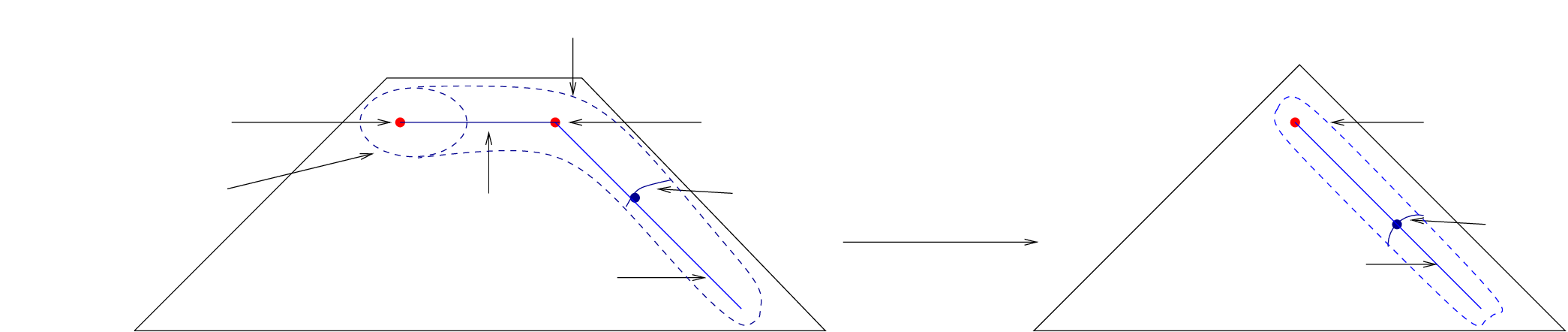}%
\end{picture}%
\setlength{\unitlength}{4144sp}%
\begingroup\makeatletter\ifx\SetFigFont\undefined%
\gdef\SetFigFont#1#2#3#4#5{%
  \reset@font\fontsize{#1}{#2pt}%
  \fontfamily{#3}\fontseries{#4}\fontshape{#5}%
  \selectfont}%
\fi\endgroup%
\begin{picture}(15912,3354)(-1814,-7723)
\put(271,-7486){\makebox(0,0)[lb]{\smash{{\SetFigFont{20}{24.0}{\rmdefault}{\mddefault}{\updefault}{\color[rgb]{0,0,0}$T$}%
}}}}
\put(9226,-7576){\makebox(0,0)[lb]{\smash{{\SetFigFont{20}{24.0}{\rmdefault}{\mddefault}{\updefault}{\color[rgb]{0,0,0}$\Spec \mathbf{k}[\mathcal{M}]$}%
}}}}
\put(-539,-5596){\makebox(0,0)[lb]{\smash{{\SetFigFont{20}{24.0}{\rmdefault}{\mddefault}{\updefault}{\color[rgb]{0,0,0}$t^{E_i}=0$}%
}}}}
\put(2971,-6631){\makebox(0,0)[lb]{\smash{{\SetFigFont{20}{24.0}{\rmdefault}{\mddefault}{\updefault}{\color[rgb]{0,0,0}$T_0$}%
}}}}
\put(3916,-4516){\makebox(0,0)[lb]{\smash{{\SetFigFont{20}{24.0}{\rmdefault}{\mddefault}{\updefault}{\color[rgb]{0,0,0}$\hat{T}$}%
}}}}
\put(-1799,-6496){\makebox(0,0)[lb]{\smash{{\SetFigFont{20}{24.0}{\rmdefault}{\mddefault}{\updefault}{\color[rgb]{0,0,0}$\mathrm{Spf} \,\mathbf{k}\lfor NE(X) \rfor$}%
}}}}
\put(5851,-6451){\makebox(0,0)[lb]{\smash{{\SetFigFont{20}{24.0}{\rmdefault}{\mddefault}{\updefault}{\color[rgb]{0,0,0}$\mathrm{Spf}\, \mathbf{k}\lfor \NN^I \rfor$}%
}}}}
\put(5491,-5596){\makebox(0,0)[lb]{\smash{{\SetFigFont{20}{24.0}{\rmdefault}{\mddefault}{\updefault}{\color[rgb]{0,0,0}$t^{-E_i}=0$}%
}}}}
\put(2116,-7261){\makebox(0,0)[lb]{\smash{{\SetFigFont{20}{24.0}{\rmdefault}{\mddefault}{\updefault}{\color[rgb]{0,0,0}$\Spec \mathbf{k}[ NE(X_{\Sigma}) ]$}%
}}}}
\put(4546,-6181){\makebox(0,0)[lb]{\smash{{\SetFigFont{10}{12.0}{\rmdefault}{\mddefault}{\updefault}{\color[rgb]{0,0,0}$1$}%
}}}}
\put(12331,-6496){\makebox(0,0)[lb]{\smash{{\SetFigFont{10}{12.0}{\rmdefault}{\mddefault}{\updefault}{\color[rgb]{0,0,0}$1$}%
}}}}
\put(13366,-6586){\makebox(0,0)[lb]{\smash{{\SetFigFont{20}{24.0}{\rmdefault}{\mddefault}{\updefault}{\color[rgb]{0,0,0}$\mathrm{Spf}\, \mathbf{k}\lfor \NN^I \rfor$}%
}}}}
\put(12826,-5596){\makebox(0,0)[lb]{\smash{{\SetFigFont{20}{24.0}{\rmdefault}{\mddefault}{\updefault}{\color[rgb]{0,0,0}$S_{(X_{\Sigma},H)}$}%
}}}}
\put(7471,-7081){\makebox(0,0)[lb]{\smash{{\SetFigFont{20}{24.0}{\rmdefault}{\mddefault}{\updefault}{\color[rgb]{0,0,0}$p$}%
}}}}
\put(9676,-7081){\makebox(0,0)[lb]{\smash{{\SetFigFont{20}{24.0}{\rmdefault}{\mddefault}{\updefault}{\color[rgb]{0,0,0}$\Spec \mathbf{k}[ NE(X_{\Sigma}) ]$}%
}}}}
\end{picture}%